\documentclass[11pt]{amsart}
\usepackage[margin=1in]{geometry}
\usepackage[T1]{fontenc}
\usepackage{lmodern,mathtools,amssymb,amsthm,microtype,longtable,booktabs,array,mathrsfs}
\usepackage[hidelinks]{hyperref}
\newtheorem{theorem}{Theorem}[section]
\newtheorem{proposition}[theorem]{Proposition}
\newtheorem{lemma}[theorem]{Lemma}
\newtheorem{corollary}[theorem]{Corollary}
\theoremstyle{definition}

\newtheorem{remark}[theorem]{Remark}
\newtheorem{conjecture}[theorem]{Conjecture}
\newtheorem*{ac}{Acknowledgements}

\newcommand{\bZ}{\mathbb Z}
\newcommand{\bN}{\mathbb N}
\newcommand{\bR}{\mathbb R}
\newcommand{\bC}{\mathbb C}
\newcommand{\bE}{\mathbb E}
\newcommand{\PP}{\mathscr P}
\newcommand{\tr}{\operatorname{tr}}

\newcommand{\U}{\mathcal U}
\newcommand{\cR}{\mathcal R}
\newcommand{\cN}{\mathcal N}
\newcommand{\cH}{\mathcal H}
\newcommand{\cM}{\mathcal M}
\newcommand{\cQ}{\mathcal Q}
\newcommand{\KNM}[2]{\mathbf K^{ #1 \subset #2}}
\newcommand{\KI}{\mathbf K^{ \mathscr I}}
\newcommand{\KIf}{\mathbf K^{ \mathscr I_f}}

\newcommand{\cb}{\mathrm cb}

\newcommand{\sgn}{\operatorname{sgn}}
\newcommand{\id}{\operatorname{id}}
\newcommand{\op}{\mathrm{op}}

\title{Bimodule Grothendieck Inequalities}

\author{Jincheng Wan}
\address{Jincheng Wan, Tsinghua University, Beijing}
\email{wanjc23@mails.tsinghua.edu.cn }

\author{Jinsong Wu}
\address{Jinsong Wu, Beijing Institute of Mathematical Sciences and Applications, Beijing, 101408, China}
\email{wjs@bimsa.cn}

\date{}

\begin{document}

	\begin{abstract}
		In this paper, we systematically investigate Grothendieck inequalities in the framework of quantum Fourier analysis.
		We obtain a sharper Grothendieck inequality for irreducible inclusions of finite factors.
		In particular, we obtain a sharp Grothendieck inequality for positive bimodule operators.
		For finite-index inclusions, we improve the upper bound by an explicit index-dependent bound.
		We also obtain multilinear Grothendieck inequalities depending on the Jones index, as well as index-free Grothendieck inequalities for irreducible inclusions of finite factors.
	\end{abstract}
	
	\maketitle
	% \tableofcontents
	% \clearpage
	
	\section{Introduction}

	The remarkable Grothendieck inequalities were initially studied by Grothendieck \cite{Gro53} in 1953.
	There are Grothendieck inequalities for the real and complex fields, and the optimal constants are denoted by \(K_{\mathrm{Gr}}^{\mathbb{R}}\) and \(K_{\mathrm{Gr}}^{\mathbb{C}}\), respectively.
	In 1977, Krivine \cite{Kri77} improved the upper bound for the real Grothendieck constant:
	\[
	K_{\mathrm{Gr}}^{\mathbb{R}}\leq \frac{\pi}{2\log(1+\sqrt{2})}.
	\]
	In 1987, Haagerup \cite{HaagerupC} improved the upper bound for the complex Grothendieck constant:
	\[
	K_{\mathrm{Gr}}^{\mathbb{C}}\leq \gamma_{H}^{-1}\simeq 1.40490913.
	\]
	There are fruitful results on sharpening the bounds of the Grothendieck constants.
	Recently, Heilman, Jones, and Malavolta \cite{HJM} improved the bound for the complex Grothendieck constant:
	\[
	1.4<K_{\mathrm{Gr}}^{\mathbb{C}}<1.404898554746.
	\]
	For further references on the Grothendieck constant, we refer to \cite{HJM} and the references therein.
	The Grothendieck inequalities have been extended to many other fields of mathematics, such as graph theory \cite{Alon2006}, \(L^p\) spaces \cite{Guruswami2012}, and \(C^*\)-algebras.
	
	The Grothendieck inequalities for \(C^*\)-algebras were studied by Pisier \cite{Pisier78} and Haagerup \cite{Haagerup}.
	The Grothendieck inequalities for operator spaces were studied by Pisier--Shlyakhtenko \cite{PS}, Xu \cite{Xu06}, among others.
	We refer to \cite{Pis12} for the development of the noncommutative Grothendieck inequalities.
	
	The modern subfactor theory initiated by Jones \cite{Jones1983} is an important subject in operator algebras, low-dimensional topology, category theory, and other areas.
	Jiang, Liu, and the second author \cite{JLW16, JLW19} studied a series of Fourier inequalities in subfactor theory.
	Liu, Palcoux, and the second author \cite{LPW} showed that these inequalities are useful in obtaining effective criteria for obstructions to unitary fusion categorification.
	All of this work initiated the program of quantum Fourier analysis \cite{JJL20} by Jaffe, Jiang, Liu, Ren, and the second author.
	
	In this paper, we shall systematically study the Grothendieck inequalities for irreducible inclusions of finite factors.
	We improve the upper bound for the Grothendieck constant for irreducible inclusions to \(\pi/2\), compared with the general bound \(2\), by using Popa's free independence theorem \cite{PopaFree95}.
	We obtain a sharper Grothendieck inequality for finite-index inclusions with fixed index.
	We prove that the sharp Grothendieck constant for bimodule maps that are positive operators is $\frac{9\pi^2}{64}.$
	We show that the Grothendieck constant for group--subgroup inclusions is one, while the Grothendieck constant for their dual inclusions is bounded by the classical complex Grothendieck constant.
	We also show that the lower bound of the Grothendieck constant for irreducible inclusions is larger than the complex Grothendieck constant.
	Hence the Grothendieck constant can be considered an interesting invariant for subfactors and fusion categories.
	We also obtain multilinear Grothendieck inequalities for irreducible finite-index inclusions of finite factors, where the constants depend on the Jones index.
	Inspired by the work of Christensen--Pop--Sinclair--Smith \cite{CPSS}, we obtain Jones-index-free Grothendieck inequalities without requiring hyperfiniteness of the subfactors.
	
	The paper is organized as follows.
	In Section 2, we recall fundamentals from subfactor theory, bimodule maps on subfactors, and basics in free probability.
	In Section 3, we improve Haagerup's Grothendieck inequality for irreducible inclusions of finite factors and obtain a sharp Grothendieck inequality for operator-positive bimodule maps.
	We also show that the Grothendieck constant can be attained at an irreducible inclusion.
	In Section 4, we characterize the Grothendieck constant for group--subgroup inclusions and their dual inclusions, subfactor intervals, tensor products, and free compositions.
	In Section 5, we obtain bilinear Grothendieck inequalities for irreducible inclusions of finite factors.
	In Section 6, we obtain multilinear Grothendieck inequalities for irreducible inclusions of finite factors.
	In Section 7, we obtain index-free Grothendieck inequalities for irreducible inclusions of finite factors.
	
	\begin{ac}
		J. Wu was supported by grants from Beijing Institute of Mathematical Sciences and Applications.
		J. Wu was supported by NSFC (Grant no. 12371124). 
	\end{ac}
	
	\section{Preliminaries}
	
	\subsection{Subfactors}
	
	Let $\cN\subset\cM$ be a unital inclusion of finite factors, with normalized trace $\tau$, trace-preserving conditional expectation $\bE_{\cN}$, and Jones index $\delta^2=[\cM:\cN]$ when the index is finite. 
	The Hilbert space $L^2(\cM, \tau)$ is the Gelfand-Naimark-Segal Hilbert space of $\tau$, with cyclic separating vector $\Omega$ and modular conjugation $J$ given by $Jx\Omega = x^*\Omega$, \ $x\in\cM$. 
	We denote by $e_{\cN}$ the Jones projection.  
	The basic construction $\cM_1=\langle \cM, e_{\cN}\rangle$ is the von Neumann algebra generated by $\cM$ and $e_{\cN}$. 
	The inclusion $\cN\subset\cM$ is called finite if $\cM_1$ is a finite factor, and irreducible if $\cN'\cap \cM=\mathbb{C}$.  
	For a finite inclusion, we have $\cM_1 = J\cN'J$, where $\cN'$ is the commutant of $\cN$ on $L^2(\cM,\tau)$. 
	For a finite inclusion, there exists a finite set $\{\eta_j\}_{j=1}^m$ of operators in $\cM$ called Pimsner-Popa basis for $\cN\subset \cM$, which satisfies $\displaystyle x=\sum_{j=1}^m  \bE_{\mathcal{N}}(x\eta^*_j)\eta_j$, for all $x\in\cM$\cite{PPindex}.

	We write $\cR(a)$ for the range projection of an operator $a$. 
	At finite index, if $p\in\cM$ is a projection and $q=\cR(\bE_{\cN}(p))$, then
	\begin{equation}\label{eq:support}
		p\le q,\qquad \tau(q)\le\delta^2\tau(p).
	\end{equation}
	Indeed, $\|(1-q)p\|_2^2=\tau((1-q)\bE_{\cN}(p))=0$, so $p\le q$. 
	Let $\tau_1$ be the normalized trace in the first basic construction. 
	Since $e_{\cN}pe_{\cN}=\bE_{\cN}(p)e_{\cN}$, the left support of $e_{\cN}p$ is $qe_{\cN}$ and its right support is dominated by $p$. 
	Polar decomposition in the basic construction and the Markov trace identity give
	\[
	\delta^{-2}\tau(q)=\tau_1(qe_{\cN})\le \tau_1(p)=\tau(p).
	\]
	We also use the Pimsner--Popa inequality $x\le\delta^2\bE_{\cN}(x)$ for $x\ge0$ in $\cM$ \cite{PPindex}.
	Throughout this paper, we denote by $C_n$ the finite cyclic group of order $n$.
	
	\subsection{Norms of bimodule maps}
	
	For the linear inequalities, a bimodule map $\Phi$ is a bounded operator on $L^2(\cM,\tau)$ commuting with the left and right actions of $\cN$. We use the inner product $\langle y,x\rangle=\tau(y^*x)$, linear in the second slot. For $1\le p,q\le\infty$ write
	\[
	\|\Phi\|_{p\to q}=\sup_{0\ne x\in\cM}\frac{\|\Phi(x)\|_q}{\|x\|_p},
	\]
	with value $+\infty$ when the indicated bound fails, including when an image does not belong to $L^q(\cM)$. In particular, $\|\Phi\|_{2\to2}$ is the operator norm on $L^2(\cM)$. When $\Phi$ maps $\cM$ boundedly to itself, its usual completely bounded norm is denoted by $\|\Phi\|_{\cb}$.
	The multilinear maps are specified separately below.
	
	\begin{lemma}\label{lem:norm}
		Let $\cN\subset\cM$ be an inclusion of finite von Neumann algebras with a faithful normal tracial state, and let $\Phi$ be a bounded bimodule operator on $L^2(\cM)$.
		Then
		\begin{equation}\label{eq:1norm}
			\|\Phi\|_{\infty\to1}=\sup_{u,v\in\U(\cM)}|\tau(v^*\Phi(u))|.
		\end{equation}
		If the inclusion is irreducible, then
		\begin{equation}\label{eq:22norm}
			\|\Phi\|_{2\to2}=
			\sup\left\{
			\left|\sum_j\tau(y_j^*\Phi(x_j))\right|:
			\sum_j x_j^*x_j,\ \sum_j x_jx_j^*,\ \sum_j y_j^*y_j,\ \sum_j y_jy_j^*\le1
			\right\},
		\end{equation}
		where the supremum runs over finite lists in $\cM$.
	\end{lemma}
	
	\begin{proof}
		Trace duality gives $\|\Phi\|_{\infty\to1}$ as the supremum of $|\tau(y^*\Phi(x))|$ over contractions $x,y$. The Russo--Dye theorem \cite{RussoDye}, applied in each slot, gives \eqref{eq:1norm}.
		
		For \eqref{eq:22norm}, the Cauchy--Schwarz inequality in the direct sum of copies of $L^2(\cM)$ gives the upper bound: 
		each of the two sums of squared $L^2$ norms is at most one. 
		For the reverse bound, fix $x,y\in\cM$ with $\|x\|_2=\|y\|_2=1$ and $\displaystyle |\langle y,\Phi(x)\rangle|>\|\Phi\|_{2\to2}-\varepsilon$. 
		Under $\U(\cN)$ conjugation the invariant subspace of $L^2(\cM)$ is $\bC1$: 
		the bounded truncations of any invariant affiliated operator belong to $\cN'\cap\cM$. 
		The projection onto the invariant subspace belongs to the closed convex hull of each orbit, by the unique vector of minimum norm in that hull. 
		Applying this argument to the tuple $(xx^*,x^*x,yy^*,y^*y)$ gives unitaries $u_j\in\cN$ and weights $\beta_j\ge0$ summing to one such that the averages of all four entries are within $\varepsilon$ of $1$ in $L^2(\cM)$. 
		Put
		\[
		x_j=\sqrt{\beta_j}u_jxu_j^*,\qquad
		y_j=\sqrt{\beta_j}u_jyu_j^*.
		\]
		Bimodularity and traciality give $\displaystyle \sum_j\langle y_j,\Phi(x_j)\rangle=\langle y,\Phi(x)\rangle$.
		
		For either list $(a_j)_j$ let $\displaystyle \mathbf a_{+}=\sum_j a_ja_j^*$, $\displaystyle \mathbf a_{-}=\sum_j a_j^*a_j$.
		Suppose that $1\vee t=\max(1,t)$. 
		The corrected list $\widetilde a_j=(1\vee \mathbf a_{+})^{-1/2}a_j(1\vee \mathbf a_{-})^{-1/2}$ satisfies
		\[
		\sum_j\widetilde a_j\widetilde a_j^*\le (1\vee \mathbf a_{+})^{-1/2} \mathbf a_{+} (1\vee \mathbf a_{+})^{-1/2} \le1,
		\qquad
		\sum_j\widetilde a_j^*\widetilde a_j\le (1\vee \mathbf a_{-})^{-1/2} \mathbf a_{-} (1\vee \mathbf a_{-})^{-1/2}\le1.
		\]
		The identities
		\begin{align*}
			\sum_j\|(1-(1\vee \mathbf a_{+})^{-1/2})a_j\|_2^2
			=& \tau\left((1-(1\vee \mathbf a_{+})^{-1/2})^2\mathbf a_{+} \right)=\| (\sqrt \mathbf a_{+}-1)_+\|_2^2, \\
			\sum_j\|a_j(1-(1\vee \mathbf a_{-})^{-1/2})\|_2^2=& \| (\sqrt \mathbf a_{-}-1)_+\|_2^2
		\end{align*}
		show, by writing 
		\begin{align*}
			(1\vee \mathbf a_{+})^{-1/2} a_j(1\vee \mathbf a_{-})^{-1/2}-a_j=((1\vee \mathbf a_{+})^{-1/2}-1)a_j (1\vee \mathbf a_{-})^{-1/2}+a_j((1\vee \mathbf a_{-})^{-1/2}-1),
		\end{align*}
		that the correction has direct sum $L^2$ norm at most $\|\mathbf a_{+} -1\|_2+\|\mathbf a_{-}-1\|_2$. 
		Thus both lists change by at most $2\varepsilon$ in direct sum norm. 
		Boundedness of $\Phi$ on $L^2(\cM)$ and Cauchy--Schwarz inequality imply that the pairing changes by at most $4\varepsilon\|\Phi\|_{2\to2}$. 
		Letting $\varepsilon$ tend to zero proves the formula.
	\end{proof}
	
	We write $\Phi\ge0$ for positivity of $\Phi$ as an operator on $L^2(\cM)$. 
	A bimodule operator is bi-Hermitian if it is self-adjoint on $L^2(\cM)$ and $\Phi(x^*)=\Phi(x)^*$, where the involution extends to $L^2(\cM)$. 
	We denote the space of such operators by $\mathbb H_{sa}$.
	
	\begin{remark}
		The Christensen--Sinclair matrix amplification of a bilinear form $u:\mathcal A\times\mathcal B\to\bC$ gives the norm \cite{CS}
		\[
		\|u\|_{\mathrm{CS}}=
		\sup\left\{\left|\sum_j u(a_j,b_j)\right|:
		\sum_j a_ja_j^*\le1,\quad \sum_j b_j^*b_j\le1\right\}.
		\]
		For $u_{\Phi}(a,b)=\tau(a\Phi(b))$, Cauchy--Schwarz inequality gives
		$\|u_\Phi\|_{\mathrm{CS}}\le\|\Phi\|_{2\to2}$ using
		$\displaystyle \sum_j\|a_j\|_2^2$, $\displaystyle \sum_j\|b_j\|_2^2\le1$. 
		When the inclusion is irreducible, taking $a_j=y_j^*$ and $b_j=x_j$ in \eqref{eq:22norm} gives the reverse bound. This norm pertains to the bilinear form $u_\Phi$.
	\end{remark}
	
	\begin{lemma}\label{lem:map-cb-norm}
		Suppose that $\cN\subset\cM$ is an irreducible inclusion of finite factors of finite index and $\Phi:\cM\to\cM$ is a bounded bimodule map. Then $\|\Phi\|_{\infty\to\infty}=\|\Phi\|_{\cb}$.
	\end{lemma}
	
	\begin{proof}
		If $\cM$ is finite dimensional, then irreducibility forces $\cN=\cM$ and $\Phi$ is scalar. 
		For $\mathrm{II}_1$ factors, a finite-index subfactor $\cN$ norms $\cM$ \cite[Corollary~3.3]{PSS}. 
		Explicitly, the norm of $X\in M_n(\cM)$ is the supremum of $\|RXC\|$ over contractive rows $R\in M_{1,n}(\cN)$ and columns $C\in M_{n,1}(\cN)$. 
		Since $R\Phi^{(n)}(X)C=\Phi(RXC)$, this identity gives
		\[
		\|\Phi^{(n)}(X)\|\le\|\Phi\|_{\infty\to\infty}\|X\|.
		\]
		Take the supremum over $n$. The level $n=1$ gives the reverse inequality.
	\end{proof}
	
	The Grothendieck constant for an inclusion $\cN \subset \cM$ is defined as follows
	\[
	\KNM{\cN}{\cM} =\sup_{\Phi\ne0}\frac{\|\Phi\|_{2\to 2}}{\|\Phi\|_{\infty\to 1}}, \quad
	\KNM{\cN}{\cM}_+ =\sup_{0\leq \Phi\ne0}\frac{\|\Phi\|_{2\to 2}}{\|\Phi\|_{\infty\to 1}}.
	\]
	Let $\mathscr I$ be the class of all irreducible inclusions of finite factors and $\mathscr I_f$ be the class of all irreducible finite-index inclusions of finite factors.
	The Grothendieck constant for $\mathscr I$ is defined as
	\begin{align*}
		\KI=\sup_{\mathscr I} \KNM{\cN}{\cM}, \quad
		\KI_+=\sup_{\mathscr I} \KNM{\cN}{\cM}_+.
	\end{align*}
	The Grothendieck constant for $\mathscr I_f$ is defined as
	\begin{align*}
		\KIf=\sup_{\mathscr I_f} \KNM{\cN}{\cM}, \quad
		\KIf_+=\sup_{\mathscr I_f} \KNM{\cN}{\cM}_+.
	\end{align*}
	The completely bounded Grothendieck constant for an inclusion $\cN \subset \cM$ is defined as follows
	\[
	\KNM{\cN}{\cM}_{\cb} =\sup_{\Phi\ne0}\frac{\|\Phi\|_{\cb}}{\|\Phi\|_{\infty\to 1}},
	\]
	where for this supremum $\Phi$ is required to map $\cM$ boundedly into itself. 
	The completely bounded Grothendieck constant for $\mathscr I$ is defined by
	\begin{align*}
		\KI_{\cb}=\sup_{\mathscr I} \KNM{\cN}{\cM}_{\cb}.
	\end{align*}

	\subsection{Bimodule Multilinear maps}
	
	Suppose that $\Phi:\cM\times \cM\to \cM$ be bounded, separately normal, and $\cN$-multimodular in the standard sense:
	\begin{align}\label{eq:balancing}
		\Phi(a x,y)=a \Phi(x,y),\quad \Phi(xa,y)=\Phi(x,ay),\quad
		\Phi(x,ya)=\Phi(x,y)a,
	\end{align}
	where $x, y\in \cM$ and $a\in \cN$.
	We also call the $\cN$-multimodularity as the $\cN$-balancing condition or balance for simplicity.
	Define
	\[
	\|\Phi\|_{\infty, \infty \to 1}=\sup_{\|x\|_\infty,\|y\|_\infty\le1}\|\Phi(x,y)\|_1,
	\qquad
	\|\Phi\|_{2, 2\to 1}=\sup_{\|x\|_2,\|y\|_2\le1}\|\Phi(x,y)\|_1,
	\]
	where $x,y\in\cM$ and the value may be infinite.
	Trace duality and the Russo--Dye theorem in each slot give $\displaystyle \|\Phi\|_{\infty, \infty \to 1}=\sup_{u, v, w\in \mathcal U (\cM)}|\tau(w^*\Phi(u,v))|$.
	We shall denote $\|\Phi\|_{\infty, \infty \to 1}$ by $\|\Phi\|_{\infty^2 \to 1}$ and $\|\Phi\|_{2, 2\to 1}$ by $\|\Phi\|_{2^2\to 1}$ for simplicity.
	The Grothendieck constant for these bimodule maps is defined as
	\begin{align*}
		\KNM{\cN}{\cM}_2=\sup_{\Phi\neq 0}\frac{\|\Phi\|_{2^2 \to 1}}{\|\Phi\|_{\infty^2 \to 1}}
	\end{align*}
	
	For $n\ge2$ let $\Phi:\cM^{\times n}\to \cM$ be bounded, separately normal and $\cN$-multimodular.
	For $1\le p_i\le\infty$ put
	\[
	\|\Phi\|_{\boldsymbol p\to1}=\|\Phi\|_{p_1,\ldots,p_n\to1}
	=\sup_{\substack{x_i\in\cM\\ \|x_i\|_{p_i}\le1}}\|\Phi(x_1,\ldots,x_n)\|_1.
	\]
	We use $\|\Phi\|_{\infty^n\to1}$ when all input exponents are $\infty$ and write
	$\gamma_{\boldsymbol p}(\Phi)=\|\Phi\|_{\boldsymbol p\to1}$ for a map denoted by $\Phi$.

	\section{Grothendieck Inequalities}
	
	Suppose that $u:\mathcal A\times \mathcal B\to\mathbb C$ is a bounded complex bilinear form  on $C^*$-algebras $\mathcal{A}$, $\mathcal {B}$, with norm $\displaystyle \|u\|=\sup_{\|a\|,\|b\|\le1}|u(a,b)|$.
	The noncommutative Grothendieck theorem, obtained by Pisier initially and by Haagerup \cite{Pisier78,Haagerup}, supplies states $f_1,f_2$ on $\mathcal A$ and $g_1,g_2$ on $\mathcal B$ such that
	\begin{equation}\label{eq:haagerup}
		|u(a,b)|\le\|u\|
		\left(f_1(a^*a)+f_2(aa^*)\right)^{1/2} \left( g_1(b^*b)+g_2(bb^*)\right)^{1/2}.
	\end{equation}
	Consequently,
	\begin{equation}\label{eq:ordinary-list}
		\left|\sum_i u(a_i,b_i)\right|\le2\|u\|
		\max\left\{\left\|\sum_i a_i^*a_i\right\|^{1/2},
		\left\|\sum_i a_i a_i^*\right\|^{1/2}\right\}
		\max\left\{\left\|\sum_i b_i^*b_i\right\|^{1/2},
		\left\|\sum_i b_i b_i^*\right\|^{1/2}\right\}.
	\end{equation}
	The coefficient $2$ comes from summing two states of norm one in each square root.
	
	\subsection{Bimodule Grothendieck inequality}
	For irreducible inclusions we improve the factor two using a free circular limit.
	For the basics of the free probability theory initiated by Voiculescu, we refer to \cite{NS} and the references therein.
	The centered semicircle law of variance one has density $\displaystyle \frac{1}{2\pi}\sqrt{4-s^2}$ on $[-2,2]$. 
	Its absolute value has the quarter-circle law
	\[
	d\nu(s)=\pi^{-1}\sqrt{4-s^2}\,ds\qquad (0\le s\le2).
	\]
	The pushforward by $s\mapsto s/2$ has density $\displaystyle \frac{4}{\pi}\sqrt{1-r^2}$ for $0\le r\le1$.
	
	\begin{theorem}\label{gb:theorem}
		Suppose that $\cN\subset \cM$ is an irreducible inclusion of finite factors.
		Then
		\begin{align*}
			1\leq  \KNM{\cN}{\cM}\leq \frac{\pi}{2}.
		\end{align*}
		If $\Phi\ge0$ is a bimodule map as an operator on $L^2(\cM)$, then
		\[
		\|\Phi\|_{2\to 2}\le\frac{9\pi^2}{64}\|\Phi\|_{\infty \to 1}.
		\]
	\end{theorem}

	\begin{proof}
		The lower bound follows from $\|\Phi\|_{\infty\to1}\le\|\Phi\|_{2\to2}$. 
		For comparison, applying \eqref{eq:ordinary-list} to $u(a,b)=\tau(b\Phi(a))$ in \eqref{eq:22norm} gives the elementary upper bound $2$. 
		We prove the sharper conclusions by the circularization argument.
		
		Normalize $\|\Phi\|_{2\to2}=1$. If the factors are finite dimensional, irreducibility forces $\cN=\cM$, and $\Phi$ is scalar. 
		In the $\mathrm{II}_1$ case, $\Phi(1)$ and $\Phi^*(1)$ are $\cN$-central $L^2$ vectors, hence $\Phi(1)=\beta1$ and $\Phi^*(1)=\overline\beta1$. 
		The scalar line and its orthogonal complement reduce $\Phi$. 
		If $|\beta|=1$, the test $u=v=1$ gives $\|\Phi\|_{\infty\to1}=1$. 
		Otherwise, given $\varepsilon>0$, choose centered $x,y\in\cM$ of $L^2$ norm one such that $\operatorname{Re}\langle y,\Phi(x)\rangle>1-\varepsilon^2/2$. 
		Expanding squared norms and using $\|\Phi\|_{2\to2}=1$ gives
		\[
		\|\Phi(x)-y\|_2<\varepsilon,\qquad
		\|\Phi^*(y)-x\|_2<\varepsilon .
		\]
		If $\Phi\ge0$, we may take $y=x$ by choosing $\langle x,\Phi(x)\rangle>1-\varepsilon^2/2$.
		
		Let $\omega$ be a free ultrafilter on $\bN$. 
		We compress the coordinatewise operator in the Hilbert ultrapower of $L^2(\cM)$ to $L^2(\cM^\omega)$. 
		This subspace reduces both left and right multiplication by $\cN^\omega$, so the compression $\Phi^\omega$ is bimodular and contractive. 
		It preserves positivity of $\Phi$ and restricts to $\Phi$ on constant $L^2$ vectors. 
		Its adjoint similarly extends $\Phi^*$. 
		Unitaries in $\cM^\omega$ lift to sequences of unitaries, by taking unitary polar extensions of uniformly bounded representatives. 
		Thus \eqref{eq:1norm} yields
		\[
		\|\Phi^\omega\|_{\infty\to1}
		=\sup_{u,v\in\U(\cM^\omega)} |\langle v,\Phi^\omega(u)\rangle|
		=\|\Phi\|_{\infty\to1}.
		\]
		Indeed the pairings are ultralimits of the original pairings, and constant sequences give the reverse inequality. 
		If $\|\Phi\|_{\infty\to\infty}\le\kappa$, then on bounded representatives the compression is the coordinatewise operator, so $\|\Phi^\omega\|_{\infty\to\infty}\le\kappa$.
		
		Popa's free-independence theorem \cite[Theorem~0.1(a)]{PopaFree14} provides a Haar unitary $w\in\cN^\omega$ free from the separable algebra generated by $x,y$ and bounded truncations of the residuals $\Phi(x)-y$, $\Phi^*(y)-x$. 
		The relative commutant of $\cN^\omega$ in $\cM^\omega$ is scalar: the convex-hull argument of Lemma~\ref{lem:norm} gives
		\[
		\|z-\tau(z)1\|_2\le \sup_{u\in\U(\cN)}\|uzu^*-z\|_2,
		\]
		and choosing a near maximizing unitary separately in each coordinate proves the claim. 
		The diffuse factor $\cN$ cannot intertwine into its scalar relative commutant, as required by the theorem. 
		Consecutive distinct powers in an alternating product of conjugated centered elements leave nonzero powers of $w$ between the elements. 
		Freeness with $w$ then shows that the conjugated algebras $w^j\mathcal A w^{-j}$ are free. 
		In particular, conjugates with distinct $j$ of any centered residual are orthogonal in $L^2$, by bounded truncation and continuity.
		
		For even $n$ put
		\[
		\mathbf x_n=\frac1{\sqrt n}\sum_{j=1}^n\varepsilon_jw^j x w^{-j},
		\qquad
		\mathbf y_n=\frac1{\sqrt n}\sum_{j=1}^n\varepsilon_jw^j y w^{-j},
		\qquad
		\varepsilon_j=\begin{cases}1,&j\text{ odd},\\ i,&j\text{ even}.\end{cases}
		\]
		Bimodularity and the residual orthogonality give
		\[
		\|\Phi^\omega(\mathbf x_n)-\mathbf y_n\|_2<\varepsilon,
		\qquad
		\| (\Phi^\omega)^*(\mathbf y_n)-\mathbf x_n\|_2<\varepsilon .
		\]
		Each marginal converges in $*$-moments to a normalized circular element. 
		In fact, for any fixed word of length $k$ in the element and its adjoint, the $k$th free cumulant of the normalized sum is the sum of $n$ terms of size $O(n^{-k/2})$. 
		The first cumulants vanish, the mixed second cumulants equal one, and the second cumulants with equal signs vanish by $\displaystyle \sum_j\varepsilon_j^2=0$. 
		All cumulants of order $k\ge3$ tend to zero, which proves the assertion by the moment--cumulant formula \cite{NS}.
		
		The distributions of $|\mathbf x_n|$ and $|\mathbf y_n|$ therefore converge weakly to the quarter-circle law $\nu$ on $[0,2]$, since the squared absolute values converge in moments to a compactly supported law. 
		Clip both polar moduli at $2$. 
		The squared clipping error is the integral of $(s-2)_+^2$ against the corresponding law and tends to zero:
		the fourth moments are uniformly bounded, so above a threshold $\alpha$ the integral is bounded by $\alpha^{-2}$ times the fourth moment.
		A second tracial ultrapower of the clipped sequences gives $X=U|X|$ and $Y=V|Y|$ with $|X|,|Y|$ of law $\nu$. 
		Continue to write $\Phi^\omega$ for the resulting compression, with the same bounds and unitary tests. 
		The two residuals still have norm at most $\varepsilon$. The polar parts $U,V$ are unitaries, since both moduli have zero kernel and the algebra is finite. In the positive case $Y=X$.
		
		Set $f(s)=s/2$, $g(s)=\sqrt{1-s^2/4}$, and put
		\[
		F_X=Uf(|X|),\qquad G_X=Ug(|X|),\qquad U_\pm=F_X\pm iG_X.
		\]
		Define $F_Y,G_Y,V_\pm$ from $Y,V$ in the same way. 
		The identity $f^2+g^2=1$ makes these four test elements unitary. 
		Since $F_X=X/2$ and $F_Y=Y/2$, the two residual bounds give
		\[
		\langle F_Y,\Phi^\omega(G_X)\rangle
		=\frac12\int_0^2 s g(s)\,d\nu(s)+O(\varepsilon),
		\qquad
		\langle G_Y,\Phi^\omega(F_X)\rangle
		=\frac12\int_0^2 s g(s)\,d\nu(s)+O(\varepsilon).
		\]
		The integral is $\displaystyle (2\pi)^{-1}\int_0^2 s(4-s^2)\,ds=2/\pi$. 
		If $\beta_F$ denotes the difference of the $F$--$F$ and $G$--$G$ terms and $\beta_G$ the sum of the two mixed terms, then
		\[
		\langle V_+,\Phi^\omega(U_-)\rangle=\beta_F-i\beta_G,
		\qquad
		\langle V_-,\Phi^\omega(U_+)\rangle=\beta_F+i\beta_G,
		\qquad \beta_G=2/\pi+O(\varepsilon).
		\]
		The maximum of the two absolute values is at least $|\beta_G|$. 
		Let $\varepsilon\downarrow0$.
		
		When $\Phi\ge0$, use the unitary polar part $U$ of the approximate eigenvector $X$. 
		Its overlap is
		\[
		\langle X,U\rangle=\frac1\pi\int_0^2s\sqrt{4-s^2}\,ds
		=\frac8{3\pi}.
		\]
		The denominator in the following expression lies in $(0,1]$ for $\varepsilon<8/(3\pi)$, and the numerator tends to the square of this overlap.
		The Cauchy--Schwarz inequality for the positive form defined by $\Phi^\omega$ gives
		\[
		\|\Phi\|_{\infty\to1}\ge\langle U,\Phi^\omega(U)\rangle
		\ge\frac{|\langle X,\Phi^\omega(U)\rangle|^2}{\langle X,\Phi^\omega(X)\rangle}
		\ge\left(\frac8{3\pi}-\varepsilon\right)^2 .
		\]
		This proves the claimed constant.
	\end{proof}

	\begin{remark}
		The same numerical factor $64/(9\pi^2)$ appears as the limiting squared mean singular value in Bandeira--Kennedy--Singer's rounding analysis for the little Grothendieck problem over unitary groups \cite{BKS}.
	\end{remark}

	\begin{remark}
		The irreducibility of the inclusion is crucial here.
		Bannink--Bri\"et--Labib--Maassen \cite{BBLM} prove, for irreducibly covariant matrix maps with normalized matrix traces,
		\[
		\|\Phi\|_{\infty\to1}\le\|\Phi\|_{2\to2} \le2\|\Phi\|_{\infty\to1}.
		\]
		Their paper proves the optimality of two for that class.
	\end{remark}
	
	\begin{remark}
		Suppose that $\Phi:\cM\to\cM$ is a completely positive bimodule map on an irreducible inclusion.
		Bimodularity gives $\Phi(1)\in\cN'\cap\cM=\bC1$. Complete positivity gives
		$\|\Phi\|_{\cb}=\|\Phi(1)\|$. 
		Testing at $1$ gives $\|\Phi\|_{\infty\to1}\ge\tau(\Phi(1))=\|\Phi(1)\|$ and the opposite inequality follows from $\|\Phi\|_{\infty\to1}\le\|\Phi\|_{\cb}$.
	\end{remark}
	
	\begin{remark}
		By Theorem \ref{gb:theorem}, we have that
		\begin{align*}
			1\leq \KI\leq \frac{\pi}{2}, \quad 1\leq \KI_+\leq \frac{9\pi^2}{64}.
		\end{align*}
	\end{remark}

	\begin{remark}
		Let $\Phi$ be a trace-symmetric unital completely positive $\cN$-bimodule map on an irreducible inclusion of $\mathrm{II}_1$ factors. 
		The bimodule max bisection problem is to estimate
		\[
		\mathsf B(\Phi)=\sup_{\substack{p=p^*=p^2\\\tau(p)=1/2}}2\tau(p\Phi(1-p)).
		\]
		For the random walk on a finite regular graph, the corresponding objective at $p=1_S$ is the fraction of edges crossing $(S,S^c)$, and $\tau(p)=1/2$ is the bisection constraint.
		The operator $I-\Phi$ is positive on $L^2(\cM)$:
		the Schwarz inequality for $\Phi$ and trace preservation imply contractivity on $L^2(\cM)$, and trace symmetry gives self-adjointness. 
		Since $(I-\Phi)(1)=0$, the trace-zero symmetry $u=2p-1$ satisfies
		\[
		\tau(u(I-\Phi)(u))=4\tau(p\Phi(1-p)).
		\]
		This gives the upper estimate below.
		
		For the lower estimate, choose matrix units in $\cN$ with $e=e_{11}$ of trace $1/2$. 
		Let $\Psi=I-\Phi$.
		Under the associated identification $\cM=M_2(e\cM e)$, bimodularity gives $\Psi=\mathrm{id}_2\otimes \Psi_e$, where $\Psi_e$ is the restriction to the corner. 
		The corner inclusion is irreducible and $\|\Psi_e\|_{2\to2}=\|\Psi\|_{2\to2}$. 
		For $u\in\U(e\cM e)$,
		\[
		s=\begin{pmatrix}0&u\\u^*&0\end{pmatrix},
		\qquad
		\tau(s\Psi(s))=\frac12\bigl(\tau_e(u\Psi_e(u^*))+\tau_e(u^*\Psi_e(u))\bigr)
		=\tau_e(u^*\Psi_e(u)),
		\]
		where $\tau_e=2\tau|_{e\cM e}$. 
		The last equality follows from $*$-preservation and positivity of the $L^2$ pairing. 
		This $s$ has trace zero. 
		By the Cauchy--Schwarz inequality for the positive form defined by $\Psi_e$, the supremum of $\tau_e(u^*\Psi_e(u))$ over corner unitaries equals $\|\Psi_e\|_{\infty\to1}$. 
		Theorem~\ref{gb:theorem} therefore gives
		\[
		\frac{32}{9\pi^2}\,\|I-\Phi\|_{2\to2}
		\le\mathsf B(\Phi)
		\le\frac12\|I-\Phi\|_{2\to2}.
		\]
	\end{remark}

	\begin{conjecture}
		In Section \ref{sec:positive}, we shall show that the sharp Grothendieck constant for positive bimodule operator is $\displaystyle \KI_+=\frac{9\pi^2}{64}$.
		For the general case, we would like to conjecture that
		\begin{align*}
			\KI=\frac{\pi}{2}.
		\end{align*}
	\end{conjecture}
	
	\subsection{Completely bounded Grothendieck constant}
	In this section, we shall obtain a completely bounded version of Grothendieck inequality.
	
	\begin{theorem}\label{thm:cb-algebra}
		Suppose $\cN\subset\cM$ is an irreducible inclusion of finite factors with index $\delta^2<\infty$. Then
		\begin{equation}\label{eq:gbcb}
			\KNM{\cN}{\cM}\le\KNM{\cN}{\cM}_{\cb}\le\frac{\pi\delta}{2}.
		\end{equation}
		For an operator-positive bimodule map,
		\begin{equation}\label{eq:gbcbp}
			\|\Phi\|_{\cb}\le\frac{9\pi^2\delta}{64}\|\Phi\|_{\infty\to1}.
		\end{equation}
		There are outer crossed-product inclusions of arbitrarily large finite index with a nonzero bimodule map satisfying
		$\|\Phi\|_{\cb}=\delta\|\Phi\|_{2\to2}=\delta\|\Phi\|_{\infty\to1}$.
	\end{theorem}
	
	\begin{proof}
		Let $\alpha=\|\Phi\|_{2\to2}$. 
		For $x\in\cM$ and $a\in\cN$, bimodularity gives
		\[
		\|\Phi(x)a\|_2^2=\|\Phi(xa)\|_2^2\le \alpha^2\|xa\|_2^2.
		\]
		Trace duality on $\cN$, testing against all $aa^*$, gives
		\begin{align*}
			\bE_{\cN}(\Phi(x)^*\Phi(x))\le \alpha^2\bE_{\cN}(x^*x)
		\end{align*}
		in $L^1(\cN)$. 
		For each bounded spectral truncation $h$ of $\Phi(x)^*\Phi(x)$, the Pimsner--Popa inequality gives
		\[
		h\le\delta^2\bE_{\cN}(h)
		\le\delta^2\alpha^2\bE_{\cN}(x^*x)
		\le\delta^2\alpha^2\|x\|_\infty^2 1.
		\]
		Taking the increasing limit shows $\Phi(x)\in\cM$ and $\|\Phi\|_{\infty\to\infty}\le\delta \alpha$. Lemma~\ref{lem:map-cb-norm} gives $\|\Phi\|_{\cb}=\|\Phi\|_{\infty\to\infty}$.
		
		Conversely, for any lists in \eqref{eq:22norm}, the column bound
		\[
		\left(\sum_j\|\Phi(x_j)\|_2^2\right)^{1/2}
		\le \left\|\sum_j\Phi(x_j)^*\Phi(x_j)\right\|_\infty^{1/2}
		\le\|\Phi\|_{\cb}
		\]
		and $\displaystyle \sum_j\|y_j\|_2^2\le1$ give $\alpha\le\|\Phi\|_{\cb}$. 
		Combining both estimates with Theorem~\ref{gb:theorem} gives \eqref{eq:gbcb} and \eqref{eq:gbcbp}.
		
		For sharpness of the index factor, take $\cM=\cN\rtimes G$ for an outer action of $G=(\bZ_2)^r$. 
		Let $e_j$ be the coordinate generators, $u_g$ the implementing unitaries, and $|g|$ the number of nonzero coordinates. 
		Define $\Phi(a_gu_g)=i^{|g|}a_gu_g$. 
		The Fourier summands are orthogonal, so $\|\Phi\|_{2\to2}=1$. 
		Testing at $1$ gives $\|\Phi\|_{\infty\to1}=1$.
		As $u_{e_j}$ are commuting self-adjoint unitaries,
		\[
		u=\delta^{-1}\sum_g i^{-|g|}u_g
		=\prod_{j=1}^r\frac{1-iu_{e_j}}{\sqrt2}
		\quad\text{is unitary},\qquad
		\Phi(u)=\delta e,\quad e=|G|^{-1}\sum_g u_g.
		\]
		The group average $e$ is a projection of trace $|G|^{-1}$, so $\|\Phi\|_{\cb}\ge\delta$. 
		The upper bound $\|\Phi\|_{\cb}\le\delta \alpha$ gives equality.
	\end{proof}
	
	\begin{remark}
		For $G=(\bZ_2)^r\times(\bZ_2)^r$, take $\sigma(a,b)=(-1)^{a\cdot b}$ and let $S_\sigma$ have Fourier coefficients $\sigma$. 
		The Fourier sum
		\[
		\sum_{a,b}(-1)^{a\cdot b+k\cdot a+\ell\cdot b}=2^r(-1)^{k\cdot\ell}
		\]
		shows that $\displaystyle u_\sigma=\delta^{-1}\sum_g\sigma(g)u_g$ takes values in $\{1,-1\}$ under every character of the group algebra. 
		Hence $u_\sigma$ is a self-adjoint unitary and $eu_\sigma=e$. 
		As $S_\sigma(u_\sigma)=\delta e$, we have $\|S_\sigma\|_{\cb}=\delta$. 
		The map $\Phi_+=(\id+S_\sigma)/2$ is an orthogonal projection on $L^2(\cM)$ and $\|\Phi_+\|_{\infty\to1}=1$. 
		Evaluating on $u_\sigma$ and compressing to $e$ gives
		\[
		\|\Phi_+\|_{\cb}\ge\frac{1+\delta}{2}.
		\]
		The triangle inequality gives equality. 
		Thus the usual completely bounded norm can grow linearly in $\delta$ even for operator-positive bimodule maps.
	\end{remark}
	
	For the downward basic construction for irreducible finite-index inclusion of $\mathrm{II}_1$ factors, we refer to  \cite{Jones1983,PPindex}.
	
	\begin{proposition}
		Let $\cN\subset\cM$ be an irreducible finite-index inclusion of
		$\mathrm{II}_1$ factors.
		Suppose that $\mathcal P\subset\cN\subset\cM=\langle\cN,e_{\mathcal P} \rangle$ is a downward basic construction.
		For any bounded $\cN$-bimodule operator $\Phi$ on $L^2(\cM)$, let
		$\Phi^\flat(x)=(\Phi^*(x^*))^*$ for $x\in \cM$.
		Then
		\begin{equation}\label{eq:downward-cb-trace}
			\|\Phi\|_{1\to1}=\|\Phi^\flat\|_{\infty\to\infty}
			=\|\Phi\|_{\cb}=\|\Phi\|_{\infty\to\infty}=\delta^2 \tau(|\Phi(e_{\mathcal{P}})|).
		\end{equation}
	\end{proposition}
	\begin{proof}
		Take a finite right Pimsner--Popa basis $(v_j)_j$ for $\mathcal P\subset\cN$, so that $\displaystyle \sum_j v_j e_{\mathcal{P}} v_j^*=1$. 
		Let $c_{ij}(x)\in\mathcal P$ be determined by
		$e_{\mathcal{P}} v_i^* x v_j e_{\mathcal{P}}=c_{ij}(x)e_{\mathcal{P}}$. 
		The matrix $C(x)=(c_{ij}(x))$ is a $*$-homomorphism into a corner of $M_k(\mathcal P)$, where $k$ is the number of basis elements. 
		For multiplication, insert the frame identity between $x$ and $z$ in $e_{\mathcal{P}} v_i^* xz v_je_{\mathcal{P}}$ to get $\displaystyle c_{ij}(xz)=\sum_l c_{il}(x)c_{lj}(z)$. 
		The adjoint identity follows by interchanging $i$ and $j$.
		
		For $y\in\mathcal P'\cap\cM$ define $\displaystyle \Phi_y(x)=\sum_{i,j}v_i c_{ij}(x)y v_j^*$. 
		On $x=a e_{\mathcal{P}} b$ with $a,b\in\cN$, the coefficients are $c_{ij}(x)=\bE_{\mathcal P}(v_i^*a)\bE_{\mathcal P}(b v_j)$.
		Since $y$ commutes with $\mathcal P$, the two basis expansions give
		$\Phi_y(a e_{\mathcal{P}} b)=a y b$. 
		The frame expresses every $x$ as a finite sum of such elements. 
		Hence $\Phi_y$ is a normal bimodule map, and
		$\Phi=\Phi_y$ for $y=\Phi(e_{\mathcal{P}})$, which lies in $\mathcal P'\cap\cM$.
		The same statements apply to $\Phi^\flat$, which is bimodular and bounded on $\cM$ by the finite-index norm estimate in Theorem~\ref{thm:cb-algebra}.
		
		By irreducibility $\Phi_y(1)$ is scalar. 
		The functional $\tau\circ \Phi_y$ is $\cN$-central by bimodularity, so its $L^1$ density is scalar as well. 
		Evaluation at $e_{\mathcal{P}}$ gives $\tau\circ \Phi_y=\delta^2\tau(y)\tau$, and evaluating at $1$ then gives $\Phi_y(1)=\delta^2\tau(y)1$.
		Let $V=(v_1,\ldots,v_k)$ and write $\operatorname{diag}(b)$ for the
		$k$-by-$k$ diagonal matrix with every diagonal entry $b$. 
		For the polar decomposition $y=w|y|$ we have
		\[
		\Phi_y(x)=V\operatorname{diag}(|y^*|^{1/2})\,C(x)\,
		\operatorname{diag}(w)\operatorname{diag}(|y|^{1/2})V^*.
		\]
		The projection $C(1)$ commutes with each displayed diagonal matrix
		and $C(x)=C(1)C(x)C(1)$. 
		The squared norms of the two outside factors, with the corresponding $C(1)$ inserted, are therefore $\|\Phi_{|y^*|}(1)\|$ and $\|\Phi_{|y|}(1)\|$. 
		This holds at each matrix amplification, while $C$ is completely contractive and $\|w\|\le1$.
		It follows that $\|\Phi\|_{\cb}\le \delta^2\tau(|\Phi(e_{\mathcal{P}})|)$.
		
		For contractions $x\in\cM$, trace duality and $\tau(e_{\mathcal{P}})=\delta^{-2}$ give
		\[
		|\tau(\Phi^\flat(e_{\mathcal{P}})x)|=|\tau(e_{\mathcal{P}} \Phi(x))|\le \delta^{-2}\|\Phi\|_{\infty\to\infty}.
		\]
		Taking the supremum gives $\delta^2\tau(|\Phi^\flat(e_{\mathcal{P}})|)\le\|\Phi\|_{\infty\to\infty}$.
		Apply both bounds to $\Phi$ and to $\Phi^\flat$. 
		The resulting chain is
		\begin{align*}
			\delta^2 \tau(|\Phi^\flat(e_{\mathcal{P}})|)&\le \|\Phi\|_{\infty\to\infty}\le \|\Phi\|_{\cb}
			\le \delta^2\tau(|\Phi(e_{\mathcal{P}})|)\\
			&\le \|\Phi^\flat\|_{\infty\to\infty}\le \|\Phi^\flat\|_{\cb}
			\le \delta^2\tau(|\Phi^\flat(e_{\mathcal{P}})|).
		\end{align*}
		Thus all inequalities are equalities, and the $L^1$ equality follows by trace duality.
	\end{proof}
	
	\subsection{Sharp Grothendieck constant for positive bimodule operators}\label{sec:positive}
	
	In this section, we shall show that the Grothendieck constant for positive bimodule operators is $\displaystyle \KI_+=\frac{9\pi^2}{64}$.
	
	\begin{proposition} \label{prop:tracial-free-index}
		Suppose that $\cN^{(j)}\subset \cM^{(j)}$ are irreducible inclusions  of $\mathrm{II}_1$ factors, where $j=1,2$.
		Then the inclusion of the reduced tracial free product
		\[
		\cN=\cN^{(1)}*\cN^{(2)}\subset \cM=\cM^{(1)}*\cM^{(2)}
		\]
		is irreducible.
		Its index is infinite as soon as one constituent inclusion is proper.
	\end{proposition}
	\begin{proof}
		Write $H_j^\circ=L^2(\cM^{(j)})\ominus\bC1$. 
		In the reduced-word decomposition of $L^2(\cM)$, the orthogonal complement of $L^2(\cM^{(1)})$ is, as a bimodule over $\cM^{(1)}$,
		\[
		L^2(\cM^{(1)})\otimes\mathcal K\otimes L^2(\cM^{(1)}),
		\]
		where $\mathcal K$ is the direct sum of the alternating tensor words beginning and ending with $H_2^\circ$. 
		The outside tensor factors include the scalar components to account for the optional first and last letters.
		If $u\in\cN^{(1)}$ is a Haar unitary, its powers tend weakly to zero. 
		On the displayed tensor product, $u^n$ acts on the first factor and $u^{-n}$ on the last, so the conjugation operators tend weakly to zero, first on simple tensors and then on the whole complement by contractivity. 
		No nonzero vector in this complement can be $\cN^{(1)}$-central. 
		An element of $\cN'\cap\cM$ therefore lies in $\cM^{(1)}$, and then in $\bigl(\cN^{(1)}\bigr)'\cap\cM^{(1)}=\bC1$. 
		This also implies that both free products are factors, since their centers lie in $\cN'\cap\cM$.
		
		Assume $\cN^{(1)}\ne\cM^{(1)}$. 
		Choose $0\ne x\in\ker\bE_{\cN^{(1)}}$ and a Haar unitary $v\in\cN^{(2)}$. 
		Put $w_j=v^jx$ for $j\ge1$. 
		On reduced words, $\bE_{\cN}$ applies the two constituent expectations letter by letter. 
		This follows since the inclusion of the free product $L^2$ spaces is the inclusion of the alternating tensor words of the two smaller centered $L^2$ spaces, and the orthogonal projection on each such word is the tensor of the two constituent projections. 
		Thus
		\[
		\bE_{\cN}(w_j^*w_k)=
		\begin{cases}
			a, & j=k,\\
			0, & j\ne k,
		\end{cases}
		\qquad a=\bE_{\cN^{(1)}}(x^*x)\ne0.
		\]
		For $j\ne k$, the word $x^*v^{k-j}x$ is centered in each letter, and its outside letters have zero constituent expectation. 
		Hence the closed spaces $\overline{w_j\cN}^{\,L^2}$ are mutually orthogonal right $\cN$-modules. 
		On $L^2(\cN)$ the square modulus of the bounded map $b\mapsto w_jb$ is left multiplication by $a$, by the displayed expectation identity. 
		Its polar decomposition identifies the module $\overline{w_j\cN}^{\,L^2}$ with $\cR(a)L^2(\cN)$, of dimension $\tau(\cR(a))>0$. 
		Summing their dimensions gives infinite index. 
		The case of a proper second constituent follows by interchanging the factors.
	\end{proof}

	\begin{lemma}\label{lem:semicircle-first-moment-product}
		Let $\mathcal A$ be a von Neumann subalgebra of a finite tracial von Neumann algebra $(\mathcal P,\tau)$ with normalized trace $\tau$ and let $y=y^*$ be a bounded, centered $\mathcal A$-valued semicircular element.
		Then
		\[
		\tau(|y|)\le\frac8{3\pi}\tau(y^2)^{1/2}.
		\]
	\end{lemma}
	\begin{proof}
		Put $\alpha=\tau(y^2)$. 
		We may assume $\alpha>0$. 
		The covariance $E_y(a)=\bE_{\mathcal A}(yay)$ is completely positive and trace-symmetric, since
		$\tau(aE_y(b))=\tau(ayby)=\tau(E_y(a)b)$ for $a,b\in\mathcal A$. 
		For $t>0$ set
		\[
		y_t=\bE_{\mathcal A}\bigl(t(t^2+y^2)^{-1}\bigr).
		\]
		It is positive and bounded below by $t/(t^2+\|y\|_\infty^2)$. 
		The operator-valued semicircular Cauchy-transform equation \cite{HeltonFarSpeicher} gives
		\[
		y_t^{-1}=t1+E_y(y_t).
		\]
		In detail, $G(z)=\bE_{\mathcal A}((z-y)^{-1})$ obeys $(z-E_y(G(z)))G(z)=1$ by the noncrossing pairing recursion, with the equality on the upper half-plane following by analytic continuation. 
		Evenness of the operator-valued law gives $G(it)=-iy_t$, which yields the claimed identity.
		
		For positive invertible $a\in\mathcal A$ define
		\[
		\Phi_t(a)=t\tau(a)+\frac12\tau(aE_y(a))-\tau(\log a).
		\]
		In a self-adjoint direction $h$, the derivative of the last term is $-\tau(a^{-1}h)$, as is seen by differentiating the resolvent integral for $\log a$ and using traciality. 
		The full derivative is $\tau(h[t1+E_y(a)-a^{-1}])$. 
		Thus $y_t$ is stationary. 
		Let $\beta_t>0$ solve $\beta_t^{-1}=t+\alpha\beta_t$, and consider the path $a_s=\beta_t^s y_t^{1-s}$, $0\le s\le1$. 
		Along it, $\tau(a_s)$ is convex by the scalar spectral theorem, and $\tau(\log a_s)$ is affine. 
		For the quadratic term write
		\[
		\tau(a_s E_y(a_s))=\langle y,L_{a_s}R_{a_s}y\rangle,
		\]
		with $L$ and $R$ denoting left and right multiplication. 
		The joint spectral measure of $L_{y_t},R_{y_t}$ against $y$ expresses this quantity as an integral of $\beta_t^{2s}(uv)^{1-s}$ against a positive measure on $(0,\infty)^2$. 
		Each integrand is convex in $s$. 
		Stationarity at $s=0$ now gives $\Phi_t(y_t)\le\Phi_t(\beta_t1)$.
		
		Differentiating along the stationary solutions gives
		\[
		\frac{d}{dt}\Phi_t(y_t)=\tau(y_t)=\frac{d}{dt}\frac12\tau(\log(t^2+y^2)).
		\]
		As $t\to\infty$, $y_t=t^{-1}1+O(t^{-3})$ in operator norm. 
		The two expressions differ by a constant, which this asymptotic identifies as $1$. 
		The same calculation for a scalar semicircular variable $y_0$ of variance $\alpha$, for which the solution of the resolvent equation is $\beta_t$, shows
		\begin{equation}\label{eq:product-semicircle-log}
			\tau\bigl(\log(1+y^2/t^2)\bigr)
			\le \tau\bigl(\log(1+y_0^2/t^2)\bigr)
			\qquad(t>0).
		\end{equation}
		For $a>0$, integration by parts gives
		\[
		\int_0^\infty\log(1+a^2/t^2)\,dt
		=2a^2\int_0^\infty(t^2+a^2)^{-1}\,dt=\pi a .
		\]
		The formula also holds for $a=0$. 
		Applying Tonelli theorem to \eqref{eq:product-semicircle-log} gives $\displaystyle \tau(|y|)\le\tau(|y_0|)=\frac{8}{3\pi}\sqrt\alpha$.
	\end{proof}
	
	\begin{theorem}\label{thm:free-product-positive-sharp}
		Suppose $\cN\subset\cM$ is an irreducible inclusion of $\mathrm{II}_1$ factors and $\cM$ contains a centered semicircular element $a$ of nonzero variance free from $\cN$. 
		Then the orthogonal projection $P_a$ onto $\overline{\operatorname{span}(\cN a\cN)}^{\,L^2(\cM)}$ satisfies
		\[
		\|P_a\|_{2\to2}=1,\qquad \|P_a\|_{\infty\to1}=\left(\frac8{3\pi}\right)^2.
		\]
		Consequently, $\KI_+=9\pi^2/64$.
	\end{theorem}
	\begin{proof}
		Let $\mathcal H_a=\overline{\operatorname{span}(\cN a\cN)}^{\,L^2(\cM)}$. 
		For a finite sum $\displaystyle x=\sum_i b_i a c_i$ with $b_i,c_i\in\cN$, set
		$\displaystyle Y= \begin{pmatrix}0&x\\x^*&0\end{pmatrix}$. 
		Scalar freeness and the Wick recursion give the $\cN$-valued semicircular law of $a$, with covariance $b\mapsto\tau(a^2)\tau(b)1$. 
		The same recursion entry by entry gives the $M_2(\cN)$-valued law of $a\otimes 1_2$. 
		Its cumulants of orders other than two vanish, also after inserting coefficients on either side. 
		Since the upper right block of $Y$ consists of the terms $(e_{11}\otimes b_i)(1_2\otimes a)(e_{12}\otimes c_i)$ and the lower left block is its adjoint, $Y$ is centered and semicircular over $M_2(\cN)$. 
		Normalized matrix trace gives
		\[
		Y^2=\begin{pmatrix}xx^*&0\\0&x^*x\end{pmatrix},
		\qquad
		\tau_{M_2(\cM)}(|Y|)=\|x\|_1,\qquad
		\tau_{M_2(\cM)}(Y^2)=\|x\|_2^2.
		\]
		Lemma~\ref{lem:semicircle-first-moment-product} gives $\displaystyle \|x\|_1\le\frac{8}{3\pi}\|x\|_2$. 
		The inequality extends to $\mathcal H_a$, and equality holds for $x=a$. 
		The same hermitianization in \eqref{eq:product-semicircle-log} gives the logarithmic comparison used below.
		
		The subspace $\mathcal H_a$ reduces both $\cN$ actions, so $P_a$ is a bimodule projection. 
		The Cauchy--Schwarz inequality and then Hilbert-space and trace duality give
		\begin{align*}
			\|P_a\|_{\infty\to1}
			&=\sup_{\|z\|_\infty\le1}\left\langle z,P_a z\right\rangle
			=\left(\sup_{\|z\|_\infty\le1}\|P_a z\|_2\right)^2\\
			&=\left(\sup_{\substack{x\in\mathcal H_a\\\|x\|_2\le1}}\|x\|_1\right)^2
			=\left(\frac8{3\pi}\right)^2.
		\end{align*}
		
		For an example, take a constituent crossed product $\cM^{(1)}=\cN^{(1)}\rtimes\bZ_2$ by an outer action, with implementing unitary $u$, in Proposition~\ref{prop:tracial-free-index}. 
		The algebra $u\cN^{(2)}u^*$ is free from $\cN^{(1)}*\cN^{(2)}$ inside $\cM^{(1)}*\cM^{(2)}$. 
		To verify this, write centered letters in $u\cN^{(2)}u^*$ as $ubu^*$ with $\tau(b)=0$, and use centered reduced words for the intervening elements of $\cN^{(1)}*\cN^{(2)}$. 
		At an end adjacent to $u$ or $u^*$, a letter $c\in\cN^{(1)}$ merges to $cu$ or $u^*c$, whose trace vanishes because $\bE_{\cN^{(1)}}(u)=0$.
		If the whole intervening word is one letter $c\in\cN^{(1)}$, the merged letter is $u^*cu$, centered by traciality. 
		Adjacent letters from $\cN^{(2)}$ are separated by $u$ or $u^*$ or by one of these centered products. 
		Thus the resulting word is reduced and centered. 
		The alternating trace vanishing extends to the generated von Neumann algebras by bounded strong approximation. 
		A diffuse abelian subalgebra of $u\cN^{(2)}u^*$ supplies a semicircular variable of variance one. Proposition~\ref{prop:tracial-free-index} gives irreducibility, and Theorem~\ref{gb:theorem} gives the matching universal upper bound.
	\end{proof}
	
	\begin{corollary}\label{cor:free-product-indefinite}
		Under the hypotheses of Theorem~\ref{thm:free-product-positive-sharp}, the bimodule operator $\displaystyle \Phi=\frac{10}{9}P_a-\frac19\id$ satisfies
		\[
		\|\Phi\|_{2\to2}=1,\qquad \|\Phi\|_{\infty\to1}<\frac{3531}{5000}.
		\]
		In particular, $\sqrt2<5000/3531<\KI\le\pi/2$.
	\end{corollary}
	\begin{proof}
		Let $\displaystyle \alpha=\frac{\pi}{10(2-\log 3)}$ and define, for $t\ge0$,
		\[
		f(t)=\frac1{10}+\frac{\alpha}{\pi}\log(1+3t/\alpha^2)
		+\frac1\pi\int_\alpha^\infty\log(1+t/r^2)\,dr .
		\]
		The function is concave, since each logarithm is concave in $t$. 
		By the integral formula for $\sqrt{t}$ from the proof of Lemma~\ref{lem:semicircle-first-moment-product},
		\[
		\Delta(t):=f(t)-\sqrt{t}
		=\frac1\pi\int_0^\alpha\log\frac{t+\alpha^2/3}{t+r^2}\,dr .
		\]
		For this identity, use
		$\displaystyle \pi^{-1}\int_0^\alpha\log[(\alpha^2/3)/r^2]dr =\alpha(2-\log 3)/\pi=1/10$.
		Jensen's inequality for $\log$ gives $\Delta(t)\ge0$ since $\alpha^2/3$ is the average of $r^2$ on $[0,\alpha]$. 
		For $t>0$, Jensen inequality applied to the reciprocal gives $\Delta'(t)\le0$, and hence $\Delta(t)\le\Delta(0)=1/10$. 
		Furthermore,
		\[
		\int_0^\infty\Delta(t^2)\,dt
		=\int_0^\alpha(\alpha/\sqrt3-r)\,dr
		=\alpha^2\left(\frac1{\sqrt3}-\frac12\right).
		\]
		The Fubini theorem applies because for each $r$ the logarithmic ratio has a constant sign in $t$ and its absolute integral is $\pi|\alpha/\sqrt3-r|$. 
		For a finite sum $x\in\operatorname{span}(\cN a\cN)$, the logarithmic comparison \eqref{eq:product-semicircle-log} applied to its hermitianization gives
		\begin{equation}\label{eq:indefinite-log-f}
			\tau(f(|x|^2))\le\int_0^2 f(\|x\|_2^2r^2)\,d\nu(r).
		\end{equation}
		Here, each scalar comparison is at the variance $\|x\|_2^2$, so the absolute scalar semicircle has law $\|x\|_2\,\nu$ in the sense of scaling the variable.
		
		Given unitaries $u,v$, multiply $v$ by a scalar phase so that $\tau(v^*\Phi(u))\ge0$, and put $w=(u+v)/2$, $\beta=\|P_a(w)\|_2$ and $t=\|w\|_2^2$. 
		The real part of the Hilbert pairing yields
		\begin{align}
			|\tau(v^*\Phi(u))|
			&=\frac{10}{9}\left(\beta^2-\frac14\|P_a(v-u)\|_2^2\right)
			-\frac19(2\|w\|_2^2-1)\notag\\
			&\le \frac{10}{9}\left(\beta^2-\frac {\|w\|_2^2}{5}\right)+\frac19 .
			\label{eq:indefinite-wedge-midpoint}
		\end{align}
		For $r\ge0$ let $\displaystyle \psi(r)=\max_{-1\le s\le1}(2rs-s^2/5)$. 
		Then $\psi(r)=5r^2$ for $0\le r\le1/5$ and $\psi(r)=2r-1/5$ for $r\ge1/5$. 
		On the latter interval, $f(r^2)\ge r$. 
		On the former interval, concavity gives $f(r^2)\ge 1/10+(5/2)r^2$ from $f(0)=1/10$ and $f(1/25)\ge 1/5$. 
		Thus $\psi(r)\le2f(r^2)-1/5$ on both intervals.
		
		For any finite sum $x$ as above, write $x=u_x|x|$ using a unitary polar extension.
		The self-adjoint element $h=\operatorname{Re}\bE_{W^*(|x|)}(u_x^*w)$ satisfies $\|h\|_\infty\le1$ and $\|h\|_2^2\le \|w\|_2^2$. 
		Since $h$ commutes with $|x|$,
		\[
		2\operatorname{Re}\tau(w^*x)-\|w\|_2^2/5
		\le\tau(2|x|h-h^2/5)
		\le\tau(\psi(|x|))
		\le 2\int_0^2 f(\|x\|_2^2r^2)\,d\nu(r)-1/5.
		\]
		The final step uses \eqref{eq:indefinite-log-f}. 
		Approximate $P_a(w)$ by such sums in $L^2(\cM)$. 
		The integral converges by continuity of $f$ on the bounded range of its arguments, and $\tau(w^*P_a(w))=\beta^2$. 
		Substituting in \eqref{eq:indefinite-wedge-midpoint} bounds the pairing by
		\begin{equation}\label{eq:free-indefinite-scalar}
			\frac{10}{9}\left(2\int_0^2 f(\beta^2r^2)\,d\nu(r)-\beta^2\right)-\frac19.
		\end{equation}
		By Theorem~\ref{thm:free-product-positive-sharp}, $\displaystyle \beta\le \frac{8}{3\pi}$ because $\|w\|_\infty\le1$.
		
		Here are explicit rational bounds for \eqref{eq:free-indefinite-scalar}. 
		Put $\displaystyle \alpha_0=400/471>\frac{8}{3\pi}$, using $\pi>157/50$. 
		If $0\le \beta \le2/5$, use $\Delta\le1/10$ and the monotonicity of $2\alpha_0\beta-\beta^2$ there to obtain the bound $14587/21195$. 
		If $2/5\le \beta \le \alpha_0$, the density of $\nu$ is at most $2/\pi$. 
		The integral of $\Delta$ calculated above bounds \eqref{eq:free-indefinite-scalar} by
		\[
		F(\beta)=\frac{10}{9}\left(2\alpha_0\beta-\beta^2+\frac{Q}{\beta}\right)-\frac19,
		\qquad Q=\frac{286}{23625}.
		\]
		Indeed, $4\alpha^2(1/\sqrt3-1/2)/\pi<Q$ follows from $\pi<22/7$, $\log3<11/10$, and $1/\sqrt3-1/2<78/1000$. 
		The logarithm bound follows by keeping the terms through degree five in the exponential series for $e^{11/10}$, whose sum is $36015101/12000000>3$.
		For the square-root bound, $(289/500)^2>1/3$. 
		Because $F''(\beta)=(10/9)(-2+2Q/\beta^3)<0$ on the interval, its tangent at $\beta_0=21/25$ is an upper bound. 
		Its values at $2/5$ and $\alpha_0$ are, respectively,
		\[
		\frac{10386295181}{14721517125},\qquad
		\frac{4896641957771}{6933834565875}.
		\]
		Both are below $3531/5000$, as is $14587/21195$. 
		These uniform strict bounds and \eqref{eq:1norm} prove the estimate for $\|\Phi\|_{\infty\to1}$. 
		Since $P_a\ne0$ is an orthogonal projection, $\Phi$ has norm one.
		Finally, $(5000)^2>2(3531)^2$ and the universal upper bound is Theorem~\ref{gb:theorem}.
	\end{proof}
	
	\subsection{Sharper Grothendieck constant for finite-index inclusions}
	
	For $\kappa\ge1$ put
	\begin{equation}\label{eq:bcap}
		\mathfrak b_\kappa=\frac2\pi+\frac{1}{2\kappa}
		\int_0^2\left(\frac{s}{2}-\kappa\sqrt{1-\frac{s^2}{4}}\right)_+^2\,d\nu(s).
	\end{equation}
	
	\begin{theorem}\label{thm:cap}
		Suppose that $\cN\subset\cM$ is an irreducible inclusion of finite factors and $0\ne\Phi$ is a bimodule operator with
		$\|\Phi\|_{\infty\to\infty}\le\kappa\|\Phi\|_{2\to2}$ for some $\kappa\ge1$. Then
		\begin{equation}\label{eq:closed}
			\frac{\|\Phi\|_{\infty\to1}}{\|\Phi\|_{2\to2}}
			\ge \mathfrak b_\kappa
			\ge \frac2\pi+\frac{1}{15\pi\kappa(1+\kappa^2)^{3/2}} .
		\end{equation}
		At finite index $\delta^2$ one can take $\kappa=\delta$ for every bimodule operator, and thus
		\begin{equation}\label{eq:closedfinite}
			\frac{\|\Phi\|_{\infty\to1}}{\|\Phi\|_{2\to2}}\ge \mathfrak b_\delta
			\ge\frac2\pi+\frac{1}{15\pi\delta(1+\delta^2)^{3/2}}.
		\end{equation}
	\end{theorem}
	
	\begin{proof}
		Normalize $\|\Phi\|_{2\to2}=1$. 
		When the norm is attained on $\bC1$, the ratio is one. 
		This suffices since $\mathfrak b_\kappa\le 2/\pi+1/8<1$, using $\displaystyle \int_0^2(s/2)^2\,d\nu(s)=1/4$. 
		Otherwise use the circular singular pair $X=U|X|$, $Y=V|Y|$ in the proof of Theorem~\ref{gb:theorem}, with residuals at most $\varepsilon$. Write $h_X=f(|X|)$, $k_X=g(|X|)$, and analogously $h_Y,k_Y$. 
		With the same complementary unitaries as in that proof, set
		\[
		Z_+=-iV_-^*\Phi^\omega(U_+),\qquad
		Z_-=iV_+^*\Phi^\omega(U_-),\qquad
		d_\pm=\|Z_\pm\|_1-\operatorname{Re}\tau(Z_\pm).
		\]
		The mixed terms give
		$\operatorname{Re}\tau(Z_++Z_-)=4/\pi+O(\varepsilon)$.
		The operator-norm cap passes to the ultrapowers, so
		$\|Z_\pm\|_\infty\le\kappa$ and $\|Z_\pm\|_1\le\|\Phi\|_{\infty\to1}$.
		If $Z=W|Z|$ with $W$ a unitary polar extension, then
		\[
		\|Z-|Z|\|_2^2
		=\tau\bigl((W-1)^*(W-1)|Z|^2\bigr)
		\le 2\kappa\bigl(\|Z\|_1-\operatorname{Re}\tau(Z)\bigr).
		\]
		Here $|Z|^2\le\kappa|Z|$ and $(W-1)^*(W-1)=2-W-W^*$.
		
		Averaging the two input unitaries and substituting the definitions of $Z_\pm$ gives
		\[
		V^*\Phi^\omega(U h_X)
		=\tfrac12\bigl(k_Y(Z_++Z_-)+i h_Y(Z_+-Z_-)\bigr).
		\]
		Replacing $Z_\pm$ by $|Z_\pm|$ changes the right side by at most
		$\displaystyle \frac12\sum_\pm\sqrt{2\kappa d_\pm}\le\sqrt{\kappa(d_++d_-)}$ in $L^2(\cM)$, since the factors $k_Y\pm ih_Y$ are unitaries. 
		Apply the conditional expectation onto $W^*(h_Y)$ and take the real part.
		The expression with $|Z_\pm|$ then lies between zero and $\kappa k_Y$, while the original left side differs from $h_Y$ by at most $\varepsilon/2$ in $L^2(\cM)$. 
		In this abelian algebra the positive part is distance to the functions bounded above by $\kappa k_Y$, so
		\[
		\|(h_Y-\kappa k_Y)_+\|_2 \le \varepsilon/2+\sqrt{\kappa(d_++d_-)}.
		\]
		Using $d_++d_-\le2\|\Phi\|_{\infty\to1}-4/\pi+O(\varepsilon)$ and taking $\varepsilon\downarrow0$ proves the lower bound $\mathfrak b_\kappa$.
		
		To estimate the integral in \eqref{eq:bcap}, substitute $t=\sqrt{1-s^2/4}$ and set $t_0=(1+\kappa^2)^{-1/2}$. 
		The positive part is supported on $0\le t\le t_0$, and its integral is
		\[
		\frac4\pi\int_0^{t_0}
		\frac{t^2}{\sqrt{1-t^2}}
		\bigl(\sqrt{1-t^2}-\kappa t\bigr)^2\,dt
		\ge \frac{4t_0^3}{\pi}\int_0^1 v^2(1-v)^2\,dv
		=\frac{2t_0^3}{15\pi}.
		\]
		For the inequality use $\sqrt{1-t^2}\le1$ in the denominator and the chord bound $\sqrt{1-t^2}-\kappa t\ge1-t/t_0$ from concavity.
		This proves \eqref{eq:closed}. 
		The finite-index operator-norm estimate in the proof of Theorem~\ref{thm:cb-algebra} gives \eqref{eq:closedfinite}.
	\end{proof}
	
	\begin{corollary}
		Each irreducible finite-index inclusion $\cN\subset\cM$ satisfies
		\[
		\KNM{\cN}{\cM}
		\le \mathfrak b_\delta^{-1}
		\le\frac{\pi}{2+\bigl(15\delta(1+\delta^2)^{3/2}\bigr)^{-1}}
		<\frac\pi2 .
		\]
	\end{corollary}
	\begin{proof}
		For the fixed inclusion, Theorem~\ref{thm:cap} bounds every nonzero map using the same index parameter $\delta$. Take the supremum over maps and apply \eqref{eq:closedfinite}.
	\end{proof}

	\subsection{\texorpdfstring{$L^p$}{Lp} Grothendieck inequalities}
	
	In this section, we obtain the $L^p$-version of Grothendieck inequalities.
	
	\begin{corollary}\label{cor:pgrothen}
		Suppose $\cN\subset\cM$ is an irreducible inclusion of finite factors and $\Phi$ is a bimodule operator. For $1\le q\le2\le p\le\infty$,
		\[
		\|\Phi\|_{\infty\to1}\le\|\Phi\|_{p\to q}
		\le\|\Phi\|_{2\to2}\le\frac\pi2\|\Phi\|_{\infty\to1}.
		\]
		For $p<q$ an $L^p$ to $L^q$ bound can fail even for identity.
	\end{corollary}
	\begin{proof}
		The upper estimate follows from $\|\Phi(x)\|_q\le\|\Phi(x)\|_2$ and $\|x\|_2\le\|x\|_p$. 
		For the lower estimate use $\|\Phi(x)\|_1\le\|\Phi(x)\|_q$ and $\|x\|_p\le\|x\|_\infty$.
		If $p<q$, projections $e$ of arbitrarily small positive trace in a $\mathrm{II}_1$ factor have $\|e\|_q/\|e\|_p=\tau(e)^{1/q-1/p}\to\infty$. 
		The final inequality is Theorem~\ref{gb:theorem}.
	\end{proof}
	
	\begin{proposition}
		Suppose $\cN\subset\cM$ is an irreducible inclusion of finite factors and $\Phi$ is a bimodule operator.
		If $2\le p\le\infty$ with conjugate exponent $p'$, then
		\[
		\|\Phi\|_{\infty\to1}\le\|\Phi\|_{p\to p'}
		\le\|\Phi\|_{\infty\to1}^{\,1-2/p}\|\Phi\|_{2\to2}^{\,2/p}
		\le \left(\frac\pi2\right)^{2/p}\|\Phi\|_{\infty\to1}.
		\]
	\end{proposition}
	\begin{proof}
		The lower bound follows from the preceding corollary \ref{cor:pgrothen}. 
		In the tracial interpolation scale \cite{KosakiLp}, interpolate the operator bounds $L^\infty\to L^1$ and $L^2\to L^2$ at parameter $\vartheta=2/p$.
		The interpolated inverse exponents are $\vartheta/2=1/p$ in the domain and $1-\vartheta/2=1/p'$ in the range. 
		Theorem~\ref{gb:theorem} bounds the resulting $L^2$ norm factor.
	\end{proof}
	
	\begin{proposition}
		Suppose $\cN\subset\cM$ is an irreducible inclusion of finite factors of finite index $\delta^2$ and $\Phi$ is a bimodule operator. 
		For $2\le p\le\infty$ and conjugate exponent $p'$,
		\[
		\|\Phi\|_{p\to p'}\le \mathfrak b_\delta^{-2/p}\|\Phi\|_{\infty\to1},
		\]
		where $\mathfrak b_\delta$ is given by \eqref{eq:bcap}.
	\end{proposition}
	\begin{proof}
		Use $\|\Phi\|_{2\to2}\le \mathfrak b_\delta^{-1}\|\Phi\|_{\infty\to1}$ from \eqref{eq:closedfinite} in the preceding interpolation bound.
	\end{proof}
	
	\begin{proposition}
		Suppose that $\cN\subset\cM$ is an irreducible inclusion of finite factors and that $\Phi$ is an $\cN$-bimodule map which is positive as an operator on $L^2(\cM)$.
		Let $2<p<\infty$, let $q=p'$ and write $\beta=64/(9\pi^2)$. 
		For any $\beta<\theta<1$, set
		\[
		q_\theta=2-q+(q-1)\theta-\theta^{q-1}.
		\]
		Then
		\[
		\|\Phi\|_{p\to q}
		\le \beta^{-2/p}
		\left[1-q_\theta\left(1-\sqrt{1-(1-\beta/\theta)^4}\right)\right]^{2/q}
		\|\Phi\|_{\infty\to1}.
		\]
	\end{proposition}
	
	\begin{proof}
		For $\Phi\ne0$ we may normalize $\|\Phi\|_{2\to2}=1$. Put
		$\Psi=\Phi^{1/2}$ on $L^2(\cM)$. The left and right $\cN$-actions commute with $\Phi$ and hence with $\Psi$.
		For $r=p$ or $r=\infty$, positivity and the Cauchy--Schwarz inequality give
		\[
		\|\Phi\|_{r\to r'}
		=\sup_{\|v\|_r,\|w\|_r\le1}
		|\langle \Psi(v),\Psi(w)\rangle|
		=\|\Psi\|_{r\to2}^2
		=\|\Psi\|_{2\to r'}^2 .
		\]
		The equality with the squared norm follows by taking $v=w$, and the last step is trace duality and self-adjointness of $\Psi$.
		% Write $s=\|\Psi\|_{2\to1}$. 
		Theorem~\ref{gb:theorem} and $\|\Psi\|_{2\to2}=1$ imply $\sqrt \beta\le \|\Psi\|_{2\to1}\le1$.
		
		Fix a unit vector $x\in L^2(\cM)$ with $A=|\Psi(x)|$ and $\gamma_\Psi=\tau(A)>0$. 
		Write $\Psi(x)=wA$ with a unitary polar extension $w$. 
		Set
		\[
		e=1_{[\theta/\gamma_\Psi,\infty)}(A),\qquad
		e_w^\perp=1-wew^*\qquad e^\perp =1-e
		\]
		Since $\gamma_\Psi\ge(\theta/\gamma_\Psi)\tau(e)$, we have $\tau(e^\perp)\ge(1-\gamma_\Psi^2/\theta)_+$.
		We claim that
		\begin{equation}\label{eq:positive-lp-two-part}
			\tau(Ae)^2+\tau(e^\perp)^4\gamma_\Psi^2\le \|\Psi\|_{2\to1}^2.
		\end{equation}
		For unitaries $u,v\in\cN$, the contraction $e_w^\perp uwv e^\perp $ is a trace-duality test for $e_w^\perp u\Psi(x)v e^\perp $. 
		Their pairing is
		\[
		\tau\bigl((e_w^\perp uwv e^\perp )^*e_w^\perp uwA v e^\perp \bigr)
		=\tau\bigl(w^*u^*e_w^\perp u w A v e^\perp v^*\bigr).
		\]
		By irreducibility, the scalar $\tau(e^\perp)1$ is in the $L^2$-closed convex hull of the $\cN$-unitary conjugates of each of $e_w^\perp $, $e^\perp$.
		Indeed, the unique vector of minimum norm in either closed convex hull is conjugation invariant, and conditional expectation onto $\cN'\cap\cM=\mathbb C1$ identifies it as $\tau(e^\perp) 1$. 
		Average the last pairing over $u$ and $v$ independently.
		Both averaged factors are contractions. 
		Their $L^2$ distances from $\tau(e^\perp)1$ tend to zero, so Cauchy--Schwarz inequality, with $A\in L^2(\cM)$, shows the averaged pairing tends to $\tau(e^\perp)^2\gamma_\Psi$. 
		Thus
		\[
		\sup_{u,v\in\U(\cN)}\|e_w^\perp u\Psi(x)v e^\perp \|_1\ge \tau(e^\perp)^2\gamma_\Psi .
		\]
		For each $u,v$ choose a phase $e^{i\phi}$ so $x_1=e^{i\phi} uxv$ satisfies $\operatorname{Re}\langle x,x_1\rangle=0$. 
		Let $v_0=we$ and let $v_1$ be the polar partial isometry of $e_w^\perp \Psi(x_1)e^\perp $.
		The left and right supports of $v_0,v_1$ are disjoint, so $v_0\pm v_1$ are contractions. 
		If $\gamma_{\perp}=\|e_w^\perp \Psi(x_1)e^\perp \|_1$, trace-duality testing of $\Psi(x\pm\varepsilon x_1)$ against $v_0\pm v_1$ and averaging the two real parts gives
		\[
		\tau(Ae)+\varepsilon \gamma_{\perp}\le \|\Psi\|_{2\to1}\sqrt{1+\varepsilon^2}\qquad(\varepsilon\ge0).
		\]
		Taking the supremum over $\varepsilon$ gives $\tau(Ae)^2+\gamma_{\perp}^2\le \|\Psi\|_{2\to1}^2$, including the cases $\tau(Ae)=0$ or $\gamma_{\perp}=0$ by limits. 
		Taking the supremum over $u,v$ proves \eqref{eq:positive-lp-two-part}.
		
		Let $\alpha=q-1\in(0,1)$.
		On the probability space with measure $A\, d\tau/\gamma_\Psi$ in the spectral calculus of $A$, the variable $z=\gamma_\Psi A$ has mean $\tau(A^2)\le1$. 
		The tangent deficit $1-\alpha+\alpha z-z^\alpha$ is nonnegative on $[0,\infty)$ and at least $q_\theta$ on $[0,\theta]$. 
		Since the probability of $z\ge\theta$ is $\tau(Ae)/\gamma_\Psi$, it follows that
		\[
		\tau(A^q)\le \gamma_\Psi^{2-q}\bigl(1-q_\theta+q_\theta\,\tau(Ae)/\gamma_\Psi\bigr).
		\]
		With $t=\gamma_\Psi^2$ and $h=\|\Psi\|_{2\to1}^2$, define $f(t)=(1-t/\theta)_+$ and
		\[
		Z(t,h)=\min\left\{1,\sqrt{h/t-f(t)^4}\right\}.
		\]
		Then \eqref{eq:positive-lp-two-part} gives
		\begin{equation}\label{eq:positive-lp-scalar}
			\frac{\tau(A^q)}{\|\Psi\|_{2\to1}^q}
			\le t^{(2-q)/2}h^{-q/2}\bigl(1-q_\theta+q_\theta Z(t,h)\bigr),
			\qquad 0<t\le h,\quad \beta\le h\le1.
		\end{equation}
		
		For clarity, we maximize the scalar expression explicitly.
		Taylor's formula gives
		\[
		0<q_\theta
		=\alpha(1-\alpha)\int_\theta^1 (r-\theta)r^{\alpha-2}\,dr
		\le\min(\alpha,1-\alpha)\frac{(1-\theta)^2}{2\theta^2}
		<\frac{\min(\alpha,1-\alpha)}8,
		\]
		since $\theta>\beta>2/3$.
		Fix $h$.
		The expression in \eqref{eq:positive-lp-scalar} is increasing in $t$ wherever $Z(t,h)=1$. 
		On a region with $Z(t,h)<1$, one has $h/t<2$, so $t>\beta/2>1/3$, $f(t)<2/3$ and $Z(t,h)>8/9$. 
		Differentiation there gives
		\[
		t\partial_t Z
		=\frac{-h/t+4(t/\theta)f(t)^3}{2Z(t,h)}
		\ge-\frac1{Z(t,h)}.
		\]
		The logarithmic derivative with respect to $t$, multiplied by $t$, is therefore at least
		\[
		\frac{1-\alpha}{2}-\frac{9q_\theta}{8(1-q_\theta)}
		>\left(\frac12-\frac9{56}\right)(1-\alpha)>0 .
		\]
		The expression is continuous at the boundary between the two regions, hence it is maximized at $t=h$.
		
		On this diagonal the expression is $h^{-\alpha}(1-q_\theta+q_\theta Z(h,h))$. 
		For $h\ge\theta$ it is decreasing, because $Z(h,h)=1$. 
		For $\beta\le h<\theta$, one has $f(h)<1/3$, $Z(h,h)>8/9$, and
		\[
		h\frac{d}{dh} Z(h,h)=\frac{2(h/\theta)f(h)^3}{Z(h,h)}<\frac1{12}.
		\]
		Consequently its logarithmic derivative, multiplied by $h$, is at most $-\alpha+q_\theta/(12(1-q_\theta))<0$. 
		The maximum in \eqref{eq:positive-lp-scalar} is attained at $t=h=\beta$, and equals
		\[
		\beta^{1-q}\left[1-q_\theta\left(1-\sqrt{1-(1-\beta/\theta)^4}\right)\right].
		\]
		The case $\gamma_\Psi=0$ gives zero on the left of \eqref{eq:positive-lp-scalar}. 
		Take the supremum over $x$ and raise to the power $2/q$. 
		Using $2(1-q)/q=-2/p$ and the norm identities at the beginning of the proof completes the argument.
	\end{proof}
	
	\subsection{Sharp Grothendieck constant for positive bimodule operators and a lower bound for finite-index inclusions}
	
	In this section, we shall prove the following theorem.
	\begin{theorem}\label{thm:fiber-main}
		There is a sequence of irreducible finite-index inclusions of hyperfinite
		\(\mathrm{II}_1\) factors with orthogonal bimodule projections whose
		Grothendieck ratios converge to \(9\pi^2/64\). Indefinite self-adjoint
		bimodule maps on the same sequence give
		\[
		\KIf>\frac{500000000}{350700187}=1.425719228\ldots,
		\qquad \KIf_+=\frac{9\pi^2}{64}.
		\]
	\end{theorem}
	
	We shall apply the matrix model and the joint moment estimate in the following to achieve the theorem.
	The construction passes trace estimates to finite index uniformly over
	all contractions in the upper factors.

	\subsubsection{Construction of the finite model}
	\begin{proposition}\label{prop:fiber-semicircle}
		There is a sequence of finite-dimensional matrix algebras with normalized trace, trace-orthonormal selfadjoint tuples $(E_1,\ldots,E_m)$, and projections $P$ onto their linear spans with the following properties.
		\begin{enumerate}
			\item The amplification of $P$ restricts to a nonzero operator two-box of an irreducible finite-index $\mathrm{II}_1$ inclusion, and the restricted projection kills $1$.
			\item Let $(B,\tau_B)$ be any finite tracial algebra with normalized trace $\tau_B$, let $\displaystyle \sum_i a_i^*a_i\le1$, $\displaystyle \sum_i a_i a_i^*\le1$, and let $(s_i)$ be free variance-one semicirculars with their tracial state. 
			Set $\displaystyle X=\sum_i a_i\otimes E_i$ and $\displaystyle Y=\sum_i a_i\otimes s_i$. 
			For every fixed continuous $\varphi:[0,\infty)\to\bR$ of at most polynomial growth,
			\[
			\sup_{B,(a_i)}\left|\tau\big(\varphi(|X|)\big)-
			\tau\big(\varphi(|Y|)\big)\right|\longrightarrow0 .
			\]
			Each $\tau$ uses the appropriate product trace.
		\end{enumerate}
	\end{proposition}
	\begin{proof}
		First fix a maximum word length $l$ and a tolerance $\varepsilon>0$. 
		Choose a finite group $K$ with an irreducible unitary representation $\rho$ of dimension $d$, with $d$ as large as needed. 
		The standard representation of the permutation group $S_{d+1}$ suffices. 
		Let the colors be $\{1,\ldots,m\}=K\times\{1,\ldots,L\}$ with $m=|K|L$.
		Left multiplication in the first coordinate embeds $K$ in $S_m$. 
		For Hermitian GUE seed matrices $D_j\in M_d(\bC)$ with $j=1, \ldots, L$, assign to color $(k,j)$ the matrix
		\[
		B_{(k,j)}=\rho(k)D_j\rho(k)^*\in M_{d}(\bC).
		\]
		We choose the seeds so that, for every list of colors $i_1,\ldots,i_q$, $q\le l$, the quantity
		\begin{align}\label{eq:permuave}
			\frac{1}{m!}\sum_{\sigma\in S_m}\tr_d
			(B_{\sigma(i_1)}\cdots B_{\sigma(i_q)})
		\end{align}
		differs by at most $\varepsilon$ from $\tau(s_{i_1}\cdots s_{i_q})$.
		The choices of $d,L$ need not be effective.
		To see existence, take independent variance-normalized GUE seeds.
		In dimension $d$, for independent GUE matrices $G_i$ with $\bE\tr_d(G_i^2)=1$, Gaussian Wick expansion of the expectation of each normalized trace gives one term $d^{-2g}$ for each color-compatible pairing, with $g$ the genus of the oriented polygon gluing.
		The terms with $g=0$ are the noncrossing pairings, giving precisely the semicircular moment. 
		One can equivalently use the asymptotic freeness theorem for independent GUE matrices \cite[Section~1.10]{MS}. 
		There are only finitely many equality patterns at these lengths, so first choose $d$ for simultaneous approximation. 
		At this fixed $d$, truncate the Gaussian eigenvalues at a sufficiently large common cutoff so as to preserve all these moments within a second tolerance. 
		The truncated law is still invariant under unitary conjugation. 
		Sample $L$ independent seeds from that law.
		
		When distinct color positions are sent by a uniform permutation to distinct $K$-orbits, their corresponding seed matrices are independent, and conjugations by $\rho(k)$ preserve their laws. 
		For a word with $C\geq 2$ distinct colors, the probability that two of their images are in the same orbit is at most $\binom C2 (|K|-1)/(m-1)$. 
		Consequently, the expectation of the permutation average approaches the desired truncated-Gaussian moment as $L\to\infty$. 
		For the variance, expand the covariance term by term using two independent permutations. 
		The two lists share a seed with probability $O_l(1/L)$. 
		Only such terms have a nonzero covariance, and these terms are uniformly bounded at the fixed cutoff.
		The variance therefore tends to zero. 
		A single realization works simultaneously for all the patterns. 
		Taking the intermediate tolerances smaller yields the assertion.
		
		Choose a finite $S_m$-equivariant quotient $F$ of the free product of the $m$ involutions with no nontrivial reduced relation of length $\le l$, as follows.
		On the ball of radius \(l\) in the free product, let the \(i\)-th
		involution interchange \(w\) with \(g_iw\) if both are in the ball and fix
		unpaired vertices.
		A reduced word of length at most \(l\) moves the identity to the corresponding vertex. 
		The ball and the generating permutations are equivariant under permutations of the colors, so the finite permutation image has the asserted properties. 
		Write $g_i$ for its generators, $\Gamma=F\rtimes S_m$, and $\pi:\Gamma\to S_m$ for the quotient.
		The vertices are all $t\in\Gamma$, with a copy of $\bC^d$ at each vertex. 
		Define $X_i$ on $\ell^2(\Gamma)\otimes \bC^d$ as follows:
		\[
		X_i(\delta_t\otimes v)=\delta_{g_i t}\otimes B_{\pi(t)^{-1}(i)}v,
		\]
		where $v\in \bC^d$ and $\delta_t$ is the Kronecker symbol.
		These are selfadjoint because $\pi(g_i t)=\pi(t)$.
		The trace of a word $X_{i_1}\cdots X_{i_q}$ of length at most $l$ is zero unless its colors freely reduce to the identity. 
		Otherwise it is exactly the permutation average above \eqref{eq:permuave}.
		Every word having a nonzero free semicircular moment admits a color-compatible noncrossing pairing and so freely reduces to the identity. 
		Thus the traces approximate the desired free moments on all patterns. 
		Distinct $X_i$ are exactly trace-orthogonal and have a common variance $\mathfrak v=\|X_i\|_2^2$ as close to one as desired. 
		Set $E_i=X_i/\sqrt {\mathfrak v}$. 
		By decreasing the original tolerances, the normalized moments can still approximate the semicircular ones to any prescribed accuracy.
		
		We now realize the matrix projection in an irreducible inclusion.
		Put $H=S_m\times(F\rtimes K)$, where the second factor is a subgroup of $\Gamma$.
		For $h=(\sigma,r)$ let $h(t)=\sigma t r^{-1}$ on vertices and let it act by $\rho(\pi(r))$ on the fiber. 
		For any $t\in\Gamma$, choose $\sigma=\pi(t)$ and $r=t^{-1}\sigma\in F$. 
		Then $h(1)=t$, proving transitivity. 
		The stabilizer of the identity consists of $(k,k)$, $k\in K$.
		Its action on the fiber is irreducible. 
		The group acts on grades $\Gamma$ by $\theta_h=\operatorname{Ad}\sigma$. 
		Indeed
		\[
		h(t)h(u)^{-1}=\theta_h(tu^{-1}),\qquad
		\pi(h(t))^{-1}\sigma(i)=\pi(r)\pi(t)^{-1}(i).
		\]
		The second identity and \(B_{ki}=\rho(k)B_i\rho(k)^*\) show that the
		combined action sends \(X_i\) to \(X_{\sigma(i)}\).

		Let \(\Omega=\Gamma\rtimes_\theta H\) act on the hyperfinite factor \(\displaystyle \cQ_0=\overline{\bigotimes}_{n\ge1}B(\ell^2(\Omega))\) by conjugation by the left regular representation in every tensor factor. 
		This action is outer.
		A fixed unitary in \(\cQ_0\) asymptotically commutes in \(L^2(\cQ_0)\) with rank-one diagonal projections placed in distant tensor legs, by approximation by local tensors. 
		A nonidentity group element sends each of those projections onto an orthogonal projection, with \(L^2\)-distance \((2/|\Omega|)^{1/2}\). 
		Denote the crossed-product unitaries by \(u_g\). 
		In $M_{|\Gamma|}(M_d(\cQ_0\rtimes\Gamma))$ take the graded algebra $\mathcal B$ having $(t,u)$ entry in $M_d(\cQ_0 u_{tu^{-1}})$.
		Conjugation by $\operatorname{diag}(u_t)_{t\in\Gamma}$ identifies this algebra with $M_{|\Gamma|d}(\cQ_0)$. 
		The upper algebra is $\cM=\mathcal B^H$, where $h\in H$ acts on coefficients by $a u_g\mapsto\alpha_h(a)u_{\theta_h(g)}$ along with the indicated actions on vertices and fibers. 
		The lower factor $\cN$ is the constant diagonal copy of $\cQ_0^H$. 
		If \(a_g\) is a Fourier coefficient of an element of
		\((\cQ_0^H)'\cap(\cQ_0\rtimes\Gamma)\), the intertwining identity is
		\(b a_g=a_g\alpha_g(b)\) for \(b\in \cQ_0^H\). 
		For \(g\ne1\) in \(\Gamma\), let \(p_n\) be the projection onto \(\ell^2(H)\subset\ell^2(\Omega)\) in the \(n\)-th tensor leg. 
		It belongs to \(\cQ_0^H\), is asymptotically central, and is orthogonal to \(\alpha_g(p_n)\). 
		The intertwining relation implies \(\|a_g(p_n-\alpha_g(p_n))\|_2\to0\). 
		By local tensor approximation, the square of this norm tends to \(2|H|\|a_g\|_2^2/|\Omega|\), so \(a_g=0\). 
		The identity coefficient is scalar, by irreducibility of the outer fixed-point inclusion \(\cQ_0^H\subset \cQ_0\).
		Consequently a commuting matrix in \(\mathcal B\) is diagonal in the vertex coordinate with blocks in \(M_d(\bC)\). 
		For an $H$-invariant such matrix the value at the identity commutes with $\rho(K)$ and is therefore scalar. 
		Transitivity makes the relative commutant trivial. 
		In particular, $\cM$ is a factor and the inclusion is irreducible.
		The full graded algebra has index $(|\Gamma|d)^2$ over the constant copy of $\cQ_0$, as its left $\cQ_0$ module consists of $(|\Gamma|d)^2$ standard summands.
		As $[\cQ_0:\cQ_0^H]=|H|$, the intermediate inclusion has finite index. 
		The operator $\id_{\cQ_0\rtimes\Gamma}\otimes P$ preserves the grading because every nonzero entry of a given $X_i$ has grade $g_i$, and it is $H$-equivariant and bimodular over $\cQ_0^H$. 
		Its fixed range contains the nonzero element $\displaystyle \sum_i u_{g_i}\otimes E_i$. 
		It annihilates the identity.
		
		Finally we prove uniformity in the coefficients. 
		For slots $j=1,\ldots,q$ write $a_i^{(j)}$ for $a_i$ or $a_i^*$ as needed in a product. 
		If an equality pattern $\mathcal P$ freely reduces, define
		\[
		C_{\mathcal P}=\sum_{\ker(i_1,\ldots,i_q)=\mathcal P}
		\tau_B(a_{i_1}^{(1)}\cdots a_{i_q}^{(q)}).
		\]
		Its absolute value is bounded by a constant depending only on $q$. 
		Insert group matrices $\lambda(g_i)$ of the \emph{free product} of involutions into each slot, using its group trace.
		On the specified pattern the product of the group matrices is the identity.
		Use inclusion--exclusion for the inequalities between distinct blocks of $\mathcal P$. 
		To impose a required equality $i_j=i_k$, multiply the two slots by $z_{i_j}$ and $\overline z_{i_k}$ for independent uniform unit phases $(z_i)_i$ and average. 
		Take independent phase families for all such equalities. 
		For each choice of phases the resulting unrestricted sum is the product trace of $q$ factors of the form $\displaystyle \sum_i z_i a_i^{(j)}\otimes\lambda(g_i)$, each of norm at most two. 
		In fact $\lambda(g_i)=C_i+C_i^*$ where the creation part $C_i$ prefixes $i$ to reduced words not starting in $i$ and is zero on the rest. 
		The $C_i$ have orthogonal ranges, so the two sums have norms bounded by the column and row norms, respectively. 
		Thus $|C_{\mathcal P}|\le2^q 2^{\binom q2}$. 
		Patterns which do not freely reduce contribute zero on both sides. 
		Since all trace errors are constant on patterns, the even moment differences are bounded by their maximum error times a constant depending only on the moment order. 
		Selecting $l\to\infty$ and the error successively small gives the claimed uniform convergence of each moment.
		
		In the free Fock representation $s_i=\ell_i+\ell_i^*$ and the ranges of the $\ell_i$ are orthogonal, so $\|Y\|\le2$.
		If the test-function claim failed, take a sequence of violating coefficients. The laws of $|Y|^2$ admit a weakly convergent subsequence supported on $[0,4]$.
		On that subsequence the laws of $|X|^2$ are tight, have uniformly bounded moments of each fixed order, and have the same limiting moments. 
		Their only weak limit is the same compactly supported measure. 
		Uniform integrability gives convergence also for the given continuous test function of polynomial growth. 
		This proves the claim.
	\end{proof}

	\begin{corollary}\label{cor:fiber-finite-depth}
		The factors in Proposition~\ref{prop:fiber-semicircle} may be taken hyperfinite. 
		Their index is $(d|\Gamma|)^2$, and the inclusions have finite depth.
	\end{corollary}
	\begin{proof}
		Take the hyperfinite factor $\cQ_0$. 
		Under the identification of $\mathcal B$ with $M_{d|\Gamma|}(\cQ_0)$ that sends a $(t,u)$ entry $b$ to $u_t^* b u_u$, the action of $h=(\sigma,r)$ is an inner matrix conjugation followed by an entrywise automorphism of $\cQ_0$.
		The latter is $\alpha_{(r\sigma^{-1},h)}$, since $u_{h(t)}^*u_{\theta_h(t)}=u_{r\sigma^{-1}}$ for all $t$. 
		For $h\ne e$ this is an outer automorphism since $(r\sigma^{-1},h)\ne e$ in $\Omega$. 
		Its matrix amplification remains outer.
		Thus $[\mathcal B:\mathcal B^H]=|H|$ and dividing $[\mathcal B:\cQ_0^H]=(d|\Gamma|)^2|H|$ gives the asserted index. 
		Finite-group fixed points of the hyperfinite factor are hyperfinite.
		
		Write $U_g$ for the $\cQ_0$--$\cQ_0$ correspondence of $\alpha_g$, $g\in\Omega$, with the convention $U_g\otimes_{\cQ_0} U_k=U_{gk}$. 
		Let $V_g$ denote its restriction to $\cQ_0^H$ on both sides. 
		The group fixed-point identity $\displaystyle L^2(\cQ_0)\otimes_{\cQ_0^H}L^2(\cQ_0)\cong\bigoplus_{h\in H} U_h$ gives
		\[
		V_g\otimes_{\cQ_0^H}V_k\cong\bigoplus_{h\in H}V_{ghk}.
		\]
		The $(t,u)$ matrix entries show that $L^2(\mathcal B)$, as an $\cQ_0$--$\cQ_0$ correspondence, is a finite direct sum of the $U_g$ with $g\in\Gamma$. 
		Since $L^2(\mathcal B^H)$ is an $\cQ_0^H$--$\cQ_0^H$ summand of $L^2(\mathcal B)$, all its tensor powers decompose into simples drawn from the finitely many finite-length correspondences $V_g$, $g\in\Omega$.
		This proves finite depth.
	\end{proof}

	\subsubsection{Transfer to unitary tests}
	Let $P_j$ denote the restricted orthogonal projection in Proposition~\ref{prop:fiber-semicircle}, and write $d\nu(r)=\pi^{-1}\sqrt{4-r^2}\,dr$ on $[0,2]$.
	For a contraction $W$ in the upper factor of the $j$th inclusion put $\displaystyle X=P_jW=\sum_i a_i\otimes E_i$ and $\displaystyle Y=\sum_i a_i\otimes s_i$.
	Here the ancillary algebra is $B=\cQ_0\rtimes\Gamma$.
	Writing $\bE_B$ for the partial trace over the matrix coordinate,
	orthonormality gives
	\[
	\bE_B\big((W-X)^*(W-X)\big)=\bE_B(W^*W)-\sum_i a_i^*a_i.
	\]
	Applied also to $W^*$ this gives both square-sum bounds of the
	proposition. 
	Furthermore $\mathfrak b_X^2:=\tau(|X|^2)=\tau(|Y|^2)
	=\operatorname{Re}\tau(W^*X)$.
	
	\begin{lemma}\label{lem:fiber-res}
		Let $Y$ be bounded, selfadjoint and semicircular over a finite tracial von Neumann algebra $\mathcal A$, with zero mean, covariance $\eta$, and $\mathfrak v=\tau(Y^2)$. 
		Let $Y_0$ be a scalar semicircle of variance $\mathfrak v$.
		For every $s>0$,
		\[
		\tau\left(\frac{Y^2}{s+Y^2}\right)\le\tau\left(\frac{Y_0^2}{s+Y_0^2}\right).
		\]
	\end{lemma}
	\begin{proof}
		We may assume $\mathfrak v>0$. 
		Write $t=\sqrt s$ and
		$g=\bE_{\mathcal A}[t(t^2+Y^2)^{-1}]$. 
		This is positive and invertible.
		The odd operator-valued moments vanish.
		On the imaginary axis the semicircular Dyson equation \cite{HeltonFarSpeicher} therefore reads $g^{-1}=t1+\eta(g)$.
		Let $\sigma$ be the spectral distribution of $g$. 
		The positive commuting operators of left and right multiplication by $g$ on $L^2$ have a joint spectral measure.
		In the vector $g^{1/2} Y g^{1/2}$ this gives a finite positive measure $\rho$ determined by
		\[
		\int f(x)h(y)\,d\rho(x,y) =\tau\bigl(f(g)g\,\eta(g h(g))\bigr).
		\]
		Both marginals are $(1-tx)\,d\sigma(x)$ by the Dyson equation and trace symmetry of $\eta$. 
		Its mass is $M=1-t\tau(g)>0$ and
		$\displaystyle \mathfrak v=\int(xy)^{-1}\,d\rho$. 
		Let $\mu=\rho/M$.
		Then $\mathbb E_\mu x=\mathbb E_\mu y= (\tau(g)-t\tau(g^2))/M\le \tau(g)$. 
		Joint convexity of $\displaystyle (x,y)\mapsto \frac{1}{xy}$ on the positive quadrant gives
		\[
		\mathfrak v\ge \frac{M}{(\mathbb E_\mu x)(\mathbb E_\mu y)}
		\ge \frac{1-t\tau(g)}{\tau(g)^2}.
		\]
		It follows that $\tau(g)\ge q$, where $1=tq+\mathfrak vq^2$.
		The scalar semicircular Dyson equation identifies
		$\displaystyle q=t\tau\left((t^2+Y_0^2)^{-1}\right)$, proving the assertion.
	\end{proof}
	
	For the quarter-circle law $\nu$ above, write $J(z)=1-h(z)^{-1}$ for $h(z)=(1+\sqrt{1+4z})/2$. 
	Then $\displaystyle J(z)=\int zr^2/(1+zr^2)\,d\nu(r)$, by the scalar Dyson equation.
	For the tuple in Proposition~\ref{prop:fiber-semicircle} use the self-adjoint dilation
	\[
	H=\begin{pmatrix}0&Y\\Y^*&0\end{pmatrix}
	=\sum_i \begin{pmatrix}0&a_i\\a_i^*&0\end{pmatrix}\otimes s_i .
	\]
	It is semicircular over $M_2(B)$ by the operator-valued Wick formula, with covariance $\displaystyle v\mapsto\sum_i h_i v h_i$ for the displayed coefficient matrices $h_i$. 
	The distribution of $H^2$ in the normalized trace is that of $|Y|^2$ since $YY^*$ and $Y^*Y$ have identical trace distributions. 
	Its variance parameter is $\mathfrak b_X^2$. 
	Lemma~\ref{lem:fiber-res} gives
	\begin{equation}\label{eq:fiber-comparisons}
		\tau\left(\frac{|Y|^2}{s+|Y|^2}\right)\le J(\mathfrak b_X^2/s),\qquad
		\tau(|Y|)\le \frac{8}{3\pi}\mathfrak b_X .
	\end{equation}
	For the second inequality integrate the first using
	$\displaystyle \sqrt z=\pi^{-1}\int_0^\infty \frac{z}{z+s}\,\frac{ds}{\sqrt s}$ and $ \displaystyle \int r\,d\nu(r)=\frac{8}{3\pi}$. 
	All integrands are nonnegative.
	
	\begin{corollary}\label{cor:finite-positive-transfer}
		The projections \(P_j\) on the sequence of finite-index inclusions satisfy
		\[
		\lim_{j\to\infty}\|P_j\|_{\infty\to1}=\frac{64}{9\pi^2}.
		\]
		Consequently \(\KIf_+=9\pi^2/64\).
	\end{corollary}
	\begin{proof}
		The transfer of the test function $x\mapsto x$ gives $\displaystyle \mathfrak b_X^2\le\tau(|X|)\le \frac{8}{3\pi}\mathfrak b_X+o(1)$, uniformly in $W$, hence $\displaystyle \mathfrak b_X\le \frac{8}{3\pi}+o(1)$.
		Positivity and Cauchy--Schwarz inequality identify $\|P_j\|_{\infty\to1}$ with $\displaystyle \sup_{U\in\U(\cM)}\|P_j U\|_2^2$.
		For the other direction put $\displaystyle Z=m^{-1/2}\sum_i u_{g_i}\otimes E_i$ in the fixed upper factor, as in the realization proof of Proposition \ref{prop:fiber-semicircle}. 
		Then $P_jZ=Z$, $Z=Z^*$ and $\|Z\|_2=1$. 
		On replacing the $E_i$ by $s_i$, one obtains an exact
		scalar semicircular distribution because the operator-valued covariance $v\mapsto m^{-1}\sum_i u_{g_i}v u_{g_i}$ is unital. 
		The noncrossing recursion then gives scalar Catalan even moments and zero odd moments. 
		Transfer implies $\displaystyle \tau(|Z|)\to \frac{8}{3\pi}$.
		Choose $U=\sgn(Z)$, with value $1$ at zero. 
		Then $\|P_jU\|_2\ge |\langle P_jU,Z\rangle|=\tau(|Z|)$.
		Combining the two estimates gives the limit. 
		Theorem~\ref{gb:theorem} gives the matching upper bound on \(\KIf_+\).
	\end{proof}

	\subsubsection{Joint sixth- and eighth-moment estimate}
	\begin{proposition}\label{prop:fiber-eighth-joint}
		For the finite-index sequence of Proposition~\ref{prop:fiber-semicircle},
		let $d=0.1291$ and $T_j=(P_j-d\,\id)/(1-d)$. Then
		\[
		\limsup_j\|T_j\|_{\infty\to1}<0.701400374,\qquad
		\KIf>\frac{10^9}{701400374}>1.42571922.
		\]
	\end{proposition}
	\begin{proof}
		We use the hermitianization $H$ of $\displaystyle Y=\sum_i a_i\otimes s_i$ from the transfer to unitary tests. 
		Write $\eta$ for its covariance and set
		$\gamma_1=\tau (\eta(1))$, $\gamma_2=\tau (\eta(1)^2)=\tau(\eta^2(1))$, $\gamma_3=\tau (\eta(1)^3)$, $\gamma_{1,2}=\tau(\eta(1)\eta^2(1))$ and $\Delta=\gamma_2-\gamma_1^2$.
		The row and column contraction inequalities give $0\le \eta(1)\le1$.
		The free Wick recursion for $M_j=\tau( H^j)$ gives
		\begin{align*}
			M_4=& 2\tau( \eta(1)^2),\qquad 
			M_6=2\tau( \eta(1)^3)+3\tau(\eta(1) \eta^2(1)),\\
			M_8=& 2\tau (\eta(1)^4)+8\tau(\eta(1)^2\eta^2(1))+4\tau (\eta^2(1)^2).
		\end{align*}
		Indeed, if $V_k$ is the conditional moment of $H^{2k}$, then $V_0=1$ and $\displaystyle V_k=\sum_{j=0}^{k-1}\eta(V_j)V_{k-1-j}$.
		In detail, \(V_2=\eta(1)^2+\eta^2(1)\) and \(V_3=\eta(1)V_2+\eta^2(1)\eta(1)+\eta(V_2)\).
		Trace symmetry of \(\eta\) gives
		\[
		\tau(V_3)=2\tau(\eta(1)^3)+3\tau(\eta(1)\eta^2(1)),\qquad
		\tau(V_4)=2\tau(\eta(1)V_3)+2\tau(\eta^2(1) V_2).
		\]
		In the second formula,
		\(\tau(\eta(1)V_3)=\tau(\eta(1)^4)+3\tau(\eta(1)^2\eta^2(1))+\tau(\eta^2(1)^2)\) and
		\(\tau(\eta^2(1)V_2)=\tau(\eta(1)^2\eta^2(1))+\tau(\eta^2(1)^2)\), giving the stated eighth moment.
		
		In the vector $H$, left and right multiplication by $\eta(1)$ have a joint spectral measure of mass $\gamma_1$ with identical marginals $x\,d\sigma_{\eta(1)}(x)$.
		Integrating $xy\le(x^2+y^2)/2$ gives $\gamma_{1,2}\le \gamma_3$. 
		Also $\gamma_3\le \gamma_1^3+(1+2\gamma_1)\Delta$ by integrating 
		\begin{align*}
			x^3-\gamma_1^3-3\gamma_1^2(x-\gamma_1)=(x+2\gamma_1)(x-\gamma_1)^2
		\end{align*}
		and using $0\le \eta(1)\le1$. 
		For $\gamma_1>0$, the Cauchy--Schwarz inequality gives
		\[
		\tau( \eta(1)^4)\ge \gamma_3^2/\gamma_2,\qquad \tau (\eta^2(1)^2)\ge \gamma_{1,2}^2/\gamma_2,\qquad
		\tau(\eta(1)^2\eta^2(1))\ge \gamma_{1,2}^2/\gamma_2.
		\]
		The last inequality uses the positive functional $a\mapsto\tau(\eta^2(1)a)$ and $\tau (\eta^2(1))=\gamma_2$. 
		Thus 
		\begin{align*}
			M_6=& 2\gamma_3+3\gamma_{1,2}, \\
			M_8\ge & (2\gamma_3^2+12\gamma_{1,2}^2)/\gamma_2\ge 4\gamma_1\gamma_3+24\gamma_1\gamma_{1,2}-14\gamma_1^2\gamma_2,
		\end{align*}
		using $(x-\gamma_1\gamma_2)^2\ge0$. 
		Whenever $\alpha_4<0$ and $\alpha_3+8\gamma_1 \alpha_4\ge0$, it follows that
		\[
		\alpha_3 M_6+\alpha_4 M_8\le
		5\alpha_3\gamma_1^3+14\alpha_4\gamma_1^4+\bigl(5\alpha_3(1+2\gamma_1)+\alpha_4(28\gamma_1+42\gamma_1^2)\bigr)\Delta .
		\]
		Indeed the coefficients of $\gamma_3,\gamma_{1,2}$ are $2\alpha_3+4\gamma_1 \alpha_4$ and $3\alpha_3+24\gamma_1 \alpha_4$, respectively. 
		Both are nonnegative. 
		The assertion at $\gamma_1=0$ follows from $\eta(1)=0$. 
		Set
		\[
		\begin{aligned}
			(\alpha_0,\alpha_1,\alpha_2)&=(0.0000002,\ 0.6125154119,\ -0.11979471),\\
			(\alpha_3,\alpha_4,h_0)&=(0.02091142647,\ -0.001649701762,\ 0.44757658117).
		\end{aligned}
		\]
		For $q(z)=\alpha_0+\alpha_1z+\alpha_2 z^2+\alpha_3 z^3+\alpha_4 z^4+h_0z/(z+0.247)$, we have
		$\displaystyle q(x^2)\ge H_d(x):=\sup_{|t|\le1}(xt-dt^2)$ on $0\le x\le2$, where $t$ is real. 
		For the exact certificate form $(x^2+0.247)(q(x^2)-H_d(x))$ on the intervals $[0,2d]$ and $[2d,2]$,
		using $H_d(x)=x^2/(4d)$ and $x-d$, respectively.
		For each of the two polynomials \(P_1,P_2\) so obtained, start with
		\(P_i,P_i'\) and adjoin the negatives of successive Euclidean remainders. 
		With zeros omitted, the Sturm variation counts are
		\[
		\begin{array}{c|cc|cc}
			&\text{left endpoint}&\text{right endpoint}&\text{left variations}&\text{right variations}\\
			\hline
			P_1&0&2d&5&5\\
			P_2&2d&2&5&5
		\end{array}
		\]
		The values \(P_1(0)=247/5000000000\) and
		\(P_2(2)=63585761/62500000000000\) are positive, while
		\(P_1(2d)=P_2(2d)>10^{-6}\) by substitution. 
		The Sturm counts rule out interior roots, proving positivity on each interval.
		
		For unitaries \(U,V\), rotate one by a scalar phase so that the tested pairing is real and nonnegative. 
		Set \(W=(U+V)/2\), \(W_-=(U-V)/2\), and \(X=P_j(W)\).
		For \(\gamma_1=\|X\|_2^2\), the midpoint identities give
		\[
		(1-d)|\tau(V^*T_j(U))|
		=\gamma_1-\|P_j(W_-)\|_2^2+d-2d\|W\|_2^2\le \gamma_1+d-2d\|W\|_2^2.
		\]
		Write \(X=J_X|X|\) with a unitary polar extension, and put
		\(a=\operatorname{Re}\bE_{W^*(|X|)}(J_X^*W)\). 
		It is a self-adjoint contraction with \(\tau(a|X|)=\gamma_1\) and \(\tau(a^2)\le \|W\|_2^2\).
		Consequently \(\gamma_1-d \|W\|_2^2\le\tau(H_d(|X|))\). 
		Apply Proposition~\ref{prop:fiber-semicircle} to \(H_d\), and then use
		\(H_d(|Y|)\le q(|Y|^2)\), since \(\|Y\|\le2\). 
		The midpoint inequality becomes
		\[
		|\tau(V^*T_j(U))|
		\le \frac{d+2\tau(q(|Y|^2))-\gamma_1}{1-d}+o(1),
		\]
		with the error uniform over the unitary pair.
		The first-moment argument in Corollary~\ref{cor:finite-positive-transfer} gives $\gamma_1\le(8/(3\pi))^2+o(1)<0.721$ uniformly for all sufficiently large $j$.
		On this interval \(\alpha_3+8\gamma_1\alpha_4>0\), by evaluating it at the right endpoint.
		The derivative \(10\alpha_3+\alpha_4(28+84\gamma_1)\) of
		\(2\alpha_2+5\alpha_3(1+2\gamma_1)+\alpha_4(28\gamma_1+42\gamma_1^2)\) is positive there, and the latter function takes the value \(-0.053583551746269364\) at \(\gamma_1=0.721\).
		Thus the coefficient of \(\Delta\) is negative. 
		The resolvent comparison and the joint moment estimate bound the limiting objective by the maximum of the concave function
		\[
		G(\gamma_1)=\frac{d+
			2[\alpha_0+\alpha_1\gamma_1+2\alpha_2 \gamma_1^2+5\alpha_3 \gamma_1^3+14\alpha_4\gamma_1^4+h_0J(\gamma_1/0.247)]-\gamma_1}{1-d}.
		\]
		The quadratic \(4\alpha_2+30\alpha_3\gamma_1+168\alpha_4 \gamma_1^2\) is negative because its discriminant and leading coefficient are both negative. 
		Also \(J\) is concave, as is seen from its integral definition.
		At $s_0=0.701327218$, let $l=22576636337458/10^{13}$ and $r=l+10^{-13}$.
		Then $l<(1+\sqrt{1+4s_0/0.247})/2<r$ by squaring.
		Substitution in $J$ and $J'$ yields \(G(s_0)<0.701400373\) and \(|G'(s_0)|<10^{-9}\), using \(J'(z)=[h(z)^2(2h(z)-1)]^{-1}\).
		The tangent at \(s_0\) bounds \(G\) strictly below \(0.701400374\) on \([0,1]\). 
		All terminating decimals here denote rational numbers.
	\end{proof}
	
	\begin{proof}[Proof of Theorem~\ref{thm:fiber-main}]
		Corollary~\ref{cor:fiber-finite-depth} gives the irreducible finite-index inclusions of hyperfinite factors. 
		Corollary~\ref{cor:finite-positive-transfer} gives the limiting ratio for their positive projections. 
		For the maps \(T_j\) of Proposition~\ref{prop:fiber-eighth-joint}, the range of \(P_j\) is nonzero and \(T_j\) acts as the identity there. 
		On the orthogonal complement, which contains \(1\), it acts by \(-1291/8709\). 
		Hence \(\|T_j\|_{2\to2}=1\). 
		For all sufficiently large \(j\), \(\|T_j\|_{\infty\to1}<0.701400374\), which proves the strict lower bound after taking reciprocals.
	\end{proof}

	\begin{conjecture}\label{conj:finite-index-supremum}
		The previous investigation does not imply that $\KIf=\KI$.
		Hence we would like to conjecture that
		\begin{align*}
			\KIf=\KI.
		\end{align*}
	\end{conjecture}
	
	\subsection{Free-product reduction and attainment}
	
	In this section, we show that $\KI$ is attainable.
	
	Suppose that $\Phi$ is a bimodule map and $\Phi^{(k)}=\id_{M_k}\otimes \Phi$ is an amplification of $\Phi$.
	
	\begin{theorem}\label{thm:universal-free-product}
		Let $\cQ_0$ be the hyperfinite $\mathrm{II}_1$ factor and let $a$ be a variance-one circular element free from $\cQ_0$. 
		Then
		\[
		\KI=\KNM{\cQ_0}{\cQ_0*W^*(a)}.
		\]
		There is a norm-one $\cQ_0$-bimodule operator $S$ on $L^2(\cQ_0*W^*(a))$
		with $S(a)=S^*(a)=a$ and $\|S^{(k)}\|_{\infty\to 1}=1/\KI$ for every $k\ge1$.
	\end{theorem}
	
	\begin{proof}
		The inclusion $\cQ_0\subset \cQ_0*W^*(a)$ is irreducible.
		Indeed, the
		reduced-word decomposition identifies the orthogonal complement of $L^2(\cQ_0)$ with a multiple of the coarse $\cQ_0$-bimodule. 
		Conjugation by successive powers of a Haar unitary in $\cQ_0$ tends weakly to zero on this complement, as is seen on elementary tensors using the weak convergence of the powers to zero. 
		Thus every central vector for the two $\cQ_0$-actions lies in $L^2(\cQ_0)$ and is scalar.
		
		Choose irreducible inclusions $\cN_j\subset\cM_j$ and norm-one bimodule operators $\Phi_j$ with $\|\Phi_j\|_{\infty\to 1}\to1/\KI$.
		By Theorem~\ref{thm:fiber-main}, $\KI>1$, so we can assume $\|\Phi_j\|_{\infty\to 1}<1$. 
		In particular these inclusions are of $\mathrm{II}_1$ factors, since an irreducible inclusion of finite matrix factors is the identity inclusion. 
		Choose $p_j\in\cN_j$ with $\tau(p_j)=1/j!$ and let $\Psi_j$ be the restriction of $\Phi_j$ to the corner. 
		Matrix units in $\cN_j$ identify $\Phi_j$ with $\id_{M_{j!}}\otimes\Psi_j=\Psi_j^{(j!)}$, with normalized traces on the corner and on the matrix algebra. 
		The corner inclusion is irreducible, $\|\Psi_j\|_{2\to2}=1$, and whenever $k\mid j!$,
		\[
		\|\Psi_j^{(k)}\|_{\infty\to 1}\le \|\Psi_j^{(j!)}\|_{\infty \to 1}=\|\Phi_j\|_{\infty\to 1},
		\]
		by tensoring unitary tests with the identity.
		
		Use the circular singular-pair construction in the proof of Theorem~\ref{gb:theorem} for each $\Psi_j$. 
		Start with centered approximate singular vectors whose two residual norms are at most $1/j$.
		Include a trace-preserving copy $\cQ_j$ of $\cQ_0$ from the lower corner in the separable test algebra $\mathcal A_j$ for Popa's theorem \cite[Theorem~0.1(a)]{PopaFree14}. 
		If $w$ is the resulting Haar unitary, the algebras $w^n\mathcal A_jw^{-n}$ for $n\ge0$ are jointly free. 
		Use $n\ge1$ in the circular sums and reserve $n=0$ for $\cQ_j$. 
		The sums have uniformly bounded $L^p$ norms for every fixed finite $p$, by moment convergence.
		Clipping their moduli at two changes them by a term tending to zero in each such $L^p$, by the compact support of the limiting modulus law and uniform higher moments. 
		H\"older's inequality applied to telescoping differences of mixed moments with $\cQ_j$, then gives variance-one circular elements $X_j,Y_j$ in the iterated tracial ultrapower, each free from $\cQ_j$. 
		The compressed operator $\widetilde\Psi_j$ satisfies
		\[
		\|\widetilde\Psi_j(X_j)-Y_j\|_2\le1/j,\qquad
		\|\widetilde\Psi_j^*(Y_j)-X_j\|_2\le1/j .
		\]
		At each ultrapower, Hilbert-space compression and unitary lifting,
		also at fixed matrix size, give
		$ \|\widetilde\Psi_j^{(k)}\|_{\infty\to 1}\le \|\Psi_j^{(k)}\|_{\infty\to 1}$.
		
		Take a final tracial ultraproduct over $j$, use the diagonal copy
		of $\cQ_0$ supplied by the embeddings into the $\cQ_j$, and compress the Hilbert-space ultraproduct of the operators to the resulting tracial $L^2$ space. 
		The bounded sequences $X_j,Y_j$ give variance-one circular elements $X,Y$, each free from the diagonal copy of $\cQ_0$. 
		The resulting contraction $T$ is $\cQ_0$-bimodular and satisfies $TX=Y$ and $T^*Y=X$. 
		Since $k\mid j!$ for all large $j$, unitary lifting gives $\|T^{(k)}\|_{\infty\to 1}\le1/\KI$ for every fixed $k$.
		
		Use the trace-preserving isomorphisms $\pi_X:\cQ_0*W^*(a)\to \cQ_0*W^*(X)$ and $\pi_Y:\cQ_0*W^*(a)\to \cQ_0*W^*(Y)$ that fix $\cQ_0$ and send $a$ to $X$ and
		$Y$, respectively. 
		Put
		\[
		S=\pi_Y^{-1}\bE_{\cQ_0*W^*(Y)}\,T\,\pi_X
		\]
		on $L^2(\cQ_0*W^*(a))$, with the expectation denoting orthogonal projection.
		The two singular identities give $S(a)=S^*(a)=a$, so $S$ has norm one. 
		An amplified pairing for $S$ against unitaries $u,v$ is the pairing for $T$ against $\pi_X(u),\pi_Y(v)$, with the maps applied entrywise. 
		Hence $\|S^{(k)}\|_{\infty\to 1}\le1/\KI$ for each $k$. 
		By irreducibility of the free-product inclusion and the definition
		of $\KI$, $1/\KI\le \|S\|_{\infty\to 1}\le \|S^{(k)}\|_{\infty\to 1}$.
	\end{proof}

	\section{Structure of Grothendieck constants}
	
	In this section, we study the Grothendieck constants for group-subgroup inclusions, subfactor intervals, tensor products and free compositions.
	
	\subsection{Group-subgroup inclusions}
	
	\begin{theorem}\label{thm:group-subgroup}
		Suppose that a finite group $G$ acts outerly on a $\mathrm{II}_1$
		factor $\cQ$, and let $H\le G$. For the irreducible inclusion
		$\cN=\cQ\rtimes H\subset\cM=\cQ\rtimes G$,
		\[
		\KNM{\cN}{\cM}=1 .
		\]
	\end{theorem}
	\begin{proof}
		Write $u_g$ for the canonical group unitaries. 
		The summands $L^2(\cQ)u_g$ are pairwise inequivalent irreducible $\cQ$--$\cQ$ bimodules. 
		Indeed the space of intertwiners between the summands indexed by $g$ and $k$ is the space of intertwiners between the automorphisms $\alpha_g$ and $\alpha_k$ in $\cQ$. 
		A nonzero intertwiner is a scalar multiple of a unitary, contradicting
		outerness when $g\ne k$. 
		Thus an $\cN$--$\cN$ bimodule operator $\Phi$ acts by scalars $\lambda_g$ on these summands, with $\lambda_{h g k}=\lambda_g$ for $h,k\in H$. 
		For an index $g$ where $|\lambda_g|$ is maximal,
		\[
		\|\Phi\|_{2\to2}=|\lambda_g|
		=|\tau(u_g^*\Phi(u_g))|\le \|\Phi\|_{\infty\to 1}\le\|\Phi\|_{2\to2}.
		\]
		The last bound is followed by the Cauchy--Schwarz inequality on unitary tests. 
		Irreducibility follows from $\cQ'\cap(\cQ\rtimes G)=\bC1$ by outerness.
	\end{proof}
	
	\begin{theorem}
		Suppose that $\cN\subset\cM$ is an irreducible inclusion of
		$\mathrm{II}_1$ factors with finite index $\delta^2\ge4$ and
		$\dim(\cN'\cap\cM_1)=2$. Then $\KNM{\cN}{\cM}=1$.
	\end{theorem}
	\begin{proof}
		Popa's relative-dimension theorem at index four \cite[Theorem 5.2]{Popa89} and
		C\'aceres's theorem at index greater than four \cite[Theorem A]{Caceres}
		give a projection $p\in\cM$ with $\bE_{\cN}(p)=1/2$. The latter
		theorem applies since irreducibility implies extremality. Thus $u=2p-1$
		is a unitary in $\ker\bE_{\cN}$.
		The two-box algebra is spanned by the Jones projection $e_{\cN}$ and $1-e_{\cN}$, so $\Phi=\alpha_1\bE_{\cN}+\alpha_2(\id-\bE_{\cN})$. 
		Testing $1$ and $u$ on both sides bounds $\|\Phi\|_{\infty\to 1}$ below by $\max\{|\alpha_1|,|\alpha_2|\}=\|\Phi\|_{2\to2}$.
	\end{proof}
	
	\begin{proposition}
		For the crossed-product inclusion in Theorem~\ref{thm:group-subgroup},
		\[
		\KNM{\cN}{\cM}_{\cb}\le\sqrt{[G:H]} .
		\]
	\end{proposition}
	\begin{proof}
		The direct sum of the right $\cN$-modules $u_g L^2(\cN)$, with $g$ ranging over representatives for $G/H$, gives $[\cM:\cN]=[G:H]$. 
		The proof of
		Theorem~\ref{thm:cb-algebra} gives $\|\Phi\|_{\cb}\le\sqrt{[\cM:\cN]}\,\|\Phi\|_{2\to2}$.
		Use $\|\Phi\|_{2\to2}=\|\Phi\|_{\infty\to 1}$ from the preceding
		crossed-product calculation.
	\end{proof}
	
	We use an equivariant matrix model for fixed-point inclusions.
	The first basic construction of $\cQ^G\subset\cQ$ is $\cQ\rtimes G$, with Jones projection $\displaystyle p=|G|^{-1}\sum_g u_g$.
	Its two-box algebra is $\bC[G]$, acting on $L^2(\cQ)$ by the automorphisms $\alpha_g$ \cite[Section 8.1]{JonesNotes}.
	If $\sigma$ is a finite-dimensional unitary representation of $G$,
	embed $B(V_\sigma)$ unitally in $\cQ^G$ and write $w_g$ for the image of $\sigma(g)$.
	Then $\displaystyle q=|G|^{-1}\sum_g w_g u_g$ and $p$ are projections of the same trace in the crossed-product factor. 
	Choose
	$v$ with $v^*v=p$ and $vv^*=q$. Since $u_gp=p$, the finite Fourier
	expansion gives $v=vp=ap$ with $a\in\cQ$. The identity coefficient
	of $apa^*=q$ gives $aa^*=1$, hence $a$ is unitary. The relation
	$w_g u_g v=v$ gives $w_g\alpha_g(a)=a$ by comparing the identity
	coefficients. Thus $b\mapsto a^*ba$ embeds
	$(B(V_\sigma),\operatorname{Ad}\sigma)$ equivariantly and preserves
	the normalized trace. In particular, the representation
	$\lambda_G\otimes 1_d$ gives a copy of $M_d$-valued functions on $G$
	with left translation and trace $|G|^{-1}\sum_g\operatorname{tr}_d$
	on its diagonal blocks.
	
	\begin{proposition}
		Suppose that $G$ is a finite abelian group acting outerly on
		a $\mathrm{II}_1$ factor $\cQ$. Then $\KNM{\cQ^G}{\cQ}=1$.
	\end{proposition}
	\begin{proof}
		Each character of $G$ occurs, up to conjugating the character, as a unitary eigenfunction for left translation in the equivariant diagonal copy of $C(G)$. 
		The abelian group algebra is diagonal on these character spaces. 
		For a two-box choose a character at which the corresponding scalar has maximal modulus and test its unitary eigenfunction in both slots.
	\end{proof}
	
	\begin{theorem}\label{thm:pos}
		Suppose that $G$ is a finite nonabelian group acting outerly on a
		$\mathrm{II}_1$ factor $\cQ$. Then
		\[
		\KNM{\cQ^G}{\cQ}\ge\frac{146}{145},
		\qquad \KNM{\cQ^G}{\cQ}_+>1 .
		\]
	\end{theorem}
	\begin{proof}
		Choose an irreducible unitary representation $\pi$ of $G$ with dimension
		$d\ge2$. 
		There are orthonormal vectors $\xi_1,\xi_2$ and a phase $e^{i\phi}$ such that the projective orbit of $(\xi_1+e^{i\phi}\xi_2)/\sqrt2$ contains $(\xi_1+z\xi_2)/\sqrt2$ for every $z$ in a rotated regular polygon $Z$ with at least three vertices. 
		If some $\pi(g)$ has an eigenvalue ratio of order at least three, use the corresponding eigenvectors and the powers of $g$, with $e^{i\phi}=1$.
		Otherwise every projective image has order at most two, so they commute projectively.
		Rephase the lifts to involutions. 
		Their pairwise scalar commutators are signs, and one pair anticommutes since commuting lifts would make $\pi$ reducible. 
		Choose $\xi_1$ fixed by the first involution and let $\xi_2$ be its image under the second. 
		The pair acts by sign change and interchange on their span. 
		Taking $e^{i\phi}=e^{i\pi/4}$ yields the four relative phases $e^{i\phi},-e^{i\phi},e^{-i\phi},-e^{-i\phi}$.
		
		Let $E_{ij}$ be the operator matrix units in this group-algebra block, and put $V=(E_{11}+e^{i\phi} E_{21})/\sqrt2$. 
		Then
		$V^*V=E_{11}$ and $\|V\|_{2\to2}=1$. 
		For unitaries $u,w\in\cQ$, absorb a phase into $w$ so that $\tau(w^*V(u))\ge0$. 
		Put $a=E_{11}(u)$, $y=E_{21}(u)$, $b=V(u)$.
		Orthogonality gives $\|u-a\|_2^2+\|w-b\|_2^2=2(1-\tau(w^*V(u)))$, and $\|a\|_2=\|y\|_2=\|b\|_2\le1$.
		For $z\in Z$, the orbit property shows that
		$b_z=(a+zy)/\sqrt2$ is a scalar multiple of an automorphic image
		of $b$. 
		Hence its distance to a unitary is at most $\|w-b\|_2$, and $\|b_z^*b_z-1\|_1\le2\|w-b\|_2$ by H\"older's inequality. 
		The zeroth and first Fourier coefficients of
		\[
		2(b_z^*b_z-1)=a^*a+y^*y-2+z a^*y+\bar z y^*a
		\]
		give $\|a^*a+y^*y-2\|_1\le4\|w-b\|_2$ and $\|a^*y\|_1\le4\|w-b\|_2$,
		since $\displaystyle \sum_{z\in Z}z=\sum_{z\in Z}z^2=0$. 
		We also have
		$\|a^*a-1\|_1\le2\|u-a\|_2$, so $\|y^*y-1\|_1\le4\|w-b\|_2+2\|u-a\|_2$.
		The scalar inequality $|\sqrt{t}-1|\le |t-1|$ for $t\ge0$ and H\"older's inequality give
		\[
		1-4\|w-b\|_2-2\|u-a\|_2\le\|y\|_1=\|u^*y\|_1\le\|a^*y\|_1+\|u-a\|_2\le 4\|w-b\|_2+\|u-a\|_2 .
		\]
		Thus $1\le 8\|w-b\|_2+3\|u-a\|_2\le\sqrt{146(1-\tau(w^*V(u)))}$. 
		Taking the supremum over $u,w$ proves $\|V\|_{\infty\to 1}\le 145/146$.
		
		For the strict positive gap, we show that some irreducible representation of a subgroup has a unit vector whose orbit projections span its whole matrix algebra.
		First, there is an abelian subgroup $H_0\le G$ and $g\in N_G(H_0)$ acting nontrivially on $H_0$, where $N_G(H_0)$ is the normalizer of $H_0$ in $G$.
		Indeed, in a nonabelian Sylow subgroup $H$ one can take
		$H_0=\langle x,Z(H)\rangle$ for $x\in Z_2(H)\setminus Z(H)$, where
		$Z_2(H)/Z(H)=Z(H/Z(H))$, $Z(H)$ is the center of $H$. 
		The subgroup $H_0$ is abelian and normal in $H$, and some element of $H$ moves $x$. 
		If no such $H_0$ exists, each Sylow subgroup $H$ is abelian and $H\subset Z(N_G(H))$.
		Burnside's normal complement theorem \cite[Theorem 1.21]{Craven} then embeds $G$ into the product of its abelian Sylow quotients by their normal complements. 
		The kernel is trivial because the orders of the normal complements have greatest common divisor one.
		This contradicts noncommutativity.
		
		Let $H_1=\langle H_0,g\rangle$. 
		Pick an $H_0$-character $\chi$ with $g$-orbit of size $d\ge2$. 
		Its stabilizer in $H_1$ is $I=\langle H_0,g^d\rangle$. 
		The character extends to $I$ by assigning to $g^d$ an $l$th root of $\chi(g^{dl})$, where $l$ is the order of $g^d H_0$ in $I/H_0$. 
		Invariance under $g^d$ verifies the multiplication relation. 
		Induction from $I$ gives a representation with $d$ distinct one-dimensional weight spaces for $H_0$, permuted transitively by $g$. 
		Its commutant consists of scalars, so it is irreducible. 
		In a weight basis, the conjugates of $\displaystyle (\sum_i e_i)e_1^*$ span all $e_i e_1^*$ by Fourier projection onto their distinct $H_0$-weights. Conjugating by powers of $g$ then gives every matrix unit. 
		Choose $d^2$ conjugates of this rank-one matrix which are independent. 
		The determinant formed by the same conjugations of $xy^*$ is a nonzero polynomial in $x,\bar y$. 
		Its restriction to $x=y=z$ is still nonzero as a polynomial in the real and imaginary parts of $z$, since $(s,t)\mapsto(s+it,s-it)$ is an invertible complex-linear change of variables. 
		Some unit vector $\xi$ therefore has orbit projections spanning $M_d(\bC)$.
		
		Put $P_\xi$ in this representation block of $\bC[H_1]\subset\bC[G]$.
		It acts as an orthogonal projection. 
		By the Cauchy--Schwarz inequality, $\displaystyle \|P_\xi\|_{\infty\to 1}=\sup_{u\in\U(\cQ)}\|P_\xi(u)\|_2^2$.
		If this were one, a maximizing sequence would give a unitary $W$
		in the tracial ultrapower $\cQ^\omega$ fixed by $P_\xi$. 
		In a basis with $e_1=\xi$, put $X_i=E_{i1}(W)$. 
		For $z=\pi(h)\xi$ with $h\in H_1$, the element $\displaystyle \sum_i z_iX_i$ is the unitary $\alpha_h^\omega(W)$. 
		Therefore
		\[
		\sum_{i,j}\overline{z_i}z_j X_i^*X_j=1 .
		\]
		Spanning by the matrices $zz^*$ extends this identity linearly to all matrices and gives $X_i^*X_j=\delta_{ij}1$. 
		The diagonal identities make each $X_i$ an isometry, hence a unitary in the finite ultrapower.
		The off-diagonal identities are then impossible.
		Thus $\|P_\xi\|_{\infty\to 1}<1$.
	\end{proof}
	
	For a finite complex matrix $A=(a_{ij})$, define
	\begin{align*}
		\mathscr S(A)&=\sup_{|z_i|=|w_j|=1}
		\Big|\sum_{ij}a_{ij}\bar z_iw_j\Big|,\\
		\mathscr U(A)&=\sup_{d\ge1}\ \sup_{V_i,W_j\in\U(d)}
		\Big|\sum_{ij}a_{ij}\operatorname{tr}_d(V_i^*W_j)\Big|,\\
		\mathscr V(A)&=\sup_{\|x_i\|=\|y_j\|=1}
		\Big|\sum_{ij}a_{ij}\langle x_i,y_j\rangle\Big|.
	\end{align*}
	Here $\operatorname{tr}_d$ is normalized and the Hilbert inner product
	is linear in its second argument. 
	Since matrix unitaries are unit vectors for the normalized Hilbert--Schmidt norm, $\mathscr S(A)\le\mathscr U(A)\le\mathscr V(A)$.
	The classical complex Grothendieck constant is $\displaystyle K_{\mathrm{Gr}}^{\bC}=\sup_{A\ne0}\mathscr V(A)/\mathscr S(A)$.
	
	\begin{theorem}\label{thm:hecke}
		Let $G$ be a finite group acting outerly on a $\mathrm{II}_1$
		factor $\cQ$, and let $H\le G$. Then
		\[
		1\le \KNM{\cQ^G}{\cQ^H}\le
		\sup_{A\ne0}\frac{\mathscr V(A)}{\mathscr U(A)}
		\le K_{\mathrm{Gr}}^{\bC}<1.405,
		\qquad
		1\le \KNM{\cQ^G}{\cQ^H}_+\le \frac4\pi .
		\]
	\end{theorem}
	\begin{proof}
		Let $\displaystyle e_H=|H|^{-1}\sum_{h\in H}u_h$. 
		Compression to $e_H L^2(\cQ)=L^2(\cQ^H)$ identifies the operator two-box algebra with $e_H\bC[G]e_H$. 
		Indeed $e_H$ commutes with the two actions of $\cQ^G$ and a bimodule endomorphism on its range extends by zero.
		The regular representation of $\bC[G]$ is faithful, and hence compression realizes the Hecke algebra faithfully on $\ell^2(X)$, $X=H\backslash G$. 
		Give $X$ uniform probability measure and write $T=(T_{xy})_{x, y\in X}$ for the resulting matrix. 
		It commutes with right translations, which act transitively on $X$.
		
		Choose unit singular functions $\xi_1,\xi_2\in\ell^2(X)$ with $T\xi_1=\|T\|_{2\to2} \xi_2$.
		For $x\in X$ take the vectors $F_x(g)=\xi_1(xg)$ and $J_x(g)=\xi_2(xg)$ in $L^2(G)$ with normalized counting measure. 
		Transitivity gives unit norm, and
		$\displaystyle \sum_y T_{xy}F_y=\|T\|_{2\to2} J_x$ by commutation with right translations.
		For $A=T/|X|$ their pairing is $\|T\|_{2\to2}$. 
		Conversely, letting $T$ act on a Hilbert-space-valued tuple and using Cauchy--Schwarz inequality gives $\mathscr V(A)\le \|T\|_{2\to2}$. 
		By the equivariant model above, block-diagonal matrix functions on $G$ constant on left $H$-cosets belong to $\cQ^H$, with exactly the normalized trace pairings for $A$.
		Thus $\mathscr U(A)\le \|T\|_{\infty \to 1}$, where $T$ also denotes the same box acting on $L^2(\cQ^H)$.
		
		When $T\ge0$, the preceding vectors may be chosen with $\xi_1=\xi_2$.
		The scalar phase bound for $A\ge0$ follows by complex Gaussian
		rounding \cite{FL}.
		In detail, write
		\begin{align*}
			\mathbb E[\overline{\operatorname{phase}\langle Z,F_x\rangle}
			\operatorname{phase}\langle Z,F_y\rangle]
			=(\pi/4)\sum_{j\ge0}c_j z_{xy}|z_{xy}|^{2j},
		\end{align*}
		where $z_{xy}=\langle F_x,F_y\rangle$, $Z$ is a standard circular Gaussian vector, $c_0=1$ and $c_j\ge0$ \cite{HaagerupC}.
		Every term paired with $A$ is nonnegative since it is the quadratic pairing of $A\ge0$ with a tuple of tensor vectors $F_x^{\otimes(j+1)}\otimes\overline{F_x}^{\otimes j}$.
		The leading term is $(\pi/4)\|T\|_{2\to 2}$. 
		Consequently $\|T\|_{\infty \to 1}\ge\mathscr S(A)\ge (\pi/4)\|T\|_{2\to 2}$.
		The remaining bounds follow from $\mathscr S(A)\le\mathscr U(A)$,
		Haagerup's bound \cite{HaagerupC}, and the identity test.
	\end{proof}

	Let $\displaystyle \mathbf K^G=\sup_{\cQ^G\subset\cQ^H}\KNM{\cQ^G}{\cQ^H}$,
	where the supremum is over finite groups $G$, subgroups $H$, and
	outer actions on $\mathrm{II}_1$ factors.
	
	\begin{theorem}\label{thm:supremum}
		Let $\cQ_0$ be the hyperfinite $\mathrm{II}_1$ factor and choose
		one outer action of each symmetric group $S_n$ on $\cQ_0$. Then
		\[
		1.2363<\mathbf K^G
		=\sup_{n\ge2}\KNM{\cQ_0^{S_n}}{\cQ_0}
		=\sup_{A\ne0}\frac{\mathscr V(A)}{\mathscr U(A)}
		\leq K_{\mathrm{Gr}}^{\bC}.
		\]
	\end{theorem}
	\begin{proof}
		The upper bounds follow from Theorem~\ref{thm:hecke}. 
		To prove realization by $S_n$, fix optimal
		unit vectors $x_i,y_j\in\bC^r$ for $\mathscr V(A)$. We can approximate
		their Gram matrix by permutations of one vector $f$. 
		In detail, take independent vectors $Z_1,\ldots,Z_n$
		uniform on the sphere of radius $\sqrt r$, and for each vector $x$
		in the finite list take the sampled tuple
		$F_x(k)=\langle Z_k,x\rangle$. With normalized counting measure
		the sample inner products converge to $\langle x,y\rangle$ since
		$\bE(Z_k Z_k^*)=1$. The coordinate laws are identical for all
		unit vectors in the list. Partition their common compact disk of
		values into finitely many cells of diameter at most $\varepsilon$.
		The empirical proportions in these cells converge to the same
		probabilities for all tuples. Hence each tuple can be matched by
		a permutation of a reference tuple with squared error bounded
		by $\varepsilon^2$ on matched coordinates and by $4r$ on an $o(1)$
		fraction of unmatched coordinates. Normalize the reference
		tuple, whose norm tends to one. Taking $n$ large and $\varepsilon$
		small gives a unit vector $f$ in the permutation representation
		$\rho$ of $S_n$ and permutations $g_i,h_j$ with
		\[
		\langle\rho(g_i)f,\rho(h_j)f\rangle
		=\langle x_i,y_j\rangle+o(1).
		\]
		For $\displaystyle T=\sum_{ij}a_{ij}g_i^{-1}h_j\in\bC[S_n]$, the matrix element
		$\langle f,\rho(T)f\rangle$ tends to the optimal vector pairing, so
		$\|T\|_{2\to2}\ge\mathscr V(A)-o(1)$ in the faithful two-box
		representation.
		For unitaries $u,v\in\cQ_0$,
		\[
		\tau(v^*T(u))=
		\sum_{ij}a_{ij}\tau\big(\alpha_{g_i}(v)^*\alpha_{h_j}(u)\big).
		\]
		All unitaries in this finite tuple can be approximated in $L^2$
		by unitaries in one matrix subalgebra of $\cQ_0$. 
		Indeed, a trace-preserving expectation onto an increasing dense sequence of matrix subalgebras
		gives contractions $w$ converging in $L^2(\cQ_0)$ to a given unitary.
		The unitary extension of the polar part of $w$ has squared distance $\tau((1-|w|)^2)\le1-\|w\|_2^2$ from $w$.
		Therefore $\|T\|_{\infty \to 1}\le\mathscr U(A)$.
		Taking suprema gives both equalities.
		
		For the strict numerical lower bound, let $P$ be the Gaussian orthogonal projection onto the span of the coordinate functions $g_j$ on $\bC^r$ with standard circular Gaussian measure $\mu$.
		Suppose that $b_*$ is the unique root of the following equation
		\[
		b+(1-b)J(b)=\frac{\sqrt\pi}{2},
		\qquad J(b)=\frac12\int_0^\infty\frac{e^{-v}}{\sqrt{v+b}}\,dv,
		\]
		in $(0,1)$.
		By taking \(\phi(b)=b+(1-b)J(b)\), we have that \(\phi(0)=\sqrt\pi/2\), \(\phi'(0+)=-\infty\), \(\phi''(b)>0\) on \((0,1)\), and \(\phi(1)=1\). 
		These facts give the unique positive root \(b_*\) and $\phi(b)\le\sqrt\pi/2$ implies that $ b\le b_*$.
		Suppose that  $F:\bC^r\to\U(d)$ is a measurable function.
		Put \(A_j=\mathbb E(\overline{g_j}F)\), \(\displaystyle L=\sum_jA_jA_j^*\), \(\displaystyle R=\sum_jA_j^*A_j\), and \(b=\operatorname{tr}_d L=\operatorname{tr}_d R\).
		Bessel's inequality applied to $F(g)v$ for each vector $v$ gives $R\le1$. 
		Applying it to $F(g)^*v$ with the conjugate Gaussian basis gives $L\le1$.
		Then $\displaystyle PF=\sum_j g_j A_j$.
		We have that $b=\mathbb E\operatorname{tr}_d(F^*PF)\le \mathbb E\operatorname{tr}_d|PF|$.
		Set \(X_s=\sum_j A_jB_j(s)\), with independent standard complex Brownian motions $B_j(s)$. 
		In the singular bases of \(X_s\), the matrix generator has row and column sums given by the diagonal entries of \(L\) and \(R\). The two deficits \(I-L\) and \(I-R\) therefore give weighted spectral measures of mass \(1-b\).
		The first moment of the left measure is
		\[
		\mathbb E\operatorname{tr}_d((I-L)X_sX_s^*)
		=s\operatorname{tr}_d(L-L^2)\le sb(1-b).
		\]
		The right measure has the same bound.
		For a standard circular Gaussian \(g\), set
		\[
		f_t(x)=\mathbb E|\sqrt{x}+\sqrt t\,g|,\qquad h_t(x)=xf_t'(x).
		\]
		The scalar heat equation gives \(\partial_t f_t=h_t'\), and
		\[
		h_t'(x)=\frac1{4\sqrt\pi}\int_0^\infty
		\frac{u^{-1/2}}{1+tu}\exp\!\left(-\frac{xu}{1+tu}\right)\,du
		\]
		shows that \(h_t'\) is decreasing and convex. 
		For \(X_s=U\operatorname{diag}(\sqrt{\lambda_i})V^*\), differentiating the spectral trace gives
		\[
		\mathcal L\operatorname{tr}_d f_t(X_s^*X_s)
		=\frac1d\sum_{i,j}w_{ij}h_t[\lambda_i,\lambda_j],
		\qquad
		w_{ij}=\sum_\ell |(U^*A_\ell V)_{ij}|^2.
		\]
		Here \(h_t[a,b]\) is the divided difference. In the trace differentiation, the coefficient of \(|(U^*A_\ell V)_{ij}|^2\) is \(f_t'(\lambda_j)+\lambda_i f_t'[\lambda_i,\lambda_j]=h_t[\lambda_i,\lambda_j]\).
		Convexity bounds the divided difference by the average of the two endpoint derivatives. 
		Subtracting the time derivative at \(t=1-s\) leaves the two deficits. 
		The mass and first-moment bounds then give
		\[
		\frac{d}{ds}\,\mathbb E\operatorname{tr}_d f_{1-s}(X_s^*X_s)
		\le -(1-b)h'_{1-s}(sb).
		\]
		Integrate from \(0\) to \(1\). 
		The endpoint at \(1\) follows from \(|f_t(x)-\sqrt{x}|\le\sqrt t\,\mathbb E|g|\). 
		Substitution in the Laplace integral gives \(\displaystyle \int_0^1h'_{1-s}(sb)\,ds=J(b)\).
		Hence
		\[
		\mathbb E\operatorname{tr}_d(|PF|)\le \frac{\sqrt\pi}{2}-(1-b)J(b).
		\]
		This yields \(b\le b_*\).

		For finite partitions $(E_a)$ of $\bC^r$, let $\displaystyle m_{aj}=\int_{E_a}g_j\,d\mu(g)$ and form the positive matrix whose entry is $\displaystyle A_{ac}=\sum_j m_{aj}\overline{m_{cj}}$.
		For step functions with unitary values the quadratic pairing of $A$ is $\mathbb E\operatorname{tr}_d(F^*PF)$, as is seen by
		expanding in the orthonormal coordinates $g_j$. 
		The Cauchy--Schwarz inequality for this positive pairing gives $\mathscr U(A)\le b_*$ also with two independent lists of unitaries. 
		For the vector test take $Y(g)=g/\|g\|_2$. 
		Unitary invariance gives $\mathbb E[\overline{g_i}Y_j(g)]=r^{-1}\mathbb E\|g\|_2\,\delta_{ij}$,
		whence
		\[
		\mathbb E\langle Y,PY\rangle
		=\frac{(\mathbb E\|g\|_2)^2}{r}
		=\frac{\Gamma(r+1/2)^2}{r\,\Gamma(r)^2}\longrightarrow1 .
		\]
		Approximating $Y$ in $L^2$ by finite unit-vector step functions
		preserves the energy in the limit since $P$ is contractive.
		The corresponding positive matrices thus give
		$\mathbf K^G\ge 1/b_*$.
		
		For the decimal comparison, at $b=808865/10^6$ one has
		$b+(1-b)J(b)-\sqrt\pi/2>2.28\cdot10^{-6}$.
		This follows from $J(b)=(\sqrt\pi/2)e^b\operatorname{erfc}(\sqrt b)$
		using $0.89936922340048<\sqrt b<0.89936922340050$,
		$3.141592653589793<\pi<3.141592653589794$, and the Taylor
		sums through term $20$ for the exponential and the error function. 
		Their absolute remainders here are bounded by $3/21!$ and $2/(21!\cdot43)$, respectively. 
		Since the root in $(0,1)$ is unique and $b+(1-b)J(b)<\sqrt\pi/2$ for small positive $b$,
		we get $b_*<808865/10^6$ and $1/b_*>1.2363$.
	\end{proof}
	
	\begin{remark}
		Haagerup's upper bound is $K_{\mathrm{Gr}}^{\bC}<1.40491\ldots$
		\cite{HaagerupC}. 
		Heilman, Jones, and Malavolta obtain
		$1.4<K_{\mathrm{Gr}}^{\bC}<1.404898554746$ \cite{HJM}.
		The conjectured value is $1/\gamma_\dagger=1.40457593466\ldots$, where
		\[
		\gamma_\dagger=\int_0^{\pi/2}
		\frac{\cos^2u}{\sqrt{1+\sin^2u}}\,du
		=\frac\pi4\sum_{k\ge0}(-1)^k c_k
		\]
		with the coefficients $c_k=\binom{2k}{k}^2/[16^k(k+1)]$ as in \cite[Eq. (2.2)]{HJM}.  
	\end{remark}
	
	\subsection{Subfactor intervals}
	
	In this section, we characterize the Grothendieck constant for subfactor intervals.
	
	\begin{theorem}\label{thm:mono}
		Let $\cN\subset\cM$ be an irreducible finite-index inclusion of
		$\mathrm{II}_1$ factors, and let $\cQ$ be an intermediate factor. 
		Then
		\[
		\max\{\KNM{\cN}{\cQ},\KNM{\cQ}{\cM}\}\le\KNM{\cN}{\cM}.
		\]
		If $e_{\cQ}\in Z(\cN'\cap\cM_1)$ and
		\[
		(1-e_{\cQ})(\cN'\cap\cM_1)(1-e_{\cQ})
		\subseteq\cQ'\cap\langle\cM,e_{\cQ}\rangle ,
		\]
		then
		\begin{equation}\label{eq:intermediate-exact-max}
			\KNM{\cN}{\cM}=\max\{\KNM{\cN}{\cQ},\KNM{\cQ}{\cM}\}.
		\end{equation}
		Both assertions hold also for the operator-positive Grothendieck constants.
	\end{theorem}
	\begin{proof}
		For an intermediate subinterval $\cN\subseteq\mathcal A\subseteq \mathcal B\subseteq\cM$, extend a bimodule map $T$ on $\mathcal B$ by $\widetilde T=T\bE_{\mathcal B}$. 
		The extension is zero off $L^2(\mathcal B)$, so it preserves the Hilbert norm and operator positivity.
		For $u,v\in\U(\cM)$,
		\[
		\tau(v^*\widetilde T(u))
		=\tau(\bE_{\mathcal B}(v)^*T(\bE_{\mathcal B}(u))).
		\]
		As expectations are contractive, this is bounded by $\|T\|_{\infty \to 1}$.
		The reverse bound follows by restricting the tests to $\mathcal B$. 
		This gives monotonicity on every subinterval.
		
		For the upper bound, the hypotheses decompose any $T$ as $T_0=e_{\cQ}T e_{\cQ}$ and $T_1=(1-e_{\cQ})T(1-e_{\cQ})$.
		The operator $T_1$ is $\cQ$--$\cQ$ bimodular and the Hilbert norm of $T$ is the maximum of the two block norms.
		Tests in $\cQ$ give $\left\|T_0|_{\cQ}\right\|_{\infty\to 1}\le\|T\|_{\infty \to 1}$.
		
		To test the other block, average $T$ under conjugations by the left action of $\U(\cQ)$ in the finite factor $\cM_1=J\cN'J$ on $L^2(\cM)$.
		The $L^2$-closed convex hull of the orbit contains the trace-preserving
		expectation onto the operators commuting with the left action of $\cQ$. 
		The approximating convex averages are uniformly bounded and converge strongly in the normal representation on $L^2(\cM)$.
		Each conjugation preserves the norm $\|\cdot\|_{\infty\to 1}$, so the averaged operator has test norm at most $\|T\|_{\infty\to 1}$. 
		The average fixes $T_1$.
		On $L^2(\cQ)$ the averaged $T_0$ commutes with left $\cQ$ and right $\cN$, so it is scalar by $\cN'\cap\cQ=\bC1$.
		The average thus has the form $S=a e_{\cQ}+T_1$ and is $\cQ$--$\cQ$ bimodular. 
		We obtain
		\[
		\|T_1\|_{2\to2}\le\|S\|_{2\to2}
		\le\KNM{\cQ}{\cM}\,\|T\|_{\infty \to 1},\qquad
		\|T_0\|_{2\to2}\le\KNM{\cN}{\cQ}\,\|T\|_{\infty \to 1}.
		\]
		This gives the upper bound. 
		Compressions, zero extensions, and the average preserve positivity, so the same argument applies to operator-positive $T$.
	\end{proof}
	
	\begin{proposition}
		Suppose that the finite-index irreducible chain
		$\cN\subset\cQ\subset\cM$ satisfies the two corner hypotheses
		of Theorem~\ref{thm:mono}. Put $d_1=[\cQ:\cN]$, $d_2=[\cM:\cQ]$,
		and $\delta=\sqrt{d_1d_2}$. Then
		\begin{equation}\label{eq:cb-intermediate-upper}
			\begin{split}
				\KNM{\cN}{\cM}_{\cb}\le
				\min\bigl\{&
				\KNM{\cQ}{\cM}_{\cb}
				+\min\{\KNM{\cN}{\cQ}_{\cb},\sqrt{d_1-1}\,\KNM{\cN}{\cQ}\},\\
				&\KNM{\cN}{\cQ}_{\cb}+\sqrt{d_2-1}\,\KNM{\cQ}{\cM},\\
				&\delta\max\{\KNM{\cN}{\cQ},\KNM{\cQ}{\cM}\}\bigr\}.
			\end{split}
		\end{equation}
	\end{proposition}
	\begin{proof}
		Keep the operators $T_0,T_1$ and $S$ from the proof of  Theorem~\ref{thm:mono},
		using the restriction of $T_0$ when its domain is $\cQ$.
		Represent the downward construction for $\cN\subset\cM$ on $L^2(\cN)$, as $\mathcal L\subset\cN\subset\cM
		=J_{\cN}\mathcal L'J_{\cN}$. 
		Taking commutants identifies $\cQ$ with $J_{\cN}\mathcal P'J_{\cN}$ for an intermediate factor $\mathcal L\subset\mathcal P\subset\cN$.
		% Let $r=e_{\mathcal L}$ and $q=e_{\mathcal P}$ in this representation.
		Then $e_{\mathcal L}\le e_{\mathcal P}$, $\tau(e_{\mathcal P})=1/d_1$, and
		$\tau(e_{\mathcal L})=1/(d_1d_2)$.
		Since $\bE_{\cQ}(e_{\mathcal L})$ lies in
		$e_{\mathcal P}\cQ e_{\mathcal P}=\mathcal P e_{\mathcal P}$ and commutes with $\mathcal L$, irreducibility of $\mathcal L\subset\cN$ yields $\bE_{\cQ}(e_{\mathcal L})=e_{\mathcal P}/d_2$.
		
		Set $y_0=T_0(e_{\mathcal P})\in\mathcal P'\cap\cQ$. 
		The projection $e_{\mathcal P}$ onto the unique trivial $\mathcal P$--$\mathcal P$ summand of $L^2(\cN)$ commutes with $y_0$, and $e_{\mathcal P} y_0=a_0 e_{\mathcal P}$ for a scalar $a_0$.
		The scalar $a$ from the averaging of $T_0$ is the normalized basic-construction trace of its restriction on $L^2(\cQ)$.
		For a right frame $(v_i)$ for $\mathcal P\subset\cN$, the elements
		$w_i=\sqrt{d_1} v_i e_{\mathcal P}$ form a right frame for $\cN\subset\cQ$ by the frame identity and $\bE_{\cN}(e_{\mathcal P})=d_1^{-1}1$.
		Consequently
		\[
		a=d_1^{-1}\sum_i\tau(w_i^*T_0(w_i))
		=\sum_i\tau(e_{\mathcal P} v_i^*v_i y_0)
		=a_0\sum_i\tau(v_i e_{\mathcal P} v_i^*)=a_0 .
		\]
		By $\cQ$-bimodularity $T_1(e_{\mathcal L})=e_{\mathcal P} T_1(e_{\mathcal L})e_{\mathcal P}$. 
		Hence
		$T(e_{\mathcal L})=T_1(e_{\mathcal L})+d_2^{-1}y_0$ is block diagonal with respect to $e_{\mathcal P},1-e_{\mathcal P}$, and its $e_{\mathcal P}$ block equals $S(e_{\mathcal L})$.
		The trace formula \eqref{eq:downward-cb-trace} gives
		\[
		\|T\|_{\cb}=\|S\|_{\cb}+d_1\tau(|(1-e_{\mathcal P})y_0|),
		\qquad d_1\tau(|(1-e_{\mathcal P})y_0|)=\|T_0|_{\cQ}\|_{\cb}-|a|.
		\]
		The last term is at most
		$d_1\sqrt{1-1/d_1}\,\|T_0(e_{\mathcal P})\|_2
		\le\sqrt{d_1-1}\,\|T_0\|_{2\to2}$.
		Together with $\|S\|_{\infty\to 1}$ and $\|T_0|_{\cQ}\|_{\infty \to 1}$ bounded by $\|T\|_{\infty\to 1}$, this proves the first alternative.
		
		For the second, use the triangle inequality and the downward Jones projection $e_{-1}$ for $\cQ\subset\cM$, which satisfies
		$\bE_{\cQ}(e_{-1})=d_2^{-1}1$. 
		Since $T_1(1)=0$,
		\[
		\begin{split}
			\|T_1\|_{\cb}
			&=d_2\|T_1(e_{-1}-d_2^{-1}1)\|_1\\
			&\le d_2\|T_1\|_{2\to2}\sqrt{d_2^{-1}(1-d_2^{-1})}
			=\sqrt{d_2-1}\,\|T_1\|_{2\to2}.
		\end{split}
		\]
		The cb norm of the zero extension of $T_0|_{\cQ}$ is the same as its cb norm on $\cQ$, by complete contractivity of $\bE_{\cQ}$
		and restriction. 
		The last alternative follows from
		$\|T\|_{\cb}\le\delta\|T\|_{2\to2}$ in the proof of Theorem~\ref{thm:cb-algebra} and the exact maximum formula.
	\end{proof}
	
	\begin{corollary}\label{cor:three-box-upper}
		Suppose both constituent inclusions are proper and
		$\dim\PP_{2,+}^{\cN\subset\cM}=3$. Then
		\[
		\KNM{\cN}{\cM}=\max\{\KNM{\cN}{\cQ},\KNM{\cQ}{\cM}\}.
		\]
	\end{corollary}
	\begin{proof}
		The three nonzero orthogonal projections $e_{\cN}$,
		$e_{\cQ}-e_{\cN}$, and $1-e_{\cQ}$ span the two-box algebra.
		In particular $e_{\cQ}$ is central and the complementary corner is
		$\bC(1-e_{\cQ})\subseteq \cQ'\cap\langle\cM,e_{\cQ}\rangle$.
		Theorem~\ref{thm:mono} applies.
	\end{proof}
	
	\begin{corollary}\label{cor:two-by-two-upper}
		If $\cN\subset\cQ\subset\cM$ is irreducible and
		$[\cQ:\cN]=[\cM:\cQ]=2$, then
		\[
		\KNM{\cN}{\cM}=1,\qquad
		\KNM{\cN}{\cM}_{\cb}=2,\qquad
		\KNM{\cN}{\cQ}_{\cb}=\KNM{\cQ}{\cM}_{\cb}=\sqrt2.
		\]
	\end{corollary}
	\begin{proof}
		At index two the standard bimodule is the sum of the trivial
		sector and one invertible sector, since the statistical dimension
		of the standard bimodule is the index and that of a nonzero simple
		sector is at least one. 
		In an irreducible inclusion each invertible sector contains a unitary normalizer. 
		In fact a nonzero vector $\xi$ in the copy of a twisted identity correspondence satisfies
		$\xi n=\theta(n)\xi$ for an automorphism $\theta$.
		Thus $|\xi|$ is affiliated with the scalar relative commutant, and $\xi$ is a scalar multiple of a unitary. 
		Any two such unitary intertwiners differ by a scalar, so that sector has multiplicity one.
		Both constituent constants are therefore one.
		
		As an $\cN$--$\cN$ bimodule, $L^2(\cQ)$ contains two invertible
		sectors, each of multiplicity one in $L^2(\cM)$. 
		Hence $e_{\cQ}$ is central. 
		The complement has statistical dimension two and consists either of a single simple sector of dimension two or of two distinct invertible sectors. 
		In the first case the two-box corner off $L^2(\cQ)$ is scalar and Theorem~\ref{thm:mono} applies.
		In the second case tests on the four normalizing unitaries attain
		every block norm. 
		Thus $\KNM{\cN}{\cM}=1$ in both cases.
		
		For an inclusion of index $\delta^2$ the trace formula gives, for
		$T=a\bE_{\cN}+b(\id-\bE_{\cN})$,
		\[
		\|T\|_{\cb}
		=\frac{|a+(\delta^2-1)b|+(\delta^2-1)|a-b|}{\delta^2}.
		\]
		Indeed $T(e_{-1})=b e_{-1}+(a-b)\delta^{-2}1$ in the downward construction with Jones projection $e_{-1}$.
		At $\delta^2=4$ take $(a,b)=(-1,1)$, and at $\delta^2=2$ take $(a,b)=(1,i)$.
		These maps have Hilbert norm and test norm one and cb norms $2$ and $\sqrt2$, respectively. 
		The matching upper bounds are
		$\delta$ by Theorem~\ref{thm:cb-algebra} and the equalities just proved.
	\end{proof}
	
	\begin{proposition}\label{prop:disjoint-crossed-upper}
		Suppose that $\cN\subset\cQ$ is an irreducible inclusion of finite
		factors and $\cM=\cQ\rtimes_\alpha G$ for a finite group $G$ with
		$\alpha_g(\cN)=\cN$. Assume that the bimodules
		$\cH_g=L^2(\cQ)u_g$ are pairwise disjoint as
		$\cN$--$\cN$ bimodules. Then
		\[
		\KNM{\cN}{\cM}=\KNM{\cN}{\cQ}.
		\]
	\end{proposition}
	\begin{proof}
		Disjointness excludes an inner $\alpha_g$ for $g\ne e$, since such an automorphism gives an equivalence between $\cH_g$ and $\cH_e$. 
		Thus $\cM$ is a factor. 
		The trivial $\cN$--$\cN$ sector occurs only once, in $\cH_e$, so $\cN\subset\cM$ is irreducible.
		Disjointness also makes every bimodule operator $T$ block diagonal.
		Identify its block on $\cH_g$ with $T_g$ on $L^2(\cQ)$ via $x\mapsto xu_g$. 
		Right multiplication by $y\in\cN$ is transported to right multiplication by $\alpha_g(y)$. 
		As $\alpha_g(\cN)=\cN$, $T_g$ is $\cN$--$\cN$ bimodular. 
		Thus $\|T\|_{2\to2}=\max_g\|T_g\|_{2\to2}$. 
		For $x,y\in\U(\cQ)$, the test on $xu_g,yu_g$ equals $\tau(y^*T_g(x))$, so $\|T\|_{\infty\to 1}\ge \|T_g\|_{\infty\to 1}$. 
		This proves the upper bound, and zero extension from $\cQ$ proves the reverse bound as in Theorem~\ref{thm:mono}.
	\end{proof}
	
	\subsection{Tensor products}
	
	Let $m_{\cN\subset\cM}$ denote the largest matrix-block size of
	$\operatorname{End}_{\cN-\cN}L^2(\cM)$.
	
	\begin{theorem}\label{thm:tensor-intrinsic}
		Let $\cN^{(j)}\subset\cM^{(j)}$, $j=1,2$, be irreducible
		finite-index inclusions of $\mathrm{II}_1$ factors. Put
		$\cN=\cN^{(1)}\overline\otimes\cN^{(2)}$,
		$\cM=\cM^{(1)}\overline\otimes\cM^{(2)}$,
		$\mathcal S_1=\cM^{(1)}\overline\otimes\cN^{(2)}$ and
		$\mathcal S_2=\cN^{(1)}\overline\otimes\cM^{(2)}$.
		Then
		\[
		\KNM{\cN}{\mathcal S_j}\le\KNM{\cN^{(j)}}{\cM^{(j)}}\quad(j=1,2),
		\]
		and, with $m_{\min}=\min_j m_{\cN^{(j)}\subset\cM^{(j)}}$,
		\[
		\max_{j=1,2}\KNM{\cN}{\mathcal S_j}
		\le\KNM{\cN}{\cM}
		\le m_{\min}\KNM{\cN}{\mathcal S_1}\KNM{\cN}{\mathcal S_2}.
		\]
		These inequalities also hold for the operator-positive Grothendieck constants.
	\end{theorem}
	\begin{proof}
		Taking commutants of the left and right representations gives
		\[
		\operatorname{End}_{\cN-\cN}L^2(\cM)=\mathcal A_1\otimes \mathcal A_2,
		\qquad \mathcal A_j=\operatorname{End}_{\cN^{(j)}-\cN^{(j)}}
		L^2(\cM^{(j)}).
		\]
		The relative commutant is scalar by slicing, and tensoring the Jones projections gives the product of the indices. 
		The two-box algebra of $\cN\subset\mathcal S_j$ is $\mathcal A_j$ with the other leg acting as an identity.
		An original unitary test tensored with the identity has the original pairing, so its ratio cannot increase.
		This gives the first inequalities.
		The maximum lower bound follows from Theorem~\ref{thm:mono} for the intermediate factors $\mathcal S_j$.
		
		For unitary pairs $x_j,y_j\in\mathcal S_j$, let $\ell_j(T_j)=\tau(y_j^*T_j(x_j))$ on $\mathcal A_j$. 
		For $u\in\U(\cN)$, the unitaries $x=u x_1 u^*x_2$ and $y=u y_1 u^*y_2$ satisfy
		\[
		\tau(y^*(T_1\otimes T_2)(x))
		=\tau\bigl(u[y_1^*T_1(x_1)]u^*[T_2(x_2)y_2^*]\bigr).
		\]
		Indeed, each $T_j$ is bimodular over the opposite intermediate factor, and both are bimodular over $\cN$. 
		Extend the identity linearly to $T\in A_1\otimes A_2$. 
		Simultaneously average all the first brackets under $\U(\cN)$ in the trace $L^2$ norm.
		The averages converge to the scalar traces because the invariant
		vectors for conjugation are scalar. 
		It follows that
		\begin{equation}\label{eq:bgi-product-tests}
			|(\ell_1\otimes\ell_2)(T)|\le\|T\|_{\infty\to 1}.
		\end{equation}
		
		Write $\kappa_j=\KNM{\cN}{\mathcal S_j}$. For norm-one dual
		functionals $f_j$ on $\mathcal A_j$, apply $\kappa_1$ to the slice
		$(\id\otimes f_2)(T)$ and $\kappa_2$ to each further slice
		$(\ell_1\otimes\id)(T)$. Taking the suprema over $\ell_1,\ell_2$
		in \eqref{eq:bgi-product-tests} gives
		$|(f_1\otimes f_2)(T)|\le\kappa_1\kappa_2\|T\|_{\infty\to 1}$.
		In a block attaining the operator norm, take left and right unit
		singular vectors. 
		Each Schmidt expansion has at most $m_{\min}$ terms, and the sum of its coefficients is at most $\sqrt{m_{\min}}$ by the Cauchy--Schwarz inequality. 
		Expanding the singular pairing into product
		vector functionals gives the multiplier $m_{\min}$.
		
		For positive $T$, use states $f_1,f_2$ and diagonal unitary tests
		$y_j=x_j$. Both slices are then positive, and such tests give
		$\|\cdot\|_{\infty \to 1}$ for positive operators by the Cauchy--Schwarz inequality.
		The preceding slice argument therefore bounds $(f_1\otimes f_2)(T)$ by the product of the positive constants times $\|T\|_{\infty\to 1}$. 
		A top unit eigenvector has squared overlap at
		least $1/m_{\min}$ with some product vector by its Schmidt expansion. 
		The vector state of this product vector tests $T$ at least at $\|T\|_{2\to2}/m_{\min}$, proving the positive bound.
	\end{proof}
	
	\begin{remark}
		Combining with Theorem~\ref{gb:theorem}, we get
		\[
		\KNM{\cN}{\cM}\le
		\min\{\pi/2,\ m_{\min}\KNM{\cN}{\mathcal S_1}
		\KNM{\cN}{\mathcal S_2}\}.
		\]
	\end{remark}
	
	\subsection{Free compositions}
	
	In the section, we characterize the Grothendieck constant for free composition of irreducible inclusions.
	
	\begin{proposition}\label{prop:free-composition-intrinsic}
		Let $\cN\subset\cQ\subset\cM$ be an irreducible finite-index
		chain. 
		Suppose that the $\cQ$--$\cQ$ tensor categories generated by
		$L^2(\cQ)\otimes_{\cN}L^2(\cQ)$ and $L^2(\cM)$ are in free
		position, as in Bisch--Jones free composition
		\cite{BischJones,BrothierWeak}. 
		Then
		\[
		\KNM{\cN}{\cM}=\max\{\KNM{\cN}{\cQ},\KNM{\cQ}{\cM}\}.
		\]
		The same formula holds for the operator-positive Grothendieck constants.
	\end{proposition}
	\begin{proof}
		Let $F$ denote restriction of both actions from $\cQ$ to $\cN$.
		With $ A=L^2(\cQ)\otimes_{\cN}L^2(\cQ)$, Frobenius reciprocity
		on the two sides gives
		\[
		\operatorname{Hom}_{\cN-\cN}(F(X),F(Y))
		\cong\operatorname{Hom}_{\cQ-\cQ}
		(X,A\otimes_{\cQ}Y\otimes_{\cQ}A).
		\]
		The unit occurs once in $A$ by irreducibility of $\cN\subset\cQ$. 
		If $X,Y$ are nontrivial simples of the second category, only the summand using the unit on both sides of $Y$ can contribute. 
		Every other summand is a reduced alternating word involving the first category. 
		Thus for these simples restriction preserves the Hom spaces. 
		When $X$ is the unit and $Y$ is nontrivial in the second category, the right-hand side is zero, and taking adjoints gives the same vanishing with $X$ and $Y$ exchanged. 
		Decompose $L^2(\cM)$ as a $\cQ$--$\cQ$ bimodule into its single unit $L^2(\cQ)$ and the remaining summands. 
		The vanishing makes $e_{\cQ}$ central among the $\cN$--$\cN$ operators, and preservation of the other Hom spaces makes the complementary operators $\cQ$--$\cQ$ bimodular. 
		Theorem~\ref{thm:mono} now applies.
	\end{proof}
	
	\begin{proposition}
		Let $H,L$ be nontrivial finite groups, and suppose that their free
		product $G=H*L$ acts outerly on a $\mathrm{II}_1$ factor $\cQ$.
		Then
		\[
		\KNM{\cQ^H}{\cQ\rtimes L}=\KNM{\cQ^H}{\cQ}.
		\]
		Consequently the Bisch--Haagerup inclusions with $H=S_n$,
		$L=\bZ_2$, and $\displaystyle \cQ=\overline\bigotimes_{G}M_2(\bC)$ with the Bernoulli shift have supremum $\mathbf K^G$.
	\end{proposition}
	\begin{proof}
		For $\cQ^H\subset\cQ\subset\cQ\rtimes L$ the composite inclusion is irreducible by the outer Bisch--Haagerup
		construction \cite{BischHaagerup}.
		The $\cQ$--$\cQ$ bimodule $L^2(\cQ)\otimes_{\cQ^H}L^2(\cQ)$ identifies with $L^2(\cQ\rtimes H)$ by the basic construction of the fixed-point inclusion. 
		Its summands are the twists by $H$, and the summands of $L^2(\cQ\rtimes L)$ are the twists by $L$. 
		Their tensor products correspond to multiplication in $G$.
		Distinct products give inequivalent twists by outerness.
		The reduced-word decomposition of $H*L$ yields free position of the two categories.
		Apply Proposition~\ref{prop:free-composition-intrinsic} and $\KNM{\cQ}{\cQ\rtimes L}=1$ from Theorem~\ref{thm:group-subgroup}.
		
		The Bernoulli action is outer because for a nonidentity element one can choose a sequence of distinct tensor coordinates whose projections are moved to distinct coordinates at a fixed $L^2$-distance. 
		Any inner automorphism acts trivially in the limit on that central sequence. 
		The infinite tensor product is hyperfinite. 
		Theorem~\ref{thm:supremum} then proves the assertion about the supremum.
	\end{proof}
	
	\section{Bilinear Grothendieck Inequalities}
	
	For a bounded, separately normal bilinear map
	$\Phi:\cM\times\cM\to L^1(\cM)$, Haagerup's inequality \cite{Haagerup}, applied to $\tau(w^*\Phi(\cdot,\cdot))$ and then maximized over $\|w\|_\infty\le1$, gives
	\begin{align*}
		\left\|\sum_j \Phi(x_j,y_j)\right\|_1
		&\le 2\|\Phi\|_{\infty^2\to 1}
		\max\!\left\{\left\|\sum_jx_j^*x_j\right\|_\infty^{1/2},
		\left\|\sum_jx_jx_j^*\right\|_\infty^{1/2}\right\}\!\\
		&\hspace{24mm}\cdot
		\max\!\left\{\left\|\sum_jy_j^*y_j\right\|_\infty^{1/2},
		\left\|\sum_jy_jy_j^*\right\|_\infty^{1/2}\right\}.
	\end{align*}
	Balancing over an irreducible subfactor reduces the right side to
	one row norm and one column norm. 
	This argument uses the same relative averaging principle as the multimodular estimates in \cite{CPSS}.
	
	\begin{proposition}
		Suppose that $\cN\subset\cM$ is an irreducible inclusion of finite factors and $\Phi:\cM\times\cM\to L^1(\cM)$ is a bounded, separately normal bimodule bilinear map.
		Then
		\[
		\left\|\sum_j \Phi(x_j,y_j)\right\|_1
		\le 2\|\Phi\|_{\infty^2\to1}
		\left\|\sum_j x_j x_j^*\right\|_\infty^{1/2}
		\left\|\sum_j y_j^* y_j\right\|_\infty^{1/2}.
		\]
	\end{proposition}
	\begin{proof}
		Fix a contraction $w\in\cM$ and set $\phi(x,y)=\tau(w^*\Phi(x,y))$.
		Put $\beta=\|\Phi\|_{\infty^2\to1}$. 
		Haagerup's theorem gives normal
		states $f_1,f_2,g_1,g_2$ with
		\[
		|\phi(x,y)|\le\beta
		\bigl(f_1(x^*x)+f_2(xx^*)\bigr)^{1/2}
		\bigl(g_1(y^*y)+g_2(yy^*)\bigr)^{1/2}.
		\]
		Write $\displaystyle x_+=\sum_j x_j^*x_j$, $\displaystyle x_-=\sum_j x_jx_j^*$,
		$\displaystyle y_+=\sum_j y_j^*y_j$, and $\displaystyle y_-=\sum_j y_jy_j^*$.
		For $u\in\U(\cN)$, balancing gives
		$\phi(x_j,y_j)=\phi(x_j u,u^*y_j)$. 
		Summing and applying the Cauchy--Schwarz inequality gives
		\[
		\left|\sum_j\phi(x_j,y_j)\right|
		\le \beta
		\bigl(f_1(u^*x_+u)+f_2(x_-)\bigr)^{1/2}
		\bigl(g_1(y_+)+g_2(u^*y_-u)\bigr)^{1/2}.
		\]
		In the Hilbert space $L^2(\cM)\oplus L^2(\cM)$, the smallest-norm
		element in the closed convex hull of
		$\{(u^*x_+u,u^*y_-u):u\in\U(\cN)\}$ is its invariant vector.
		Irreducibility and trace preservation identify it as $(\tau(x_+)1,\tau(y_-)1)$. 
		Take finite convex combinations approaching this pair in $L^2(\cM)$. 
		The Cauchy--Schwarz inequality over the convex weights bounds the average of the product of square roots by the product of the square roots of the two averages. 
		The averaged operators are uniformly bounded, so their $L^2$ convergence implies convergence against the indicated normal states. 
		Hence
		\[
		\left|\sum_j\phi(x_j,y_j)\right|
		\le\beta\bigl(\tau(x_+)+f_2(x_-)\bigr)^{1/2}
		\bigl(g_1(y_+)+\tau(y_-)\bigr)^{1/2}
		\le 2\beta\|x_-\|_\infty^{1/2}\|y_+\|_\infty^{1/2}.
		\]
		Here $\tau(x_+)=\tau(x_-)\le\|x_-\|_\infty$ and $\tau(y_-)=\tau(y_+)\le\|y_+\|_\infty$. 
		Take the supremum over $w$.
	\end{proof}
	
	\begin{proposition}\label{cor:two-input-product}
		Suppose that $\cN\subset\cM$ is an irreducible inclusion of finite
		factors and $T:\cM\to\cM$ is a bounded normal bimodule map.
		For $\Phi_T(x,y)=T(xy)$ one has
		\[
		\|\Phi_T\|_{\infty^2\to1}=\|T\|_{\infty\to1},
		\qquad
		\|\Phi_T\|_{2^2\to1}=\|T\|_{1\to1}.
		\]
		If the inclusion has finite index $\delta^2$, then
		\[
		\|\Phi_T\|_{2^2\to1}=\|T\|_{\mathrm{cb}}
		\le\delta\,\KNM{\cN}{\cM}\,\|\Phi_T\|_{\infty^2\to1}.
		\]
		For operator-positive bimodule map $T$, one can use $\KNM{\cN}{\cM}_+$ instead.
	\end{proposition}
	\begin{proof}
		The three bimodule identities follow from those of $T$ and
		associativity. 
		Products of contractions are contractions, and every
		contraction is the product of itself and $1$. 
		Thus the first norm identity holds.
		H\"older's inequality gives $\|T(xy)\|_1\le\|T\|_{1\to1}\|x\|_2\|y\|_2$.
		Conversely, if $z=v|z|$ is the polar decomposition of $z\in\cM$,
		take $x=v|z|^{1/2}$ and $y=|z|^{1/2}$. 
		Then $xy=z$ and $\|x\|_2\|y\|_2=\tau(|z|)$, which proves the second identity.
		For finite index, trace duality gives $\|T\|_{1\to1}=\|T^\flat\|_{\infty\to\infty}$,
		which equals $\|T\|_{\mathrm{cb}}$ by \eqref{eq:downward-cb-trace}.
		The remaining bounds follow from Theorem~\ref{thm:cb-algebra} and the definitions of the two constants.
	\end{proof}
	
	\begin{theorem}\label{thm:general-two-input}
		Suppose that $\cN\subset\cM$ is an irreducible inclusion of
		$\mathrm{II}_1$ factors of finite index $\delta^2$.
		For a bounded, separately normal bimodule bilinear map $\Phi:\cM\times\cM\to\cM$ one has
		\[
		\|\Phi\|_{\infty^2\to1}\le\|\Phi\|_{2^2\to1}
		\le\min\{\delta^2,4+8\delta\}\,\|\Phi\|_{\infty^2\to1}.
		\]
	\end{theorem}
	
	\begin{proof}
		The first inequality follows from $\|x\|_2\le\|x\|_\infty$.
		Multiplication attains equality by H\"older's inequality and testing at $(1,1)$.
		Write $\beta=\|\Phi\|_{\infty^2\to1}$.
		
		\emph{Disjoint copies of a middle corner.}
		Let $x=xp$ and $y=py$ for a nonzero projection $p\in\cN$.
		Let $r_1=\cR(x)$ and $r_2=\cR(y^*)$ be the left support of $x$
		and the right support of $y$, respectively.
		The two supports of any operator have equal trace, so $\tau(r_i)\le\tau(p)$.
		The least projections in $\cN$ dominating them are $q_i=\cR(\bE_{\cN}(r_i))$, and the support estimate \eqref{eq:support} gives $\tau(q_i)\le\delta^2\tau(p)$.
		In particular, $q_1x=x$ and $yq_2=y$, so $z=\Phi(x,y)=q_1zq_2$.
		
		Given an integer $k$ with $k\delta^2\tau(p)\le1$, choose $v_i,u_i,w_i\in\cN$ with initial projections $q_1,p,q_2$, respectively, and pairwise orthogonal final projections within each family.
		Define $\displaystyle X=\sum_{i=1}^k v_ixu_i^*$ and $\displaystyle Y=\sum_{i=1}^k u_i y w_i^*$. 
		The summands of $X$ have orthogonal initial and final supports, and likewise for $Y$. 
		It follows that $\|X\|_\infty\le\|x\|_\infty$ and $\|Y\|_\infty\le\|y\|_\infty$. 
		By balancing and $u_i^*u_j=\delta_{ij}p$,
		\[
		\Phi(X,Y)=\sum_{i=1}^k v_i z w_i^*,\qquad
		\left|\sum_{i=1}^k v_i z w_i^*\right| =\sum_{i=1}^k w_i|z|w_i^* .
		\]
		For the second identity, the square of the left side is $\displaystyle \sum_i w_i z^*z w_i^*$, whose summands have orthogonal supports.
		Taking the trace gives
		\begin{equation}\label{eq:two-input-packing}
			k\|\Phi(x,y)\|_1\le\beta\|x\|_\infty\|y\|_\infty .
		\end{equation}
		
		\emph{Continuity and passage to the ultrapower.}
		For arbitrary $x,y\in\cM$ and a nonzero projection $p\in\cN$, split $p$ into $m$ equal-trace projections $p_i\in\cN$.
		By balancing, $\Phi(x,py)=\sum_i\Phi(xp_i,p_i y)$.
		Apply \eqref{eq:two-input-packing} with $k=\lfloor m/(\delta^2\tau(p))\rfloor$ to each term. 
		Since $m/k\to\delta^2\tau(p)$, one obtains
		$\|\Phi(x,py)\|_1\le\delta^2\tau(p)\beta\|x\|_\infty\|y\|_\infty$.
		Positive spectral sums and operator-norm approximation extend
		this inequality to positive $a\in\cN$ in place of $p$, with
		$\|a\|_1$ in place of $\tau(p)$. 
		For general $a=v|a|$, balance moves $v$ to the first input, yielding
		\[
		\|\Phi(x,ay)\|_1\le\delta^2\beta\|x\|_\infty\|y\|_\infty\|a\|_1.
		\]
		Put $a=\bE_{\cN}(x^*x)^{1/2}$ and $b=\bE_{\cN}(yy^*)^{1/2}$. 
		The Pimsner--Popa inequality \cite{PPindex} gives $x^*x\le\delta^2a^2$ and $yy^*\le\delta^2b^2$.
		Thus $x=x_0a$ and $y=by_0$ with $\|x_0\|_\infty,\|y_0\|_\infty\le\delta$. 
		For instance $x(a+\varepsilon 1)^{-1}$ is bounded by $\delta$ and converges strongly on the range of $a$ and its kernel to the required $x_0$.
		Balance and $\|ab\|_1\le\|a\|_2\|b\|_2=\|x\|_2\|y\|_2$ now give
		\begin{equation}\label{eq:two-input-preliminary}
			\|\Phi(x,y)\|_1\le\delta^4\beta\|x\|_2\|y\|_2 .
		\end{equation}
		
		Let $\omega$ be a free ultrafilter on $\bN$ and identify $\cM$ with its constant sequences in the tracial ultrapower $\cM^\omega$.
		The coordinatewise map $\Phi^\omega$ on $\cM^\omega$ is well defined. 
		Indeed, on bounded input sequences, changing a representative by an $L^2$-null sequence gives an $L^1$-null output by \eqref{eq:two-input-preliminary}. 
		The outputs have uniform operator-norm bounds, and
		$\|z\|_2^2\le\|z\|_\infty\|z\|_1$ then gives $L^2$-nullity.
		Truncating representatives of contractions and testing constant
		sequences gives $\|\Phi^\omega\|_{\infty^2\to1}=\beta$.
		The module identities hold over $\cN^\omega$.
		The support bound used in \eqref{eq:two-input-packing} also holds
		there. Lift a projection $r$ to projections $r_n\in\cM$, choose
		$q_n\in\cN$ with $q_n r_n=r_n$ and
		$\tau(q_n)\le\delta^2\tau(r_n)$, and take their ultraproduct $q$.
		
		\emph{The coefficient $\delta^2$.}
		Fix $x,y\in\cM$. Popa's relative freeness theorem
		\cite[Theorem~0.1(a)]{PopaFree14} supplies a Haar unitary
		$w\in\cN^\omega$ free from the algebra generated by $x^*x$ and
		$yy^*$. The relative commutant in the ultrapower is scalar.
		Indeed, relative $L^2$ averaging in $\cM$ gives
		$\|d-\tau(d)1\|_2\le\sup_{u\in\U(\cN)}\|[d,u]\|_2$ for each
		$d\in\cM$, and the unitaries on the right may be chosen coordinatewise.
		Take a partition $p_1,\ldots,p_m$ of $1$ in $W^*(w)$, with
		$\tau(p_i)=1/m$.
		
		For a centered self-adjoint $A$ and a projection $p$ free from $A$,
		the free compression formula \cite[Lecture~14]{NS} gives
		$\kappa_j^{(p)}(pAp)=t^{j-1}\kappa_j(A)$, where
		$t=\tau(p)$ and $\kappa_j^{(p)}$ uses the normalized corner trace.
		The moment--cumulant formula gives $|\kappa_j(A)|\le C_A^j$
		for some $C_A$. A noncrossing partition in an even moment of
		degree $2\ell$ that contributes has no singleton and at most
		$\ell$ blocks. Its contribution is bounded in absolute value by
		$t^{\ell}C_A^{2\ell}$. There are at most $4^{2\ell}$ such partitions,
		so the even moment is bounded by $(4C_A\sqrt t)^{2\ell}$.
		Taking $2\ell$th roots and letting $\ell\to\infty$ proves
		$\|pAp\|_\infty\le4C_A\sqrt t$. Applying this to the centered
		parts of $x^*x$ and $yy^*$ yields, uniformly in $i$,
		\[
		\|xp_i\|_\infty^2\le\|x\|_2^2+O_x(m^{-1/2}),
		\qquad
		\|p_i y\|_\infty^2\le\|y\|_2^2+O_y(m^{-1/2}).
		\]
		Use $\Phi(x,y)=\sum_i\Phi^\omega(xp_i,p_i y)$ and
		\eqref{eq:two-input-packing} in the ultrapower with
		$k_m=\lfloor m/\delta^2\rfloor$. For $m>\delta^2$ this gives
		\[
		\|\Phi(x,y)\|_1
		\le\frac{m}{k_m}\,\beta\,\max_i\|xp_i\|_\infty\max_i\|p_i y\|_\infty .
		\]
		As $m\to\infty$, the right side has upper limit
		$\delta^2\beta\|x\|_2\|y\|_2$.
		
		\emph{The coefficient $4+8\delta$.}
		Fix one $p=p_i$ and write $x_p=xp$, $y_p=py$ and
		$z=\Phi^\omega(x_p,y_p)=v|z|$. Put
		$r_1=\cR(x_p)$ and $r_2=\cR(y_p^*)$. Both have trace at most $1/m$.
		The right support of $z$ is dominated by the least
		$\cN^\omega$-projection dominating $r_2$. The support bound
		therefore gives $\tau(v^*v)\le\delta^2/m$.
		
		Choose $u_1,\ldots,u_m\in\cN^\omega$ with common initial projection
		$p$ and orthogonal final projections. By successive applications of
		Popa's theorem, choose $2m$ jointly free Haar unitaries
		$v_j,w_j\in\cN^\omega$ which are free from the algebra generated
		by $x_p,y_p,z,v,r_1,r_2,p$ and all $u_j$. Set
		\[
		X=\sum_j v_jx_p u_j^*,\qquad
		Y=\sum_j u_jy_p w_j^*,\qquad W=\sum_j v_j v w_j^*.
		\]
		Orthogonality and balancing give
		\[
		XX^*\le\|x_p\|_\infty^2\sum_j v_j r_1 v_j^*,\qquad
		Y^*Y\le\|y_p\|_\infty^2\sum_j w_jr_2w_j^*,\qquad
		\Phi^\omega(X,Y)=\sum_j v_j z w_j^* .
		\]
		For centered free elements $a_j$, Voiculescu's estimate
		\cite[Theorem~3.3]{CadilhacRicard} reads
		\begin{equation}\label{eq:voiculescu-two-input}
			\left\|\sum_j a_j\right\|_\infty
			\le\max_j\|a_j\|_\infty+2\left(\sum_j\|a_j\|_2^2\right)^{1/2}.
		\end{equation}
		The conjugated projections in each sum are free. For either sum,
		write $t=\tau(r_i)\le1/m$. After centering, each summand has
		operator norm at most $1$ and squared $L^2$ norm $t(1-t)$.
		Thus each sum has norm at most $mt+1+2\sqrt{mt(1-t)}\le4$.
		We obtain $\|X\|_\infty\le2\|x_p\|_\infty$ and
		$\|Y\|_\infty\le2\|y_p\|_\infty$.
		
		The summands of $W$ are centered and free. To check this,
		replace $v_j$ by $g_j=v_jw_j^*$ in the free generating set.
		Reduced words show that $(g_j,w_j)_j$ is again a free Haar
		family, free from the test algebra. The conjugates of the
		test algebra by the $w_j$ are free from each other and jointly
		free from the algebra of the $g_j$. Therefore the products
		$g_j(w_jvw_j^*)$ are centered and lie in free subalgebras.
		They have norm at most $1$ and squared $L^2$ norm at most
		$\delta^2/m$. Applying \eqref{eq:voiculescu-two-input} gives
		$\|W\|_\infty\le1+2\delta$. Moreover,
		\[
		\tau(W^*\Phi^\omega(X,Y))
		=\sum_{j,k}\tau(w_j v^*v_j^* v_k z w_k^*)
		=m\tau(|z|).
		\]
		The diagonal terms equal $\tau(v^*z)$, and each off-diagonal
		term has trace zero because it contains a Haar generator exactly
		once and its remaining factors are free from that generator.
		It follows that
		\[
		m\|\Phi^\omega(xp,py)\|_1
		\le(4+8\delta)\beta\|xp\|_\infty\|py\|_\infty.
		\]
		Apply this estimate to each $p_i$, sum, and use the uniform free
		compression estimates as $m\to\infty$. This proves the theorem.
	\end{proof}
	
	\begin{theorem}\label{thm:crossed-subgroup-two-input}
		Let a finite group $G$ act outerly on a $\mathrm{II}_1$ factor $\cQ$, let $H\le G$,
		and set $\cN=\cQ\rtimes H\subset\cM=\cQ\rtimes G$.
		For every bounded, separately normal bimodule bilinear map $\Phi$,
		\[
		\|\Phi\|_{2^2\to1}\le\sqrt{[G:H]}\,\|\Phi\|_{\infty^2\to1}.
		\]
		The bound is sharp in this family.
	\end{theorem}
	\begin{proof}
		Write $u_g$ for the crossed-product unitaries and $\alpha_g$ for
		the action. For $a\in\cQ$, the three module identities imply
		\[
		\Phi(u_g,u_k)a
		=\Phi(u_g,\alpha_k(a)u_k)
		=\Phi(u_g\alpha_k(a),u_k)
		=\alpha_{gk}(a)\Phi(u_g,u_k).
		\]
		Thus $u_{gk}^*\Phi(u_g,u_k)$ lies in
		$\cQ'\cap(\cQ\rtimes G)=\bC1$, and
		$\Phi(u_g,u_k)=c(g,k)u_{gk}$ for a scalar $c(g,k)$.
		Unitary testing gives $|c(g,k)|\le\|\Phi\|_{\infty^2\to1}$.
		
		Take representatives $s$ of the cosets $Hs$ and write
		$\displaystyle x=\sum_s a_su_s$ with $a_s\in\cN$. 
		Fourier orthogonality gives
		$\displaystyle \|x\|_2^2=\sum_s\|a_s\|_2^2$. 
		If $\displaystyle y=\sum_k y_ku_k$ with $y_k\in\cQ$, then balancing and outer bimodularity yield
		\[
		\Phi(x,y)=\sum_s a_su_s T_s(y),\qquad
		T_s(y)=\sum_k c(s,k)y_k u_k,\qquad
		\|T_s(y)\|_2\le\max_k|c(s,k)|\,\|y\|_2.
		\]
		H\"older's inequality on each summand and the Cauchy--Schwarz inequality over $s$ therefore give
		\[
		\|\Phi(x,y)\|_1\le
		\left(\sum_{s\in H\backslash G}\max_{k\in G}|c(s,k)|^2\right)^{1/2}
		\|x\|_2\|y\|_2
		\le \sqrt{[G:H]}\,\|\Phi\|_{\infty^2\to1}\|x\|_2\|y\|_2.
		\]
		
		For a sharp example take $H=\{e\}$ and $G=(\bZ_2)^r$.
		Let $e_j$, $1\le j\le r$, denote the coordinate
		generators. 
		The commuting factors $(1+i u_{e_j})/\sqrt2$ are
		unitaries, so their product has the form $\displaystyle U=|G|^{-1/2}\sum_{g\in G}\eta_g u_g$ with $|\eta_g|=1$.
		Let $T(\sum_g a_gu_g)=\sum_g\eta_ga_gu_g$ and $\Phi(x,y)=T(xy)$. 
		The map $T$ is $\cQ$-bimodular and an $L^2$-isometry. 
		Hence for contraction inputs,
		$\|\Phi(x,y)\|_1\le\|T(xy)\|_2=\|xy\|_2\le1$, with equality at $(1,1)$. 
		The projection $\displaystyle p=|G|^{-1}\sum_g u_g$ has trace $1/|G|$.
		Taking $x=y=\sqrt {|G|}\,p$ yields $\Phi(x,y)=\sqrt {|G|}\,U$, while both input $L^2$ norms are one.
	\end{proof}
	
	\begin{proposition}\label{thm:fixed-subgroup-two-input}
		Let a finite group $G$ act outerly on a $\mathrm{II}_1$ factor $\cQ$, let
		$H\le G$, and put $X=G/H$ and $D=[G:H]$. If an abelian subgroup
		of $G$ acts transitively on $X$, then
		\[
		\KNM{\cQ^G}{\cQ^H}_2\le\sqrt {[G:H]}.
		\]
		Under this hypothesis equality holds if there is a $G$-invariant
		complex Hadamard matrix on $X$. In particular equality holds for
		the dihedral pairs $G=C_{[G:H]}\rtimes C_2$, $H=C_2$, for $[G:H]\ge3$.
	\end{proposition}
	\begin{proof}
		Set $n=|G|$ and $\cN=\cQ^G$. 
		Since $[\cQ:\cN]=n$, the downward basic construction \cite[\S3.1]{Jones1983}
		gives a projection $p\in\cQ$ with $\bE_{\cN}(p)=1/n$.
		The projections $p_g=\alpha_g(p)$ sum to $1$ by the group averaging
		formula for $\bE_{\cN}$. They are pairwise orthogonal: compressing
		the sum by any $p_g$ shows $p_gp_hp_g=0$ for $h\ne g$.
		
		First consider a bimodule bilinear map $B:\cQ\times\cQ\to\cQ$ over $\cN$. 
		For $a,b\in G$, $B(p_a,p_b)$ commutes with $\displaystyle \sum_g\alpha_g(z)$ for every $z\in p\cQ p$, by moving this
		element from the left of the first input to the right of the
		second via the module identities. 
		The commutant in $\cQ$ of this diagonal algebra is $\operatorname{span}\{p_g:g\in G\}$.
		On each diagonal corner the commutation relation forces a scalar, since $p\cQ p$ is a factor. 
		For an off-diagonal corner $p_s\cQ p_t$, conjugate by $\alpha_{s^{-1}}$. 
		A nonzero intertwiner would have polar part $v$ with $v^*v=\theta(p)$, $vv^*=p$, and $v\theta(z)=zv$ on $p\cQ p$, where $\theta=\alpha_{s^{-1}t}$.
		Choose partial isometries $r_j$ with common initial projection $p$ and final projections forming a partition of $1$. 
		Then $\displaystyle V=\sum_j r_jv\theta(r_j)^*$ is unitary and $V\theta(x)=xV$ for all $x\in\cQ$, as follows by writing $\displaystyle x=\sum_{i,j}r_i(r_i^*xr_j)r_j^*$. 
		This contradicts outerness when $s\ne t$. 
		Thus
		\[
		B(p_a,p_b)=\sum_k c_{k;a,b}p_k .
		\]
		
		For $h\in G$, choose a partial isometry $w_h\in\cQ$ from $p$ to $p_h$. 
		Then $\displaystyle V_h=\sum_g\alpha_g(w_h)$ is a unitary of $\cN$ with $V_hp_gV_h^*=p_{gh}$. 
		Covariance of $B$ under conjugation by $V_h$ implies $c_{kh;ah,bh}=c_{k;a,b}$. 
		If $d_{s,t}=c_{e;s^{-1},t^{-1}}$, the map
		$\displaystyle (x,y)\mapsto\sum_{s,t}d_{s,t}\alpha_s(x)\alpha_t(y)$
		agrees with $B$ on each $(p_a,p_b)$: its $p_k$ coefficient is
		$d_{ka^{-1},kb^{-1}}=c_{k;a,b}$. 
		The maps therefore agree on all of $\cQ\times\cQ$, since
		\[
		x=\sum_a n\,\bE_{\cN}(xp_a)p_a,\qquad
		y=\sum_b p_b n\,\bE_{\cN}(p_b y),
		\]
		as can be checked by expanding the group averages and using orthogonality.
		
		Given $\Phi$ on $\cQ^H\times\cQ^H$, extend it to such a $B$ by applying $\bE_{\cQ^H}$ to both inputs. 
		Identify $x\in\cQ^H$ with the equivariant tuple $x_{tH}=\alpha_t(x)$, where $X$ has normalized counting measure.
		This map preserves $\|\cdot\|_p$ for $1\le p\le\infty$. 
		The formula for $B$ shows that at $k=tH$, $\Phi(x,y)$ is represented by
		\[
		b(x,y)_k=\sum_{i,j\in X} c_{kij}x_i y_j,\qquad
		c_{kij}=\sum_{\substack{taH=i\\tbH=j}}d_{a,b}.
		\]
		The coefficients are independent of the representative $t$.
		Indeed $\displaystyle f\mapsto\sum_g f(g^{-1}H)p_g$ embeds $\ell^\infty(X)$ in $\cQ^H$, and the $p_e$ coefficient of its $tH$ coordinate is $f(tH)$. 
		Taking the two inputs to be embedded coordinate indicators reads off each $c_{kij}$ at $p_e$, independently of the representative. 
		Replacing $t$ by $gt$ in the coefficient formula gives $c_{gk,gi,gj}=c_{kij}$. 
		The model $b$ therefore commutes with translations by $G$. 
		Its scalar restriction is realized on the embedded $\ell^\infty(X)$: the $p_e$ coefficient at $k$ is $b(f,f')_k$, and equivariance determines all the other diagonal coefficients from these.
		
		An abelian transitive subgroup acts regularly after dividing
		by its kernel, since its stabilizers are identical. 
		Identify $X$ with this quotient group and let $v_\chi$ denote its characters.
		The $v_\chi$ are an orthonormal unitary basis for normalized
		counting measure. 
		Translation covariance gives
		$b(v_\chi,v_\psi)=m(\chi,\psi)v_{\chi\psi}$, since the
		corresponding character eigenspace is one-dimensional.
		The scalar
		embedding gives $|m(\chi,\psi)|\le\|\Phi\|_{\infty^2\to1}$.
		For $\cQ$-valued tuples expand $\displaystyle x=\sum_\chi x_\chi v_\chi$
		and $\displaystyle y=\sum_\psi y_\psi v_\psi$.
		With
		$T_\chi(y)=\sum_\psi m(\chi,\psi)y_\psi v_\psi$, the output is
		$\displaystyle \sum_\chi x_\chi v_\chi T_\chi(y)$.
		Fourier orthogonality gives
		$\|T_\chi(y)\|_2\le\|\Phi\|_{\infty^2\to1}\|y\|_2$ and
		$\displaystyle \sum_\chi\|x_\chi\|_2^2=\|x\|_2^2$.
		Apply H\"older's inequality and sum over the $[G:H]$ characters to prove
		the asserted bound.
		
		Finally let $K$ be a $G$-invariant complex Hadamard matrix on
		$X$, meaning $KK^*=[G:H]I$ and $|K_{ij}|=1$. 
		Applying $[G:H]^{-1/2}K$ to the equivariant coordinates defines a normal
		$\cN$-bimodule map $T:\cQ^H\to\cQ^H$. 
		The constant tuple
		subspace is preserved and the scalar constant line is multiplied
		by a phase. 
		Thus $T(1)$ is scalar unitary. 
		Also $T$ is an $L^2$-isometry. 
		It follows that $\Phi(x,y)=T(xy)$ has $\|\Phi\|_{\infty^2\to1}=1$. 
		Let $q_i$ be an embedded scalar coordinate projection. 
		It has trace $1/[G:H]$ and $\sqrt {[G:H]}\,T(q_i)$ is unitary, since all entries in the $i$th column of $K$ have modulus one. 
		Testing $x=y=\sqrt {[G:H]}\,q_i$ gives $\|\Phi(x,y)\|_1=\sqrt {[G:H]}$ and $\|x\|_2=\|y\|_2=1$.
		
		For the dihedral pairs, $X=C_{[G:H]}$ has rotations and inversion.
		With integer representatives $t$, define the even, periodic phase
		\[
		f(t)=
		\begin{cases}
			\exp(\pi i t^2/[G:H]),&[G:H]\text{ even},\\
			\exp(2\pi i t^2/[G:H]),&[G:H]\text{ odd}.
		\end{cases}
		\]
		For a shift $r\ne0$ modulo $[G:H]$, the sum $\displaystyle \sum_t f(t+r)\overline{f(t)}$ is a geometric sum with ratio
		$e^{2\pi i r/[G:H]}$ in even order and $e^{4\pi i r/[G:H]}$ in odd
		order, so it is zero. 
		Hence $K_{ij}=f(j-i)$ is a $G$-invariant complex Hadamard matrix, and the rotations give the abelian transitive subgroup. 
		The same pairs attain the crossed-product bound in Theorem~\ref{thm:crossed-subgroup-two-input}.
		Indeed, the multiplier $T(u_{ah})=f(a)u_{ah}$, $a\in C_{[G:H]}$ and
		$h\in H$, is $\cQ\rtimes H$-bimodular because $f(a)=f(-a)$.
		Its Fourier symbol has modulus one, and
		$[G:H]^{-1/2}\sum_a f(a)u_a$ is unitary by the same autocorrelation
		identity.
		The preceding crossed-product test with $\displaystyle p=[G:H]^{-1}\sum_a u_a$ proves equality.
	\end{proof}
	
	\begin{remark}
		The abelian transitivity hypothesis by itself allows strict inequality. 
		For an outer $S_{m}$-action and subgroup $S_{m-1}$, consider
		the natural action on $X=\{1,\ldots,m\}$. 
		For $m\ge3$ one has
		\[
		\KNM{\cQ^{S_{m}}}{\cQ^{S_{m-1}}}_2\le 2+\sqrt2 .
		\]
		The cyclic subgroup generated by an $m$-cycle is transitive, so this gives strict inequality in Proposition~\ref{thm:fixed-subgroup-two-input} for $m \ge12$.
		
		To prove the estimate, write $E$ for averaging over $X$, viewed
		as a projection onto the constant tuples.
		The orbit of a triple $(k,i,j)$ under $S_m$ is determined by which of its entries coincide. 
		The five tensors representing
		$xy$, $xE(y)$, $E(x)y$, $E(xy)$, and $E(x)E(y)$ span these
		five orbit patterns.
		On the cyclic character basis, such a map has a symbol of the form
		\[
		\mathbf m(\chi,\psi)=\alpha_0+\alpha_1\,1_{\{\psi=1\}}+\alpha_2\,1_{\{\chi=1\}}
		+\alpha_3 \,1_{\{\chi\psi=1\}}+\alpha_4 \,1_{\{\chi=\psi=1\}}.
		\]
		Each of its five values has modulus at most $\beta=\|\Phi\|_{\infty^2\to1}$ by the scalar contraction tests above. 
		Let $\lambda_{11}$ be the value $\alpha_0$ for two nontrivial
		characters with nontrivial product, and let $\lambda_{\rm inv}$ be the value $\alpha_0+\alpha_3$ for two nontrivial characters with trivial product. 
		Let $\lambda_{10}=\alpha_0+\alpha_1$, $\lambda_{01}=\alpha_0+\alpha_2$, and
		$\lambda_{00}=\alpha_0+\alpha_1+\alpha_2+\alpha_3+\alpha_4$ be the values with a trivial input character, using $0$ for the trivial character.
		With $x=x_0+x_1$ and $y=y_0+y_1$, where $x_0=E(x)$ and $y_0=E(y)$, the formula becomes
		\[
		b(x,y)=\lambda_{11}x_1y_1+
		(\lambda_{\rm inv}-\lambda_{11})E(x_1y_1)
		+\lambda_{10}x_1y_0+\lambda_{01}x_0y_1+\lambda_{00}x_0y_0 .
		\]
		H\"older's inequality and the $L^1$ contractivity of $E$ give
		\[
		\|b(x,y)\|_1\le\beta
		\begin{pmatrix}\|x_1\|_2 & \|x_0\|_2\end{pmatrix}
		\begin{pmatrix}3&1\\1&1\end{pmatrix}
		\begin{pmatrix}\|y_1\|_2\\\|y_0\|_2\end{pmatrix}.
		\]
		The matrix has norm $2+\sqrt2$, and the decompositions are
		orthogonal in $L^2$. The tuple embedding thus proves the
		claimed estimate for $\Phi$.
	\end{remark}
	
	\subsection{Matrix algebras}
	
	Every unital inclusion of full matrix factors is unitarily equivalent to $\cN=M_m\otimes 1_n\subset\cM=M_m\otimes M_n$. 
	Its relative
	commutant is $1_m\otimes M_n$, so the only irreducible case is
	$n=1$. 
	With normalized traces, write $\mathbf K_2(m,n)$ for the supremum
	of $\|\Phi\|_{2^2\to1}/\|\Phi\|_{\infty^2\to1}$ over nonzero maps satisfying the three bimodule identities.
	
	\begin{theorem}\label{thm:matrix-two-input}
		Suppose that $s=\min(m,n)$. We have that
		\[
		n\sqrt{s}\le \mathbf K_2(m,n)\le n\min\{s,\sqrt{2m}\}.
		\]
		In particular, $\mathbf K_2(1,n)=n$ and $\mathbf K_2(m,1)=1$.
		For $m\le n$, the upper and lower bounds differ by a factor of at
		most $\sqrt2$. 
		The inequality
		$\|\Phi\|_{\infty^2\to1}\le\|\Phi\|_{2^2\to1}$ has optimal coefficient one, attained by multiplication.
	\end{theorem}
	
	\begin{proof}
		Use $\|\cdot\|_{S_p}$ for the Schatten norm with ordinary trace.
		The two factors $(mn)^{-1/2}$ in the input $L^2$ norms cancel the
		factor $(mn)^{-1}$ in the output $L^1$ norm. Thus
		\begin{equation}\label{eq:matrix-two-input-unnormalized}
			\begin{aligned}
				\|\Phi\|_{2^2\to1}
				&=\sup_{\|X\|_{S_2},\,\|Y\|_{S_2}\le1}\|\Phi(X,Y)\|_{S_1},\\
				mn\,\|\Phi\|_{\infty^2\to1}
				&=\sup_{\|X\|_\infty,\,\|Y\|_\infty\le1}\|\Phi(X,Y)\|_{S_1}.
			\end{aligned}
		\end{equation}
		
		\emph{The dimension bound.}
		The module identities give a unique bilinear map
		$f:M_n\times M_n\to M_n$ such that
		\begin{equation}\label{eq:matrix-two-input-amplification}
			\Phi(X,Y)_{ij}=\sum_{k=1}^m f(X_{ik},Y_{kj}).
		\end{equation}
		Indeed, the outer identities place
		$\Phi(e_{ik}\otimes x,e_{\ell j}\otimes y)$ in the $(i,j)$ block
		and identify that block with the corresponding $(1,1)$ block.
		Balancing with $e_{kk}\otimes1$ gives zero for $k\ne\ell$, and
		balancing with $e_{1k}\otimes1$ gives
		$\Phi(e_{1k}\otimes x,e_{k1}\otimes y)
		=\Phi(e_{11}\otimes x,e_{11}\otimes y)$.
		For a fixed $f$ denote the two norms at size $m$ by
		$\beta^{(m)}=\|\Phi\|_{\infty^2\to1}$ and
		$\gamma^{(m)}=\|\Phi\|_{2^2\to1}$. 
		Testing with
		$I_m\otimes x,I_m\otimes y$ gives $\beta^{(1)}\le\beta^{(m)}$.
		
		Let $\cH=M_n$ with its ordinary Hilbert--Schmidt inner product.
		The rule $L_f(|x\rangle\langle y^*|)=f(x,y)$ defines a linear
		map $L_f:S_1(\cH)\to S_1(\bC^n)$. 
		The left argument is bilinear
		in $(x,y)$, since both adjoint on $\cH$ and the bra operation
		are conjugate-linear.
		Taking $x,y$ from the matrix unit basis
		of $\cH$ gives a basis of $S_1(\cH)$ and determines $L_f$.
		For $\displaystyle u_k=\sum_i e_i\otimes X_{ik}$ and $\displaystyle v_k=\sum_j e_j\otimes Y_{kj}^*$, we have
		\[
		\Phi(X,Y)=(\id_m\otimes L_f)
		\left(\sum_k |u_k\rangle\langle v_k|\right),\qquad
		\left\|\sum_k|u_k\rangle\langle v_k|\right\|_{S_1}
		\le\sum_k\|u_k\|\,\|v_k\|\le\|X\|_{S_2}\|Y\|_{S_2}.
		\]
		Conversely, a rank-one input to $\id_m\otimes L_f$ is realized
		by filling just the first block column of $X$ and the first block row of $Y$. 
		Rank-one operators suffice to compute the induced trace norm by singular value decomposition. 
		Trace duality then gives
		\begin{equation}\label{eq:matrix-two-input-diamond}
			\gamma^{(m)}=\|\id_m\otimes L_f\|_{S_1\to S_1}
			=\|\id_m\otimes L_f^*\|_{\infty\to\infty}.
		\end{equation}
		
		Expand a contraction $\displaystyle W=\sum_{\ell=1}^{n^2}A_\ell\otimes U_\ell$
		in a normalized-trace orthonormal basis of Weyl unitaries on
		$\bC^n$. 
		Then
		$\displaystyle \sum_\ell A_\ell A_\ell^*=(\id_m\otimes\tr_n)(WW^*)\le I_m$.
		Factoring the amplified image as the product of the row
		$(A_\ell\otimes 1)_\ell$ and the column
		$(1\otimes L_f^*(U_\ell))_\ell$ yields
		\[
		\|(\id_m\otimes L_f^*)(W)\|_\infty
		\le \left\|\sum_\ell L_f^*(U_\ell)^*L_f^*(U_\ell)\right\|_\infty^{1/2}
		\le n\,\gamma^{(1)}.
		\]
		Each block of the image also has norm at most $\gamma^{(1)}$.
		Its action on $(\xi_j)_{j=1}^m$ therefore has squared norm at most
		$m (\gamma^{(1)})^2(\sum_j\|\xi_j\|)^2\le
		m^2(\gamma^{(1)})^2\sum_j\|\xi_j\|^2$, which gives the alternative bound $m\,\gamma^{(1)}$.
		At size $m=1$, \eqref{eq:matrix-two-input-unnormalized} and
		$\|x\|_\infty\le\|x\|_{S_2}$ give
		$\gamma^{(1)}\le n\beta^{(1)}$. 
		Consequently $\gamma^{(m)}\le n s\beta^{(m)}$.
		
		\emph{The packing bound.}
		The rank-one reduction above lets us take $X=X(e_{11}\otimes1)$
		and $Y=(e_{11}\otimes1)Y$, with ordinary Hilbert--Schmidt norms
		one. 
		Put $Z=\Phi(X,Y)$. 
		Choose an orthonormal basis of $\bC^m$
		in which the diagonal entries of $(\id_m\otimes\operatorname{Tr}_n)(XX^*)$ are $1/m$. 
		For example, apply the Fourier matrix to an eigenbasis
		of this positive matrix of trace one. 
		Take the $m$ powers $U_j$
		of the diagonal clock in that basis, as unitaries in $\cN$.
		The average of $U_jXX^*U_j^*$ retains only the diagonal blocks.
		Each is positive with ordinary trace $1/m$, hence is bounded by
		$I_n/m$. 
		This gives $\displaystyle \sum_j U_jXX^*U_j^*\le I$. 
		Likewise choose
		$V_j$ so that $\displaystyle \sum_j V_jY^*YV_j^*\le I$.
		
		Suppressing the identity in the $M_n$ leg, take scalar phases
		$\varepsilon_j$ of modulus one and set
		\[
		X'=\sum_j\varepsilon_j U_j X e_{1j},\qquad
		Y'=\sum_j e_{j1}YV_j^*.
		\]
		The products $e_{1j}e_{k1}=\delta_{jk}e_{11}$ give
		$X'X'^*=\sum_j U_jXX^*U_j^*$ and
		$Y'^*Y'=\sum_j V_jY^*YV_j^*$, so both inputs are contractions.
		The same matrix-unit identity and balancing give
		$\displaystyle \Phi(X',Y')=\sum_j\varepsilon_j U_j ZV_j^*$.
		
		For independent Steinhaus phases, the $p=1$ noncommutative Khintchine inequality \cite[Theorem~2.1 and Remark~2.8]{HMKhintchine}
		implies
		\[
		\bE\left\|\sum_j \varepsilon_j Z_j\right\|_{S_1}
		\ge \frac1{\sqrt2}\left(\sum_j\|Z_j\|_{S_1}^2\right)^{1/2}.
		\]
		To see the relation with its row-plus-column infimum, pair any
		block column $(A_j)_j$ with the block row of polar adjoints $(\lambda_j v_j^*)_j$, where $A_j=v_j|A_j|$ and $\displaystyle \sum_j\lambda_j^2=1$ with $\lambda_j\ge0$. 
		The row has norm at most one. 
		Optimizing over $\lambda_j$ gives
		$\displaystyle \left\| (\sum_j A_j^*A_j)^{1/2}\right\|_{S_1}\ge
		\left(\sum_j\|A_j\|_{S_1}^2\right)^{1/2}$. 
		The same bound holds for rows by taking adjoints, and the triangle inequality in $\ell^2$ gives the displayed consequence for every decomposition $Z_j=A_j+B_j$.
		Here $Z_j=U_j ZV_j^*$ has trace norm $\|Z\|_{S_1}$. 
		Hence $mn\beta^{(m)}\ge\sqrt{m/2}\|Z\|_{S_1}$. 
		Taking the supremum over $X,Y$ as above proves $\gamma^{(m)}\le n\sqrt{2m}\,\beta^{(m)}$.
		
		\emph{The partial-transpose lower bound.}
		Let $p=I_m\otimes P_s$, where $P_s$ projects onto the first $s$
		basis vectors, and let $q=I_m\otimes e_{11}$. 
		Define
		\begin{equation}\label{eq:matrix-two-input-pt}
			\Phi(X,Y)=(pXqYp)^\Gamma,\qquad \Gamma=\id_m\otimes t_n,
		\end{equation}
		where $t_n$ is matrix transpose. 
		The projections $p,q$ commute
		with $\cN$ and $\Gamma$ is an $\cN$-bimodule map, so the
		three identities hold. 
		For contraction inputs, $pXqYp$ is a
		contraction of rank at most $m$. 
		Partial transpose preserves
		the Hilbert--Schmidt norm and the corner $p\cM p$, whose
		underlying dimension is $ms$. 
		Thus
		\[
		\|\Phi\|_{\infty^2\to1}\le\frac{\sqrt{ms}\sqrt m}{mn}
		=\frac{\sqrt s}{n}.
		\]
		The induced trace norm of partial transpose on this corner
		is at most $s$, by trace duality and the Weyl row estimate
		above with transpose on $M_s$. 
		The rank-one projection onto
		$\displaystyle \xi=s^{-1/2}\sum_{i=1}^s e_i\otimes e_i$ transposes to
		the flip divided by $s$ on $\bC^s\otimes\bC^s$, so its
		trace norm is $s$. 
		For a unit $h$ in the range of $q$, let
		$X=|\xi\rangle\langle h|$ and $Y=X^*$. 
		These inputs have ordinary
		Hilbert--Schmidt norms one and realize that projection as
		$pXqYp$. Since $\|pXqYp\|_{S_1}\le\|X\|_{S_2}\|Y\|_{S_2}$,
		we get $\|\Phi\|_{2^2\to1}=s$ and a ratio at least $n\sqrt{s}$.
		Finally multiplication gives equality in the reverse estimate,
		by H\"older's inequality and the test at $(1,1)$.
	\end{proof}
	
	If $s$ divides $m$, take Weyl unitaries $D_0,S_0\in M_m$ with
	$D_0S_0=\omega S_0D_0$, where $\omega=e^{2\pi i/s}$, by
	repeating the $s$-dimensional representation. 
	The contractions
	\[
	X=s^{-1/2}\sum_{i=0}^{s-1}D_0^i\otimes e_{i+1,1},\qquad
	Y=s^{-1/2}\sum_{j=0}^{s-1}S_0^j\otimes e_{1,j+1}
	\]
	have $X^*X=q$ and $YY^*=q$. 
	The $(j+1,i+1)$ block of
	$(pXqYp)^\Gamma$, taken in the $M_n$ leg, is $s^{-1}D_0^i S_0^j$. 
	Its product with its adjoint therefore has $(j+1,k+1)$ block
	$\displaystyle s^{-2}\sum_i\omega^{i(j-k)}S_0^{j-k}=s^{-1}\delta_{jk}I_m$
	inside the $p$ corner. 
	Hence $|\Phi(X,Y)^*|=p/\sqrt s$, and
	$\|\Phi\|_{\infty^2\to1}=\sqrt{s}/n$. 
	In particular, for $m=n$ this map has ratio $n^{3/2}$ exactly.
	
	\section{Multilinear Grothendieck Inequalities}
	
	Christensen--Pop--Sinclair--Smith \cite[Theorem~3.3]{CPSS} consider a separately
	normal, $\cN$-multimodular $k$-linear map $\Phi:\cM^k\to B(\cH)$, where
	$\cN\subset \cM$ is an irreducible inclusion of $\mathrm{II}_1$ factors and
	$\cN$ is hyperfinite. With $\tr_m$ the normalized trace on $M_m(\cM)$ and
	\[
	\rho_m(X)=\bigl(\|X\|_{\op}^2+m\,\tr_m(X^*X)\bigr)^{1/2},
	\]
	their matrix-product amplification satisfies
	\begin{equation}\label{eq:CPSS}
		\|\Phi^{(m)}(X_1,\ldots,X_k)\|_{\op}
		\le 2^{k/2}\|\Phi\|\prod_{j=1}^k\rho_m(X_j).
	\end{equation}
	Here $\|\Phi\|$ uses operator norms on the domain and range. 
	We study the $L^1$ output norm with normalized $L^p$ input norms. 
	Throughout this section $n\ge2$, and
	$\Phi:\cM^{\times n}\to\cM$ is bounded, separately normal, and
	$\cN$-multimodular.
	
	\begin{proposition}\label{prop:multilinear-reverse}
		Suppose that $\cN\subset\cM$ is an irreducible inclusion of finite factors.
		For every $\boldsymbol p$ with $1\le p_j\le\infty$,
		\[
		\|\Phi\|_{\infty^n\to1}\le \|\Phi\|_{\boldsymbol p\to1}.
		\]
		The coefficient $1$ is optimal whenever $\displaystyle \sum_j 1/p_j\le1$.
		If $\displaystyle \sum_j 1/p_j>1$, a finite forward coefficient can fail even for multiplication on the identity inclusion of a $\mathrm{II}_1$ factor.
	\end{proposition}
	\begin{proof}
		Each operator-norm unit ball is contained in the corresponding $L^{p_j}$ unit ball because the trace is normalized. 
		When $\displaystyle s=\sum_j 1/p_j\le1$, H\"older's inequality and $\|z\|_1\le \|z\|_{1/s}$ give
		\[
		\|x_1\cdots x_n\|_1\le\prod_{j=1}^n\|x_j\|_{p_j}.
		\]
		The tuple $(1,\ldots,1)$ shows that both norms of multiplication are $1$.
		When $s>1$, take a projection $e$ of trace $t>0$ in a $\mathrm{II}_1$ factor
		and put $x_j=t^{-1/p_j}e$. 
		Each $x_j$ has $L^{p_j}$ norm one, whereas
		$\|x_1\cdots x_n\|_1=t^{1-s}$, which tends to infinity as $t\downarrow0$.
	\end{proof}
	
	On the cycle with vertices $0,1,\ldots,n$ in this cyclic order, denote
	shortest-path distance by $d(i,j)$.
	
	\begin{theorem}\label{thm:multilinear-fixed-trace}
		Suppose that $\cN\subset\cM$ is an irreducible inclusion of finite factors with Jones index $\delta^2$. 
		Let $1\le p_i\le\infty$,
		$\displaystyle \frac{1}{p}=\sum_{i=1}^n\frac{1}{p_i}\le1$, and $p_0=p'$ with the usual endpoint convention.
		Put $\displaystyle r_j=\sum_{i=0}^n\frac{1}{p_i}d(i,j)$. 
		Then
		\begin{equation}\label{eq:multilinear-fixed-trace}
			\|\Phi\|_{\boldsymbol p\to1}\le
			\min\left\{\delta^{2\left\lceil\frac{n+1}{2}\right\rceil/p},\
			2e^2(1+\sqrt2)\sum_{0\le j<k\le n}\delta^{r_j+r_k}\right\}
			\|\Phi\|_{\infty^n\to1}.
		\end{equation}
		In particular, with $\kappa_n=(2\lfloor(n+1)^2/4\rfloor-1)/n$ one has
		\begin{equation}\label{eq:multilinear-symmetric}
			\|\Phi\|_{n^n\to1}\le
			\min\left\{\delta^{2\left\lceil\frac{n+1}{2}\right\rceil},\
			2e^2(1+\sqrt2)\binom{n+1}{2}\delta^{\kappa_n}\right\}
			\|\Phi\|_{\infty^n\to1}.
		\end{equation}
		For $n=2$ and $p_1=p_2=2$, Theorem~\ref{thm:general-two-input} gives a
		stronger bound.
	\end{theorem}
	
	\begin{proof}
		Put $D=\delta^2$, $L_n=\lceil(n+1)/2\rceil$,
		$c=2e^2(1+\sqrt2)$, and $A=\|\Phi\|_{\infty^n\to1}$. It is convenient to
		put the trace variable first and use the cyclic order $0,1,\ldots,n$:
		\[
		F(x_0,x_1,\ldots,x_n)=\tau\bigl(x_0\Phi(x_1,\ldots,x_n)\bigr).
		\]
		Trace duality shows that $\|F\|=A$ on the operator-norm unit balls.
		The identities for $\Phi$ give balance between $x_i$ and $x_{i+1}$:
		a factor $a\in\cN$ can be moved from the right of $x_i$ to the left of
		$x_{i+1}$. This holds between $x_0$ and $x_1$ by the left module identity
		and between $x_n$ and $x_0$ by the right module identity and the tracial
		property. When $D=1$, multimodularity gives
		$\Phi(x_1,\ldots,x_n)=x_1\cdots x_n\Phi(1,\ldots,1)$ and
		$\Phi(1,\ldots,1)$ is central. 
		H\"older's inequality therefore proves the result.
		When $D>1$ the factors are of type $\mathrm{II}_1$, since a unital inclusion
		of matrix factors with trivial relative commutant must be equality.
		
		\smallskip
		\noindent\textit{Support bounds at every vertex.}
		Let $t_i=\tau(\cR(x_i))$, the common trace of the left and right supports
		of $x_i$. A zero input gives $F=0$, so suppose all $t_i>0$.
		For an edge $e$ with endpoints $e_-,e_+$ in cyclic order write
		$d(i,e)=\min\{d(i,e_-),d(i,e_+)\}$, and set
		\[
		\alpha_e=\min_{0\le i\le n} D^{1+d(i,e)}t_i.
		\]
		Process the edges in increasing order of $\alpha_e$. At each edge choose a
		projection $q_e\in\cN$ of trace at most $\alpha_e$ and multiply the input
		at $e_-$ on the right and the input at $e_+$ on the left by $q_e$.
		To find $q_e$, fix a vertex $i$ attaining the minimum in $\alpha_e$.
		If $i$ is incident to $e$, the support on the relevant side of the
		current $i$th input has trace at most $t_i$. The support estimate
		\eqref{eq:support} covers it by a projection of $\cN$ with trace
		at most $Dt_i$.
		
		If $i$ is not incident to $e$, let $e'$ be the preceding edge on a
		shortest path from $i$ to an endpoint of $e$. Then
		$d(i,e')=d(i,e)-1$ and
		\[
		\alpha_{e'}\le D^{1+d(i,e')}t_i=\alpha_e/D<\alpha_e.
		\]
		The edge $e'$ has already been processed. At the shared vertex, one
		support lies under $q_{e'}$, so the opposite support has trace at most
		$\tau(q_{e'})$. By \eqref{eq:support}, this opposite support can be
		covered on the $e$ side by a projection of trace at most
		$D\tau(q_{e'})\le \alpha_e$.
		
		In both cases $q_e$ acts as an identity on the input whose support it
		covers. Balance then inserts $q_e$ on the other input without changing
		$F$, and multiplication on the opposite side preserves all earlier support
		relations. Write $\tilde x_j$ for the final inputs. They satisfy
		$\|\tilde x_j\|_\infty\le\|x_j\|_\infty$ and
		\begin{equation}\label{eq:multilinear-vertex-support}
			\tau(\cR(\tilde x_j))\le \beta_j
			:=\min_{0\le i\le n}D^{d(i,j)}t_i.
		\end{equation}
		Indeed, for $i=j$ use the original support bound. For $i\ne j$, the edge
		at $j$ toward $i$ on a shortest path has $d(i,e)=d(i,j)-1$, so its
		projection bounds the support of $\tilde x_j$ by $D^{d(i,j)}t_i$.
		
		\smallskip
		\noindent\textit{Packing and the first coefficient.}
		Fix a distinguished input of support trace $t$. Since $1+d(i,e)\le L_n$
		for every vertex $i$ and edge $e$, all the $q_e$ have trace at most
		$D^{L_n}t$. Write $q_j$ for the projection at the edge from $j$ to $j+1$,
		with indices read modulo $n+1$. For an integer $k$ with $kD^{L_n}t\le1$,
		choose $u_{j,1},\ldots,u_{j,k}\in\cN$ with initial projection $q_j$ and
		orthogonal final projections for each $j$. Set
		\[
		y_j=\sum_{\ell=1}^k u_{j-1,\ell}\tilde x_j u_{j,\ell}^* .
		\]
		For fixed $j$, the summands have pairwise orthogonal initial and final
		supports. Thus $\|y_j\|_\infty\le\|\tilde x_j\|_\infty$. 
		In the expansion of $F(y_0,\ldots,y_n)$, balance produces a product
		$u_{j,\ell}^*u_{j,\ell'}$ at each edge. 
		This is zero for $\ell\ne\ell'$ and equals $q_j$ for $\ell=\ell'$. Consequently all labels agree in every surviving term, and each surviving term is $F(x_0,\ldots,x_n)$. 
		Hence
		\[
		k\,|F(x_0,\ldots,x_n)|\le A\prod_{j=0}^n\|x_j\|_\infty.
		\]
		Suppose the distinguished input is $v e$, where $v$ is a partial isometry, and  $0\neq e\leq v^*v$ is a projection. 
		Split $e$ in $\cM$ into $m$ orthogonal subprojections of trace $\tau(e)/m$. Apply the last bound to each summand in the distinguished slot with $k_m=\lfloor m/(D^{L_n}\tau(e))\rfloor$. As $m\to\infty$, $m/k_m\to D^{L_n}\tau(e)$. 
		Together with the original norm
		bound this gives, for any such distinguished input in slot $i$,
		\[
		|F(x_0,\ldots,ve,\ldots,x_n)|
		\le A\min\{1,D^{L_n}\tau(e)\}\prod_{j\ne i}\|x_j\|_\infty .
		\]
		Writing a general $x_i=v_i|x_i|$ as the integral of $v_i1_{(s,\infty)}(|x_i|)$ gives the one-slot estimate
		\begin{equation}\label{eq:multilinear-one-slot}
			|F(x_0,\ldots,x_n)|\le
			D^{L_n} A \|x_i\|_1\prod_{j\ne i}\|x_j\|_\infty
			\qquad(0\le i\le n).
		\end{equation}
		The integral can first be taken for finite-spectrum $|x_i|$, and uniform
		spectral approximation gives the general case.
		
		To obtain the mixed-exponent bound, decompose $x_1,\ldots,x_n$
		by the same layer integral, keeping $x_0$ fixed. Write
		$t_i(s_i)=\tau(1_{(s_i,\infty)}(|x_i|))$ and let $\mu_i$ be the decreasing
		singular-value function of $x_i$, extended by zero beyond $1$. The estimate
		on each layer uses $\min\{1,D^{L_n}\min_{1\le i\le n}t_i(s_i)\}$. Integrating
		this minimum by introducing a variable $r\in(0,1)$ gives
		\[
		\int_{(0,\infty)^n}\min\left\{1,D^{L_n}\min_{1\le i\le n}t_i(s_i)\right\}\,
		d\boldsymbol s
		=\int_0^1\prod_{i=1}^n\mu_i(r/D^{L_n})\,dr
		\le D^{L_n/p}\prod_{i=1}^n\|x_i\|_{p_i}.
		\]
		For the equality, $r<D^{L_n}t_i(s_i)$ precisely on a set of $s_i$
		of length $\mu_i(r/D^{L_n})$. 
		For the inequality, use H\"older's inequality with exponents $p_1,\ldots,p_n,p_0$ on $(0,1)$, with the last factor equal to one. 
		For each finite $p_i$ the change of variable $r=D^{L_n}v$ gives a factor at most $D^{L_n/p_i}$, and the infinite-exponent factors are bounded by their essential suprema. 
		Taking the supremum over $\|x_0\|_\infty\le1$ yields the first coefficient in
		\eqref{eq:multilinear-fixed-trace}.
		
		\smallskip
		\noindent\textit{Free translates and the second coefficient.}
		By \eqref{eq:multilinear-one-slot}, $F$ extends to the tracial ultrapower
		by taking ultralimits along a free ultrafilter $\omega$. In fact, the
		difference produced by changing a representative is a sum of one-slot
		differences, bounded by \eqref{eq:multilinear-one-slot}. Each tends to
		zero because the representatives are uniformly bounded in operator norm
		and $\|z\|_1\le\|z\|_2$. 
		The resulting $F^\omega$ is cyclically balanced over
		$\cN^\omega$. 
		It has norm $A$ by lifting contractions through polar truncation and testing on constant sequences.
		
		By the relative freeness argument in the proof of Theorem~\ref{thm:general-two-input} (Popa's theorem \cite[Theorem~0.1(a)]{PopaFree14}), choose a free Haar family
		$\{u_{j,\ell}:0\le j\le n,\ 1\le\ell\le k\}$ in $\cN^\omega$ free from the algebra generated by the $\tilde x_j$. 
		For independent symmetric random signs $\epsilon_{j,\ell}$, set
		\[
		z_j=\sum_{\ell=1}^k\epsilon_{j-1,\ell}\epsilon_{j,\ell}\,
		u_{j-1,\ell}\tilde x_j u_{j,\ell}^*.
		\]
		In the expansion of $F^\omega(z_0,\ldots,z_n)$, each edge contributes
		$\epsilon_{j,\ell}\epsilon_{j,\ell'}$ for the labels chosen at its
		endpoints. 
		Its expectation is zero unless $\ell=\ell'$. 
		The unitaries cancel in every equal-label term by balance. 
		Therefore
		\[
		\bE F^\omega(z_0,\ldots,z_n)=kF(x_0,\ldots,x_n).
		\]
		
		For each fixed $j$, the unsigned summands are centered and free. To see
		this, write $w_\ell=u_{j,\ell}$ and
		$g_\ell=u_{j-1,\ell}u_{j,\ell}^*$. The replacement of
		$u_{j-1,\ell}$ by $g_\ell$ is a change of free generators, so the
		$g_\ell,w_\ell$ form a free Haar family free from the input algebra.
		Write each summand as $g_\ell(w_\ell\tilde x_jw_\ell^*)$. In a product
		of centered elements from the different conjugated input algebras,
		the consecutive $w_\ell$-letters cannot cancel across different labels.
		The same reduced-word argument works when centered words in the
		$g_\ell$ are inserted. Thus the conjugated input algebras are jointly
		free from each other and from the algebra of the $g_\ell$, which proves
		the claim. With $\tilde t_j=\tau(\cR(\tilde x_j))$, the squared $L^2$ norm
		of each summand is $\|\tilde x_j\|_2^2\le \tilde t_j\|x_j\|_\infty^2$.
		The free sum estimate \eqref{eq:voiculescu-two-input} gives, uniformly
		in the signs,
		\[
		k|F(x_0,\ldots,x_n)|
		\le A\prod_{j=0}^n\|x_j\|_\infty
		\prod_{j=0}^n(1+2\sqrt{k\tilde t_j}).
		\]
		
		Choose $j_0$ with $\tilde t_{j_0}$ maximal and set
		\[
		b=\sum_{j\ne j_0}\sqrt{\tilde t_j},\qquad
		P=\sqrt{\tilde t_{j_0}}\,b,\qquad
		Q=\sum_{\substack{j<k\\j,k\ne j_0}}\sqrt{\tilde t_j\tilde t_k}.
		\]
		If a compressed input vanishes, then $F=0$. Otherwise $b>0$, and
		maximality gives $b^2\le P+2Q$. If $b<1$, choose
		$k=\lfloor b^{-2}\rfloor$. The product outside slot $j_0$ is at most
		$\exp(2\sqrt{k}b)\le e^2$. Since $k\ge b^{-2}/2$, the factor at
		$j_0$, divided by $k$, is bounded by
		\[
		k^{-1}+2\sqrt{\tilde t_{j_0}/k}
		\le 2b^2+2\sqrt2\,P
		\le 2(1+\sqrt2)P+4Q.
		\]
		If $b\ge1$, the original norm estimate suffices because
		$1\le b^2\le P+2Q$. Using \eqref{eq:multilinear-vertex-support} in both
		cases gives
		\begin{equation}\label{eq:multilinear-pair-supports}
			|F(x_0,\ldots,x_n)|
			\le cA\prod_{j=0}^n\|x_j\|_\infty
			\sum_{0\le j<k\le n}\sqrt{\beta_j\beta_k}.
		\end{equation}
		
		Now decompose all the inputs into spectral layers, including $x_0$.
		For each $j$, put $\beta_j(\boldsymbol s)=\min_i D^{d(i,j)}t_i(s_i)$.
		Tonelli's theorem and H\"older's inequality with $\displaystyle \sum_{i=0}^n1/p_i=1$ give
		\[
		\int_{(0,\infty)^{n+1}}\beta_j(\boldsymbol s)\,d\boldsymbol s
		=\int_0^\infty\prod_{i=0}^n\mu_i(r/D^{d(i,j)})\,dr
		\le D^{r_j}\prod_{i=0}^n\|x_i\|_{p_i}.
		\]
		Here $r<\beta_j(\boldsymbol s)$ if and only if
		$r/D^{d(i,j)}<t_i(s_i)$ for every $i$. This proves the identity as in the
		first layer integral. A change of variable in the $i$th $L^{p_i}$ factor
		gives $D^{d(i,j)/p_i}$.
		The Cauchy--Schwarz inequality then bounds the integral of
		$\sqrt{\beta_j(\boldsymbol s)\beta_k(\boldsymbol s)}$ by
		$D^{(r_j+r_k)/2}\prod_i\|x_i\|_{p_i}$. Apply this to
		\eqref{eq:multilinear-pair-supports}. As above, the layer calculation
		can first be performed for finite spectra and extended by uniform
		approximation. Since $\|x_0\|_{p_0}\le\|x_0\|_\infty$, taking the
		supremum over contractions $x_0$ proves the second bound.
		
		Finally, when $p_i=n$ for $1\le i\le n$, one has $p_0=\infty$.
		The sum of distances from any vertex is $\lfloor(n+1)^2/4\rfloor$:
		on a cycle of length $2m+1$ it is $2(1+\cdots+m)=m(m+1)$, and on
		a cycle of length $2m$ it is $2(1+\cdots+(m-1))+m=m^2$.
		It follows that $r_j=(\lfloor(n+1)^2/4\rfloor-d(0,j))/n$.
		For $j\ne k$ one has $d(0,j)+d(0,k)\ge1$, and therefore
		$r_j+r_k\le\kappa_n$.
	\end{proof}
	
	\begin{corollary}\label{cor:multilinear-product}
		Suppose that $\cN\subset\cM$ is irreducible with finite Jones index
		$\delta^2$, and let $T$ be a two-box operator. If $\sum_i1/p_i=1$, then
		$B_T(x_1,\ldots,x_n)=T(x_1\cdots x_n)$ satisfies
		\[
		\|B_T\|_{\infty^n\to1}=\|T\|_{\infty\to1},\qquad
		\|B_T\|_{\boldsymbol p\to1}=\|T\|_{\mathrm{cb}}
		\le\delta\,\KNM{\cN}{\cM}\,\|B_T\|_{\infty^n\to1}.
		\]
		For operator-positive bimodule map $T$ one may replace $\KNM{\cN}{\cM}$ by $\KNM{\cN}{\cM}_+$.
	\end{corollary}
	\begin{proof}
		Products of contractions are contractions, and each contraction occurs
		as a product using $1$ in the other slots. This proves the first
		equality. 
		For the mixed-exponent norm, H\"older's inequality shows $\|x_1\cdots x_n\|_1\le\prod_i\|x_i\|_{p_i}$. 
		Conversely, for $z=u|z|\in\cM$, put $x_i=|z|^{1/p_i}$ when $p_i<\infty$ and
		$x_i=1$ when $p_i=\infty$. 
		Left-multiply the first finite-exponent factor by $u$. 
		Then $x_1\cdots x_n=u|z|=z$ and
		\[
		\prod_{i=1}^n\|x_i\|_{p_i}
		=\prod_{i:p_i<\infty}\tau(|z|)^{1/p_i}
		=\|z\|_1 .
		\]
		The polar partial isometry preserves the norm of the factor to which
		it is applied, since the factor is supported on $u^*u$.
		Thus $\|B_T\|_{\boldsymbol p\to1}=\|T\|_{1\to1}$, and
		Proposition~\ref{cor:two-input-product} gives the remaining assertions.
	\end{proof}
	
	\begin{proposition}
		\label{prop:multilinear-lower}
		Let $\cN\subset\cN\rtimes(\bZ_2)^r$ be an outer crossed-product inclusion
		of $\mathrm{II}_1$ factors with index $\delta^2=2^r$.
		For each critical tuple $\displaystyle \sum_i1/p_i=1$ the best forward coefficient is at least
		\[
		\delta^{\,2-2\min_i 1/p_i}.
		\]
		In particular the $L^n$ coefficient is at least $\delta^{2-2/n}$.
		For each integer $h\ge1$ with $n\ge4h-1$ it is also at least
		$\delta^{2h(n+1-2h)/n}$.
	\end{proposition}
	
	\begin{proof}
		Let $g_1,\ldots,g_r$ be the standard generators of $G=(\bZ_2)^r$, write $u_g$ for the canonical group unitaries, and let $|g|$ be the number of nonzero coordinates of $g$. 
		Let $D=2^r$.
		Put
		\[
		e=\frac{1}{D}\sum_{g\in G}u_g,\qquad
		U=\prod_{j=1}^r\frac{1+\mathrm i u_{g_j}}{\sqrt2}
		=\frac{1}{\delta}\sum_{g\in G}\mathrm i^{|g|}u_g,\qquad
		T\Bigl(\sum_{g\in G}a_gu_g\Bigr)
		=\sum_{g\in G}\mathrm i^{|g|}a_gu_g .
		\]
		The group average $e$ is a projection with $\bE_{\cN}(e)=D^{-1}1$.
		The factors in the product for $U$ are unitaries because each $u_{g_j}$
		is a self-adjoint unitary. The Fourier summands $\cN u_g$ are orthogonal
		in $L^2$, so $T$ and $T^{-1}$ are $L^2$ isometries. Both are bounded
		normal $\cN$-bimodule maps fixing $1$, and $T(e)=\delta^{-1}U$.
		
		Choose a slot $j$ minimizing $1/p_j$ and set
		\[
		B(x_1,\ldots,x_n)=
		T\bigl(x_1\cdots x_{j-1}T^{-1}(x_j)x_{j+1}\cdots x_n\bigr).
		\]
		This map is separately normal and $\cN$-multimodular.
		For contractions $x_i$, its $L^2$ norm is bounded by one, since bounded
		left and right multipliers have $L^2$ operator norm at most their
		operator norm, and $T,T^{-1}$ are $L^2$ isometries. Also
		$B(1,\ldots,1)=1$, so $\|B\|_{\infty^n\to1}=1$. Now choose
		$x_j=U$ and $x_i=\delta^{2/p_i}e$ for $i\ne j$. Since
		$\tau(e)=\delta^{-2}$ these have the prescribed unit $L^{p_i}$ norms.
		Using $T^{-1}(U)=\delta e$ and $e^2=e$ gives
		\[
		B(x_1,\ldots,x_n)
		=\delta^{1+2\sum_{i\ne j}1/p_i}T(e)=\delta^{2-2/p_j}U.
		\]
		Taking the trace norm proves the first lower bound.
		
		For the second bound, choose $q_j\in\cN^G$ with
		$\tau(q_j)=\delta^{-2j}$ and $q_0=1$. The fixed-point algebra of the
		outer action is a $\mathrm{II}_1$ factor, so these projections exist.
		They commute with the group unitaries. Set
		\[
		P_j=q_{j-1}e,\qquad V_j=q_{j-1}U .
		\]
		Then $P_j$ is a projection of trace $\delta^{-2j}$, $V_j$ is unitary
		in the $q_{j-1}$ corner, and $T(P_j)=\delta^{-1}V_j$ by bimodularity.
		
		Build a multilinear map starting with $T^{-1}(x)$ and at each step
		replace the current output $y$ by $T^{\pm1}(a y b)$, introducing new
		inputs $a,b$. The module identities follow inductively because $T^{\pm1}$
		are bimodule maps and multiplication moves $\cN$-elements between
		adjacent factors. Each map in the construction sends tuples of
		contractions into the $L^2$ unit ball and takes the tuple of identities
		to $1$. Its $\infty^n\to1$ norm is consequently one.
		
		The input $U$ starts the construction at $\delta P_1$.
		For $1\le j<h$, take a unitary $v\in\cM$ with $vP_jv^*=q_j$
		and choose $a=vP_j$, $b=P_jv^*U$. Both are partial isometries of
		support trace $\delta^{-2j}$, and
		\[
		aP_jb=q_jU=V_{j+1},\qquad
		T^{-1}(a\,\delta^jP_j\,b)=\delta^{j+1}P_{j+1}.
		\]
		At $\delta^hP_h$, take $a=b=P_h$ and apply $T$ to get
		$\delta^{h-1}V_h$. For the descending steps, $h\ge j>1$, choose a
		unitary $v\in\cM$ with $v q_{j-1}v^*=P_{j-1}$ and use
		$a=v q_{j-1}$, $b=V_j^*v^*P_{j-1}$. Both are partial isometries of
		support trace $\delta^{-2(j-1)}$, and
		\[
		aV_jb=P_{j-1},\qquad
		T(a\,\delta^{j-1}V_j\,b)=\delta^{j-2}V_{j-1}.
		\]
		The final output is $V_1=U$. There are $4h-1$ inputs in all, including the
		initial $U$. For each $1\le j<h$ there are four inputs of support trace
		$\delta^{-2j}$, and there are two inputs of support trace $\delta^{-2h}$.
		Additional inputs can be inserted at the turning step and evaluated at
		$P_h$, since multiplying by further copies of $P_h$ fixes that step's
		value. Each partial isometry of support trace $\delta^{-2j}$ has
		$L^n$ norm $\delta^{-2j/n}$. After rescaling the inputs to unit $L^n$
		norm, the output has trace norm
		\[
		\delta^{(2/n)(4\sum_{j=1}^{h-1}j+2h+h(n-4h+1))}
		=\delta^{2h(n+1-2h)/n}.
		\]
		This completes the proof of the proposition.
	\end{proof}
	
	For $n=3$, Proposition~\ref{prop:multilinear-lower} and
	Theorem~\ref{thm:multilinear-fixed-trace} give lower and upper bounds of
	$D^{2/3}$ and $O(D^{7/6})$, respectively. In terms of powers of the
	index $D$, choosing $h=\lfloor(n+1)/4\rfloor$ gives a lower exponent
	asymptotic to $n/8$. The exponent $\kappa_n/2$ in the second upper bound
	is asymptotic to $n/4$. In particular, a coefficient bounded by a fixed
	power of $D$ uniformly in the arity is impossible.
	
	\section{An index-free Grothendieck inequality}
	
	Let $\cN\subset\cM$ be an irreducible inclusion of finite factors, and
	let $\tau$ be the normalized trace on $\cM$. For a rectangular array
	$X=[x_{ab}]$ over $\cM$, let $\|X\|_{\op}$ be its block operator norm and put
	\begin{equation}\label{eq:rho}
		s(X)^2=\sum_{a,b}\tau(x_{ab}^*x_{ab}),\qquad
		\rho(X)=\bigl(\|X\|_{\op}^2+s(X)^2\bigr)^{1/2}.
	\end{equation}
	The sum uses the unnormalized matrix trace. Both quantities are unchanged
	by padding with zero rows or columns. If $X$ is a single row then
	$s(X)^2=\tau(XX^*)\le\|X\|_{\op}^2$, and for a single column use
	$\tau(X^*X)$ instead.
	
	\begin{theorem} \label{thm:path}
		Let $n\ge2$ and let $\Phi:\cM^{\times n}\to\cM$ be a bounded, separately
		normal multilinear map. Assume that for $1\le j<n$ and $a\in\cN$,
		\[
		\Phi(\ldots,x_j a,x_{j+1},\ldots)=
		\Phi(\ldots,x_j,a x_{j+1},\ldots).
		\]
		Take $d_0=d_n=1$ and positive integers $d_1,\ldots,d_{n-1}$.
		For $X_j=[x^j_{ab}]\in M_{d_{j-1},d_j}(\cM)$ write $i_0=i_n=1$.
		Then
		\begin{align}
			\Big\|\sum_{i_1,\ldots,i_{n-1}}
			\Phi(x^1_{i_0i_1},\ldots,x^n_{i_{n-1}i_n})\Big\|_1
			&\le 2^{(n-2)/2}\|\Phi\|_{\infty^n\to1}\prod_{j=1}^n\rho(X_j)
			\label{eq:path}\\
			&\le 2^{n/2}\|\Phi\|_{\infty^n\to1}
			\|X_1\|_{\op}\|X_n\|_{\op}\prod_{j=2}^{n-1}\rho(X_j).
			\label{eq:pathrc}
		\end{align}
	\end{theorem}
	
	For example, at $n=3$ the second bound reads
	\[
	\Big\|\sum_{a,b}\Phi(x_a,y_{ab},z_b)\Big\|_1
	\le 2\sqrt2\,\|\Phi\|_{\infty^3\to1}
	\Big\|\sum_a x_ax_a^*\Big\|_\infty^{1/2}
	\left(\|[y_{ab}]\|_{\op}^2+\sum_{a,b}\|y_{ab}\|_2^2\right)^{1/2}
	\Big\|\sum_b z_b^*z_b\Big\|_\infty^{1/2}.
	\]
	At $n=2$ the product of interior factors is empty, and
	\eqref{eq:pathrc} recovers the balanced row--column inequality for
	bilinear maps. 
	Diagonal interior matrices with diagonal entries $x_{j,i}$ give a single-index sum $\displaystyle \left\|\sum_i\Phi(x_{1,i},\ldots,x_{n,i})\right\|_1$, with each interior factor equal to $\displaystyle \big(\max_i\|x_{j,i}\|_\infty^2+\sum_i\|x_{j,i}\|_2^2\big)^{1/2}$.
	
	\begin{theorem}\label{thm:cyc}
		Let $k\ge2$, and let $F:\cM^{\times k}\to\bC$ be bounded, separately normal and cyclically $\cN$-balanced. 
		Thus $a\in\cN$ can be moved from the right of one input to the left of the following input, including from the last input to the first. 
		Take $d_k=d_0$ and $X_j=[x^j_{ab}]\in
		M_{d_{j-1},d_j}(\cM)$ for $1\le j\le k$. 
		Then
		\begin{equation}\label{eq:cyc}
			\left|\sum_{i_0,\ldots,i_{k-1}}
			F(x^1_{i_0i_1},x^2_{i_1i_2},\ldots,x^k_{i_{k-1}i_0})\right|
			\le 2^{k/2}\|F\|\prod_{j=1}^k\rho(X_j).
		\end{equation}
		For the trace form $F(x_0,x_1,\ldots,x_n)=\tau(x_0\Phi(x_1,\ldots,x_n))$ of an $\cN$-multimodular map, the coefficient multiplying the input factors is $2^{(n+1)/2}\|\Phi\|_{\infty^n\to1}$.
	\end{theorem}
	
	Write $\mathcal A=M_m(\cM)$, $\mathcal P=M_m(\cN)$ and $p=e_{11}\otimes1$. 
	We identify
	each scalar matrix unit with its amplification by $1$ in $\mathcal A$. 
	With normalized trace $\tau_m=m^{-1}(\operatorname{Tr}_m\otimes\tau)$, one has $s(X)^2=m\tau_m(X^*X)$. 
	The row and column observation gives
	\begin{equation}\label{eq:pinch}
		\rho(Xp)\le\sqrt2\|X\|_{\op},\qquad
		\rho(pX)\le\sqrt2\|X\|_{\op}.
	\end{equation}
	
	The following lemma is obtained from \cite[Lemma 3.1]{CPSS} and here the common convex averages in \(L^2\) replace the hyperfinite averaging step.
	\begin{lemma}\label{lem:corner}
		Let $\psi:\mathcal A\times \mathcal A\to\bC$ be a separately normal bilinear form with
		$\psi(Xa,Y)=\psi(X,aY)$ for $a\in \mathcal P$.
		If
		$|\psi(Xp,pY)|\le C\|X\|_{\op}\|Y\|_{\op}$, then
		$|\psi(X,Y)|\le C\rho(X)\rho(Y)$.
	\end{lemma}
	\begin{proof}
		The corner form $\theta(X,Y)=\psi(Xp,pY)$ has norm at most $C$.
		Haagerup's four-state estimate for separately normal forms
		\cite{Haagerup,CPSS} gives normal states $f,h,g,\ell$ on $\mathcal  A$ with
		\[
		|\theta(X,Y)|\le C
		\bigl[f(XX^*)+h(X^*X)\bigr]^{1/2}
		\bigl[g(YY^*)+\ell(Y^*Y)\bigr]^{1/2}.
		\]
		Because $\displaystyle \sum_i e_{i1}e_{1i}=1$ and $e_{i1}p=e_{i1}$,
		$p e_{1i}=e_{1i}$, balance gives
		\[
		\psi(X,Y)=\sum_i\theta(Xe_{i1},e_{1i}Y).
		\]
		Summing the estimates by the classical Cauchy--Schwarz inequality, and multiplying out the matrix units, gives the bound
		\[
		|\psi(X,Y)|\le C
		\bigl[f(XX^*)+h(e_{11}\otimes A_X)\bigr]^{1/2}
		\bigl[g(e_{11}\otimes A_Y)+\ell(Y^*Y)\bigr]^{1/2}.
		\]
		Here $\displaystyle A_X=\sum_i(X^*X)_{ii}$ and $\displaystyle A_Y=\sum_i(YY^*)_{ii}$.
		For $u\in\U(\cN)$ substitute $X(1_m\otimes u)$ and
		$(1_m\otimes u^*)Y$. 
		This preserves $\psi(X,Y)$ and leaves $XX^*$ and
		$Y^*Y$ unchanged. The two diagonal sums become $u^*A_Xu$ and $u^*A_Yu$.
		
		Use common convex averages of these two conjugations. 
		In
		$L^2(\cM)\oplus L^2(\cM)$ the closed convex hull of the joint orbit
		has a unique element of minimal norm, which is fixed by all the
		conjugations. 
		The fixed vectors are scalar in each component by
		irreducibility: the spectral projections of a self-adjoint fixed vector are in $\cN'\cap\cM$, and taking real and imaginary parts gives the
		general case.
		Since averaging preserves both traces, the minimal
		point is $(\tau(A_X)1,\tau(A_Y)1)$.
		The averages are uniformly bounded
		in operator norm, so $L^2$ convergence implies convergence against
		normal functionals. 
		For example, one can approximate a normal
		functional's $L^1$ density by bounded densities, and use the
		Cauchy--Schwarz inequality on the bounded densities. 
		Apply the classical Cauchy--Schwarz inequality to the same averaging weights in the last bound on $\psi$.
		Passing to the limit yields
		\[
		|\psi(X,Y)|\le C
		\bigl[f(XX^*)+h(p)\,\tau(A_X)\bigr]^{1/2}
		\bigl[g(p)\,\tau(A_Y)+\ell(Y^*Y)\bigr]^{1/2}.
		\]
		Now $\tau(A_X)=s(X)^2$, $\tau(A_Y)=s(Y)^2$, and $h(p),g(p)\le1$.
		The remaining state evaluations are bounded by $\|X\|_{\op}^2$ and
		$\|Y\|_{\op}^2$, respectively. 
		This proves the estimate.
	\end{proof}
	
	\begin{proof}[Proof of Theorems~\ref{thm:path} and~\ref{thm:cyc}]
		Pad all matrices to a common size $m$ using zero rows and columns. For the path, fix a contraction $w\in\cM$, set
		$H(x_1,\ldots,x_n)=\tau(w\Phi(x_1,\ldots,x_n))$, and keep the endpoint
		indices fixed at $1$. 
		For the cycle use $H=F$.
		Let $S$ be a set of interfaces whose indices are allowed to range over
		$\{1,\ldots,m\}$, and fix every other index at $1$. 
		Denote the resulting sum by $H_S$.
		For each matrix position $j$ let $a_j$ be the number of its
		interfaces in $S$. 
		We prove the bound
		\begin{equation}\label{eq:interface-induction}
			|H_S|\le \|H\|\,2^{|\{j:a_j=2\}|/2}
			\prod_{j:a_j=0}\|X_j\|_{\op}\prod_{j:a_j>0}\rho(X_j).
		\end{equation}
		It holds when $S=\varnothing$, since only the $(1,1)$ entries contribute and their norms are bounded by the block operator norms.
		
		To enlarge $S$, let $e$ be the new interface and write $X,Y$ for the
		matrix on its left and right in the order of traversal.
		With the
		other matrices fixed, denote $H_{S\cup\{e\}}$ by $\psi(X,Y)$.
		This is separately normal and balanced over $P$. 
		Indeed, for $a=[a_{uv}]\in
		\mathcal P$, the consecutive entries in $\psi(Xa,Y)$ have the form
		$\sum_t x_{bt}a_{ti}$ and $y_{ic}$, where $i$ is the index at $e$.
		In the sum over $t,i$, move $a_{ti}$ to the left of $y_{ic}$.
		The entries now take the form $x_{bt}$ and $\sum_i a_{ti}y_{ic}$,
		which gives $\psi(X,aY)$. 
		The same expansion at the closing interface uses the closing balance identity for $F$.
		
		The form $\psi(Xp,pY)$ fixes the new index at $1$.
		The induction bound therefore applies, with $Xp,pY$ at the incident positions.
		Each of those positions has either zero or one interface in the
		old $S$. 
		In the first case the corresponding factor is at most
		$\|X\|_{\op}$ or $\|Y\|_{\op}$.
		In the second case \eqref{eq:pinch}
		bounds the factor by $\sqrt2$ times that operator norm. 
		Thus
		Lemma~\ref{lem:corner} gives \eqref{eq:interface-induction} for the enlarged
		set with $\rho(X),\rho(Y)$ in place of the operator norms.
		The $\sqrt2$ factors occur precisely at the positions whose count $a_j$ increases from one to two.
		
		When all interfaces are allowed to range, $n-2$ positions on the path and all $k$ positions on the cycle have $a_j=2$. 
		This proves
		\eqref{eq:cyc}. 
		In the path case $\|H\|\le\|\Phi\|_{\infty^n\to1}$, and
		taking the supremum over $w$ gives \eqref{eq:path}. 
		Using the row and column norm bounds for $X_1,X_n$ gives \eqref{eq:pathrc}.
	\end{proof}
	
	% \begin{remark}
		% The estimates \eqref{eq:path}--\eqref{eq:cyc} use the mixed matrix norm of
		% \eqref{eq:CPSS} and a trace pairing with the output. Their constants
		% depend on the arity and apply at every Jones index, including infinite
		% index. Cyclic balance suffices to move an $\cN$ coefficient across
		% the closing interface.
		% \end{remark}
	\section*{Declaration of generative AI use}
	
	During the preparation of this work, the authors used generative AI
	tools to improve the language and readability of the manuscript and
	to assist with developing and checking mathematical arguments.
	The authors reviewed and revised the AI-assisted content and verified
	the mathematical statements and proofs.
	The authors take full responsibility for the accuracy and integrity
	of the manuscript.
	
	\bibliographystyle{plain}
	\bibliography{young}
	
\end{document}